\documentclass[10pt]{amsart}
\usepackage{xypic,amscd,amssymb,latexsym}
\usepackage{amsmath}	
\usepackage{graphicx}
\usepackage{pdfsync}
\usepackage{cancel}

\usepackage{pstricks}
\usepackage{pst-poly}

\usepackage{mathrsfs}

\usepackage{url}
\usepackage{pifont}

\DeclareMathAlphabet{\mathpzc}{OT1}{pzc}{m}{it}

\usepackage{setspace}

\usepackage{changepage}

\usepackage{tikz}
\renewcommand*\circled[1]{\tikz[baseline=(char.base)]{
            \node[shape=circle,draw,inner sep=2pt] (char) {#1};}}

\begin{document}
\pagenumbering{arabic}
\pagestyle{plain}

 \newtheorem{theorem}{Theorem}[section]
\newtheorem{proposition}[theorem]{Proposition}
\newtheorem{lemma}[theorem]{Lemma}
\newtheorem{corollary}[theorem]{Corollary}
\newtheorem{remark}[theorem]{Remark}
\newtheorem{definition}[theorem]{Definition}
\newtheorem{question}[theorem]{Question}
\newtheorem{claim}[theorem]{Claim}
\newtheorem{conjecture}[theorem]{Conjecture}
\newtheorem{defprop}[theorem]{Definition and Proposition}
\newtheorem{example}[theorem]{Example}
\newtheorem{deflem}[theorem]{Definition and Lemma}

\newcommand{\II}[0]{{\rm \it II}}

\def\qed{{\quad \vrule height 8pt width 8pt depth 0pt}}

\newcommand{\cplx}[0]{\mathbb{C}}

\newcommand{\fr}[1]{\mathfrak{#1}}

\newcommand{\vs}[0]{\vspace{2mm}}

\newcommand{\til}[1]{\widetilde{#1}}

\newcommand{\mcal}[1]{\mathcal{#1}}

\newcommand{\ul}[1]{\underline{#1}}

\newcommand{\ol}[1]{\overline{#1}}

\newcommand{\wh}[1]{\widehat{#1}}

\address{Department of Mathematics, Ewha Womans University, 52 Ewhayeodae-gil, Seodaemun-gu, Seoul 03760, Republic of Korea}

\author{Hyun Kyu Kim}
\thanks{keywords: quantum torus, representation theory, quantum Teichm\"uller theory, 6j symbols, Clebsch-Gordan coefficients}

\email[H.~Kim]{hyunkyukim@ewha.ac.kr, hyunkyu87@gmail.com}

\numberwithin{equation}{section}

\title{Finite dimensional quantum Teichm\"uller space from the quantum torus at root of unity}

\begin{abstract}
Representation theory of the quantum torus Hopf algebra, when the  parameter $q$ is a root of unity, is studied. We investigate a decomposition map of the tensor product of two irreducibles into the direct sum of irreducibles, realized as a `multiplicity module' tensored with an irreducible representation.  The isomorphism between the two possible decompositions of the triple tensor product yields a map ${\bf T}$ between the multiplicity modules, called the 6j-symbols. We study the left and right dual representations, and correspondingly, the left and right representations on the ${\rm Hom}$ spaces of linear maps between representations. Using the isomorphisms of irreducibles to left and right duals, we construct a map ${\bf A}$ on a multiplicity module, encoding the permutation of the roles of the irreducible representations in the identification of the multiplicity module as the space of intertwiners between representations. We show that ${\bf T}$ and ${\bf A}$ satisfy certain consistency relations, forming a Kashaev-type quantization of the Teichm\"uller spaces of bordered Riemann surfaces. All constructions and proofs in the present work uses only plain representation theoretic language with the help of the notions of the left and the right dual and Hom representations, and therefore can be applied easily to other Hopf algebras for future works.
\end{abstract}

\maketitle

\vspace{-7mm}

\tableofcontents

\section{Background and introduction}

Quantum Teichm\"uller theory emerged as an approach to 3 or 2+1 dimensional quantum gravity; it was established in 1990's independently by Kashaev \cite{K98} and by Chekhov-Fock \cite{Fock} \cite{CF} (the latter methodology generalized to Fock-Goncharov quantization of cluster varieties later). The Teichm\"uller space $\mathcal{T}(S)$ of a punctured surface $S$ is the set of all isotopy classes of complex structures on $S$, with various different prescribed behaviors at the punctures, leading to different versions. A version of $\mathcal{T}(S)$ is a smooth manifold with a symplectic or Poisson structure, called the {\em Weil-Petersson} structure, which is invariant under the natural action on $\mathcal{T}(S)$ of the {\em mapping class group} $\mathrm{MCG}(S)$ of $S$, defined as the group of isotopy classes of orientation-preserving self-diffeomorphisms of $S$. One thus looks for an `equivariant deformation quantization' of $\mathcal{T}(S)$, which is to replace smooth real functions on $\mathcal{T}(S)$ by self-adjoint operators on a Hilbert space depending on a real quantum parameter $\hbar$ so that the commutators of these operators recover the classical Poisson bracket in the $1$st order terms in $\hbar$, together with a consistent assignment of unitary operators to elements of $\mathrm{MCG}(S)$ that intertwine these self-adjoint operators, where in the classical limit $\hbar \to 0$ this intertwining action of $\mathrm{MCG}(S)$ must recover the classical $\mathrm{MCG}(S)$ action. As a result, one would obtain a real-parameter family of projective unitary representations of $\mathrm{MCG}(S)$ on a Hilbert space.

\vs

The quantization problem just described is too difficult to solve, and what has been done is as follows. We first choose an extra combinatorial-topological data on $S$, namely an {\em ideal triangulation} of $S$, which is a triangulation of $S$, defined up to homotopy, whose vertices are the punctures. With the help of this extra data, one constructs a coordinate system of $\mathcal{T}(S)$, hence coordinate functions, which we replace by certain self-adjoint operators on a Hilbert space, containing the information of the classical Poisson bracket appropriately. A different ideal triangulation leads to different classical coordinate functions, related to the previous ones by certain coordinate change formulas, and we construct different quantum coordinate operators on a different Hilbert space. To each change of ideal triangulations, we construct a unitary map between the corresponding Hilbert spaces that intertwine the quantum coordinate operators, so that this intertwining formula recovers the classical coordinate change map in the limit $\hbar \to 0$, and that the composition of changes of ideal triangulations goes to the composition of unitary intertwining maps, up to multiplicative constants. This also leads to a projective representation of $\mathrm{MCG}(S)$, because each element $\mathrm{MCG}(S)$ can be realized as a change of ideal triangulations, and the above construction of quantum coordinate operators has a necessary built-in invariance under $\mathrm{MCG}(S)$.

\vs

\begin{figure}[htbp!]
$\begin{array}{llll}
\begin{pspicture}[showgrid=false,linewidth=0.5pt,unit=6mm](-1,-1.2)(2.0,1.3)
%
\psarc[arcsep=0.5pt](1.091,-3.055){2.182}{61.7}{158.0}
\psarc[arcsep=0.5pt](-7.2,0.4){6.8}{-22}{16.3}
\psarc[arcsep=0.5pt](4.8,3.9){5.7}{-163.7}{-118.3}
%
%
\rput{127.0}(1.2,1.0){\fontsize{19}{19} $\cdots$}
\rput{86.0}(-1.4,0.2){\fontsize{19}{19} $\cdots$}
\rput{20.0}(0.58,-1.65){\fontsize{17}{17} $\cdots$}
%
\rput(-0.5,-1.3){\fontsize{11}{11} $\bullet$}
%
\rput(0.1,-0.3){\fontsize{11}{11} $t$}
%
\rput[l](3.2,0){\pcline[linewidth=0.7pt, arrowsize=2pt 4]{->}(0,0)(1;0)\Aput{\,${\bf A}_t$}}
\end{pspicture}
%
%
%
%
& \begin{pspicture}[showgrid=false,linewidth=0.5pt,unit=6mm](-2.0,-1.2)(2.4,1.3)
%
\psarc[arcsep=0.5pt](1.091,-3.055){2.182}{61.7}{158.0}
\psarc[arcsep=0.5pt](-7.2,0.4){6.8}{-22}{16.3}
\psarc[arcsep=0.5pt](4.8,3.9){5.7}{-163.7}{-118.3}
%
\rput(1.1,-0.7){\fontsize{11}{11} $\bullet$}
%
\rput(0.1,-0.3){\fontsize{11}{11} $t$}
%
\rput{127.0}(1.2,1.0){\fontsize{19}{19} $\cdots$}
\rput{86.0}(-1.4,0.2){\fontsize{19}{19} $\cdots$}
\rput{20.0}(0.58,-1.65){\fontsize{17}{17} $\cdots$}
\end{pspicture}
& \begin{pspicture}[showgrid=false,linewidth=0.5pt,unit=6mm](-2.0,-1.2)(2.0,1.3)
%
\psarc[arcsep=0.5pt](0.343,-4.457){3.771}{61.7}{126.4}
\psarc[arcsep=0.5pt](-3.36,0.48){2.4}{-53.6}{37.4}
\psarc[arcsep=0.5pt](-0.185,2.862){1.569}{-142.6}{-29}
\psarc[arcsep=0.5pt](3.220,0.937){2.341}{151}{-118.3}
%
\psarc[arcsep=0.5pt](7.2,8.4){10.8}{-142.6}{-118.3}
%
\rput(-0.95,1.00){\fontsize{10}{10} $\bullet$}
\rput(0.72,1.42){\fontsize{10}{10} $\bullet$}
%
\rput(-0.5,-0.1){\fontsize{10}{10} $s$}
\rput(0.3,0.8){\fontsize{10}{10} $t$}
%
\rput{104.0}(1.7,0.5){\fontsize{15}{15} $\cdots$}
\rput{5.0}(-0.10,1.94){\fontsize{15}{15} $\cdots$}
\rput{85.0}(-1.7,0.3){\fontsize{15}{15} $\cdots$}
\rput{4.0}(0.08,-1.50){\fontsize{17}{17} $\cdots$}
%
\rput[l](3.2,0){\pcline[linewidth=0.7pt, arrowsize=2pt 4]{->}(0,0)(1;0)\Aput{${\bf T}_{st}$}}
\end{pspicture}
%
%
%
%
& \begin{pspicture}[showgrid=false,linewidth=0.5pt,unit=6mm](-2.0,-1.2)(2.4,1.3)
%
\psarc[arcsep=0.5pt](0.343,-4.457){3.771}{61.7}{126.4}
\psarc[arcsep=0.5pt](-3.36,0.48){2.4}{-53.6}{37.4}
\psarc[arcsep=0.5pt](-0.185,2.862){1.569}{-142.6}{-29}
\psarc[arcsep=0.5pt](3.220,0.937){2.341}{151}{-118.3}
%
%
\psarc[arcsep=0.5pt](-8.8,7.733){11.467}{-53.0}{-29.7}
%
\rput(-1.00,1.30){\fontsize{10}{10} $\bullet$}
\rput(0.67,1.10){\fontsize{10}{10} $\bullet$}
%
\rput(-0.50,0.8){\fontsize{10}{10} $s$}
\rput(0.3,-0.1){\fontsize{10}{10} $t$}
%
\rput{104.0}(1.7,0.5){\fontsize{15}{15} $\cdots$}
\rput{5.0}(-0.10,1.94){\fontsize{15}{15} $\cdots$}
\rput{85.0}(-1.7,0.3){\fontsize{15}{15} $\cdots$}
\rput{4.0}(0.08,-1.50){\fontsize{17}{17} $\cdots$}
\end{pspicture}
\end{array} $
\caption{The actions of ${\bf A}_t$ and ${\bf T}_{st}$ on dotted ideal triangulations}
\label{fig:action_on_dotted_ideal_triangulations}
\end{figure}
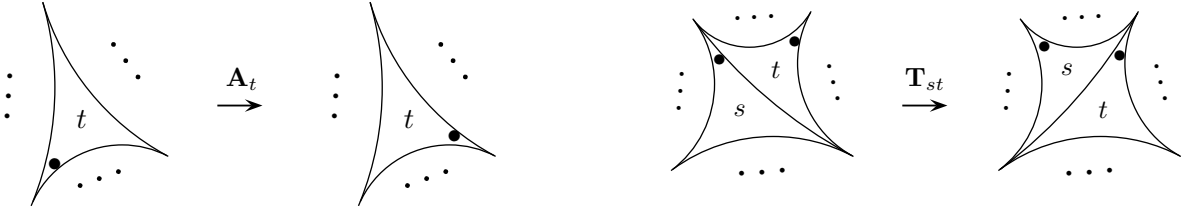

Breaking down the problem even more, one observes that any change of ideal triangulations is generated by the elementary ones called the {\em flip} along an edge, which changes an ideal triangulation only on one edge, replacing this edge by the other diagonal of the unique quadrilateral in which it was contained as a diagonal. The flips satisfy certain algebraic relations; twice-flip at the same edge is the identity, alternating sequence of five flips at two adjacent edges is the identity (called the {\em pentagon relation}), and two flips at edges not shared by a triangle commute with each other. It is a theorem that any relation is a consequence of those. Chekhov-Fock(-Goncharov) constructed suitable intertwining operator for each flip, and verified that the algebraic relations are satisfied by these operators up to constants. For Kashaev's quantization, we use {\em dotted ideal triangulations}, which are ideal triangulations together with the choice of a distinguished corner at each triangle and with the choice of labeling of triangles. Any change of dotted ideal triangulations are generated by elementary moves of three types: 1) ${\bf A}_t$, for a triangle $t$, moves the dot of the triangle $t$ counterclockwise, 2) ${\bf T}_{st}$, for two adjacent triangles $s$ and $t$, is the flip that is confined to the dot configuration of triangles $s$ and $t$ as in Fig.\ref{fig:action_on_dotted_ideal_triangulations} \footnote{I note that the basic source code for Fig.\ref{fig:action_on_dotted_ideal_triangulations} is taken from \cite{Kim16}.}, and 3) ${\bf P}_\gamma$, for a permutation $\gamma$ of triangle labels, permutes the triangle labels. 
Any algebraic relation among them is a consequence of the following ones:
\begin{align}
\label{eq:Kashaev_relations}
{\bf A}_t^3 = \mathrm{id}, \quad
{\bf T}_{st} {\bf T}_{rs} = {\bf T}_{rs} {\bf T}_{rt} {\bf T}_{st}, \quad
{\bf A}_s {\bf T}_{st} {\bf A}_t = {\bf A}_t {\bf T}_{ts} {\bf A}_s, \quad
{\bf T}_{st} {\bf A}_s {\bf T}_{ts} = {\bf A}_s {\bf A}_t {\bf P}_{(st)},
\end{align}
and `trivial' relations: any generators whose sets of subscripts are disjoint commute with one another, conjugation by ${\bf P}_\gamma$ on a generator yields the same generator with subscripts applied by $\gamma$, and the permutation group relations for ${\bf P}_\gamma$. Kashaev constructs unitary intertwining operators for each elementary change of dotted ideal triangulations, and verified that the algebraic relations are satisfied up to constants. See \cite{Kim16} for an interesting discrepancy between the $\mathrm{MCG}(S)$ representations from the Chekhov-Fock-Goncharov quantization and the Kashaev quantization.

\vs

The Chekhov-Fock and Kashaev quantizations are built on infinite dimensional separable Hilbert spaces, looking naturally like $L^2(\mathbb{R}^n)$. Meanwhile, the unitary intertwining operators ${\bf A}_t, {\bf T}_{st}, {\bf P}_\gamma$ of the Kashaev quantization are recovered purely from the representation theory of a Hopf algebra, in my joint work with Igor Frenkel \cite{FK12}. We considered one of the most basic quantum groups, called the {\em quantum plane} algebra $B_q$, defined for a complex parameter $q$ of modulus $1$ which is not a root of unity, which one can think as being related to the irrational positive real parameter $\hbar = b^2$ by $q = e^{\pi \sqrt{-1}\, \hbar}$. It is defined by
$$
\langle X^{\pm 1}, Y | XY = q^2 YX\rangle
$$
as an algebra over $\mathbb{C}$, equipped with the Hopf $*$-algebra structure given by
\begin{align*}
{\renewcommand{\arraystretch}{1.0}  \begin{array}{rll}
\mbox{coproduct} & : \quad  \Delta X = X\otimes X, & \quad \Delta Y = Y \otimes X + 1 \otimes Y, \\
\mbox{counit} & :  \quad \epsilon(X) = 1, & \quad \epsilon(Y) = 0, \\
\mbox{antipode} & : \quad S(X) = X^{-1}, & \quad S(Y) = -YX^{-1}, \\
\mbox{$*$-structure} & : \quad X^* = X, & \quad Y^* = Y.
\end{array} }
\end{align*}
A reader may recognize this as a Borel subalgebra of the famous quantum group $\mathcal{U}_q(\mathfrak{sl}_2)$, or more precisely, $\mathcal{U}_q(\mathfrak{sl}(2,\mathbb{R}))$. When considering representations, the $*$-structure stipulates that the operators for $X$ and $Y$ be self-adjoint. A certain class of nicely behaved representations called {\em integrable} representations of $B_q$ is studied in the literature \cite{Schmudgen}. An important observation is that an irreducible integrable representation is unique up to unitary isomorphisms, denoted by $\mathscr{H} \equiv L^2(\mathbb{R}, dx)$, where the action of the generators is given by $(X.f)(x) = f(x-\sqrt{-1}\,b)$ and $(Y.f)(x) = e^{2\pi b x} f(x)$ on a dense subspace (recall $\hbar = b^2$); this observation is essentially from the Stone-von Neumann theorem \cite{vonNeumann}. The tensor product $\mathscr{H} \otimes \mathscr{H}$ of Hilbert spaces becomes a representation of $B_q$ via the coproduct, as usual in the theory of Hopf algebras; it is integrable, so we expect it to decompose into direct sum of irreducibles, or in this case, direct integral of the unique irreducible $\mathscr{H}$. Such a decomposition is realized as a unitary $B_q$-intertwining isomorphism
\begin{align}
\label{eq:the_F-map}
F : \mathscr{H} \otimes \mathscr{H} \to M \otimes \mathscr{H},
\end{align}
where $M \equiv L^2(\mathbb{R})$, the `multiplicity space', is a trivial representation of the Hopf algebra $B_q$, meaning that the action is by counit; as a consequence, on the RHS of \eqref{eq:the_F-map} the algebra $B_q$ acts only on the second tensor factor $\mathscr{H}$. The triple tensor product $\mathscr{H} \otimes \mathscr{H} \otimes \mathscr{H}$ can also be decomposed into a direct integral of $\mathscr{H}$, or more precisely can be identified with $M \otimes M \otimes \mathscr{H}$. This can be done in two ways using the $F$ map, depending on whether we choose the parenthesizing $(\mathscr{H} \otimes \mathscr{H}) \otimes \mathscr{H}$ or $\mathscr{H} \otimes (\mathscr{H} \otimes \mathscr{H})$; the first means to decompose the first two factors to get $M \otimes \mathscr{H} \otimes \mathscr{H}$, and then decompose into $M \otimes M \otimes \mathscr{H}$, while the second means first to go to $\mathscr{H} \otimes M \otimes \mathscr{H}$ then to $M \otimes M \otimes \mathscr{H}$. Composition of these two decomposition maps yields a map $M \otimes M \otimes \mathscr{H} \to M \otimes M \otimes \mathscr{H}$, which is of the form ${\bf T} \otimes \mathrm{id}$ for some operator ${\bf T} : M \otimes M \to M \otimes M$, encoding the parenthesis-change $(\mathscr{H} \otimes \mathscr{H}) \otimes \mathscr{H} \leadsto \mathscr{H} \otimes (\mathscr{H} \otimes \mathscr{H})$. This operator corresponds to what is called the `6j-symbol' in the mathematical physics literature. Meanwhile, the above $M$ can be viewed as the space of intertwiners $\mathrm{Hom}_{B_q}(\mathscr{H}, \mathscr{H} \otimes \mathscr{H})$; by studying the dual representations of the unique irreducible $\mathscr{H}$ with the help of the antipode of $B_q$, we were able to construct a natural map ${\bf A} : M \to M$ which amounts to cyclically permuting the roles of the three $\mathscr{H}$'s in the intertwiner space $\mathrm{Hom}_{B_q}(\mathscr{H}, \mathscr{H} \otimes \mathscr{H})$. We then showed that these operators ${\bf T}$ (on $M\otimes M$) and ${\bf A}$ (on $M$) coincide with the operators of the same names from Kashaev's quantum Teichm\"uller theory, up to unitary transformation and a slight modification; the ${\bf P}_\gamma$ operator is just the permutation of the $M$ factors. One way of summarizing our result is to say that, in special cases, we recovered the quantum Teichm\"uller (Hilbert) space as the space of intertwiners of the quantum plane algebra $B_q$, and Kashaev's unitary operators ${\bf A}, {\bf T}, {\bf P}$ as certain maps on the space of intertwiners of $B_q$ naturally constructed using the representation theory of $B_q$.

\vs

It is striking that although the representation theory of a Hopf algebra is purely algebraic, its rigid tensor category structure yields the main result of Kashaev's quantum Teichm\"uller theory, which falls into the realm of quantum geometry. So, our work \cite{FK12} just described gives an indication that the structure of the Kashaev quantization, although it may look not so natural at the first sight (e.g. because of the usage of dotted triangulations), might also be something quite natural in mathematics. To add one more such indication, here I make a small side remark to be revisited later; the Kashaev's original quantization, rather than Chekhov-Fock's, is what is used crucially in Teschner's solution \cite{Teschner} of the `$\mathfrak{sl}_2$ modular functor conjecture'.

\vs

On the other hand, in the literature, a finite dimensional version of the Chekhov-Fock quantization has been studied \cite{BBL}. To each ideal triangulation is assigned a finite dimensional complex vector space, and for each flip an isomorphism between these vector spaces is constructed. A work in the flavor of \cite{FK12} for this finite dimensional version of Chekhov-Fock quantization is done by Bai \cite{Bai}. The same Hopf algebra as $B_q$, without the $*$-structure, is considered, which is named the {\em Weyl algebra} $\mathcal{W}$, or simply the {\em quantum torus} as is called in the present paper and by many authors. In case $q$ is a root of unity, there exist finite dimensional representations, although the irreducible representation is not unique anymore. The corresponding 6j symbol is studied by Kashaev \cite{K94} \cite{K99}, and in \cite{Bai} Bai showed that this 6j symbol coincides with the finite dimensional version of the quantum flip operator for the Chekhov-Fock quantzation studied in \cite{BBL}.

\vs

In the present paper, we continue the studies of Kashaev \cite{K94} \cite{K99} and Bai \cite{Bai}. I first studied and investigated the left and right dual representations which are one of the standard things to study in the representation theory of Hopf algebras; the action of $\mathcal{W}$ on the dual space of a representation is given as the transpose of the antipode action, or the inverse-antipode action. Although in our case the left and the right duals are isomorphic, they are not identical. Inspired by these two duals, we study `left and right' representation structures on the $\mathrm{Hom}$ spaces of linear maps between representations, and find a natural but nontrivial link between these two; such must be a canonical aspect to be studied for the representation theory of any Hopf algebra, and I find it surprising that it has not been emphasized in the literature thus far, at least written in a language of the present style. Moreover, we find a new answer for a decomposition map $F$ \eqref{eq:the_F-map} in an explicit and `compact' form to be used conveniently for our purposes. As a result, purely from the representation theory of the quantum torus Hopf algebra $\mathcal{W}$ we fully recover the Kashaev-type operators ${\bf A}, {\bf T}, {\bf P}$ on finite dimensional vector spaces, satisfying the consistency relations \eqref{eq:Kashaev_relations} up to constants; see Thm.\ref{thm:main} for a precise formulation of our main theorem. As pointed out to the author by a referee, a finite dimensional version of the Kashaev quantum Teichm\"uller theory, i.e. the operators ${\bf A}, {\bf T}, {\bf P}$, is constructed in \cite{K98} \cite{Kash00b}, and the rigid tensor category structure for the category of representations of $\mathcal{W}$ is studied in \cite{K99} \cite{GKT}, but in a somewhat esoteric language. The present paper provides a concise and clear framework of establishing these ${\bf A}, {\bf T}, {\bf P}$ using only plain representation theoretic language, hence can easily be mimicked for other Hopf algebras.

\vs

The key ingredient of the original quantization works of Kashaev and Chekhov-Fock is a special function called the {\em quantum dilogarithm} of Faddeev-Kashaev \cite{FK94}, or more precisely, the `compact' quantum dilogarithm $\Psi^q(x) = \prod_{i=1}^\infty(1 - q^{2i-1}x)^{-1}$ or the `non-compact' quantum dilogarithm $\Phi^\hbar(z) = \exp ( - \frac{1}{4} \int_\Omega \frac{e^{-ipz}}{\sinh(\pi p)\sinh(\pi \hbar p)} \frac{dp}{p})$, where $\Omega$ is the real line contour avoiding the origin by a small half circle above the origin. The corresponding key ingredient in the present paper is the {\em cyclic quantum dilogarithm} \cite{FK94}, which is a finite version, or a root of unity version. We slightly modify its definition, and study its properties, and suggest more to be studied in the future.

\vs

{\bf Acknowledgments.} This work was supported by the Ewha Womans University Research Grant of 2017. This research was supported by Basic Science Research Program through the National Research Foundation of Korea(NRF) funded by the Ministry of Education(grant number 2017R1D1A1B03030230). H.K. thanks anonymous referees and the editors for their works.

\section{Representations of the quantum torus at root of unity}

In the present section, we mostly collect from the literature basic definitions and properties of representations of the quantum torus algebra. We modify, re-define, and re-establish some of them.

\subsection{Cyclic representations}

\begin{definition}[quantum torus]
\label{def:W}
Let $q$ be either a formal symbol or a nonzero complex number. The \ul{\em quantum torus} $\mathcal{W} = \mathcal{W}^q$ is a Hopf algebra over $\mathbb{C}[q,q^{-1}]$, which is
\begin{align}
\label{eq:quantum_torus_algebra}
\langle \, X^{\pm 1}, Y ~ | ~ XY = q^2 YX \, \rangle
\end{align}
as an algebra over $\mathbb{C}[q,q^{-1}]$, and whose Hopf algebra structure is given by the following coproduct, counit, and antipode:
\begin{align}
\label{eq:W_Hopf_algebra_structure}
\left\{ 
 {\arraycolsep=15pt \begin{array}{ll}
    \Delta(X) = X\otimes X, & \Delta(Y)  = Y \otimes X + 1 \otimes Y ; \\
    \epsilon(X) = 1, & \epsilon(Y) = 0; \\
    S(X) = X^{-1}, & S(Y) = - Y X^{-1}.
  \end{array}}
\right.
\end{align}
\end{definition}

\begin{remark}
The algebra defined in \eqref{eq:quantum_torus_algebra} is often called the Weyl algebra in the literature, including \cite{Bai}. 
\end{remark}
The notation $\mathcal{W}$, and the choice of the particular Hopf algebra structure as above, are taken from \cite{K94}; in particular, the coproduct convention is different from \cite{Bai}. As usual in the theory of Hopf algebras \cite{CP}, a \ul{\em representation} of $\mathcal{W}$ (or a $\mathcal{W}$-representation) is a left $\mathcal{W}$-module, where $\mathcal{W}$ is viewed just as an algebra. We shall only consider the case when $q$ is a nonzero complex number, so a representation is understood as a $\mathbb{C}$-algebra homomorphism $\mu : \mathcal{W} \to \mathrm{End}_\mathbb{C}(V)$ for some complex vector space $V$; action of an element $W$ on a vector $v\in V$ is denoted by $\mu(W)(v)$, or by $W.v$ if $\mu$ is clear from the context. In the present paper, all representations will be assumed to be finite dimensional complex vector spaces. As pointed out in \cite{Bai}, non-trivial finite dimensional representations exist only when $q^2$ is a root of unity. So, from now on, we assume that
\begin{align*}
  \mbox{$q^2$ is a primitive $N$ th root of unity, for some positive integer $N$.}
\end{align*}
More precisely, we could assume that
$$
q = e^{\pi \sqrt{-1} / N}.
$$
In particular, we have $q^N = -1$. For our purposes, we restrict ourselves to the case when
\begin{align}
\label{eq:N_odd}
\mbox{$N$ is odd}.
\end{align}

\begin{definition}[\cite{K94}, \cite{K99}, \cite{Bai}]
A finite dimensional representation $\mu : \mathcal{W} \to \mathrm{End}_\mathbb{C}(V)$ is called \ul{\em cyclic} if $\mu(X)$ and $\mu(Y)$ are invertible.
\end{definition}

Prop.2 of \cite{Bai} says that every cyclic representation of $\mathcal{W}$ is a direct sum of cyclic irreducible representations, and it also completely classifies the cyclic irreducibles. Here, let us come up with a notation for them which is slightly different from \cite{Bai}.
\begin{definition}[cyclic irreducible representations; \cite{Bai}]
\label{def:mu_lambda}
A \ul{\em weight} is a pair $\lambda = (x,y)$ of complex numbers $x,y$. A weight $\lambda = (x,y)$ is said to be \ul{\em non-singular} if both $x$ and $y$ are nonzero.

For a non-singular weight $\lambda = (x,y)$, define the representation $V_\lambda = V_{(x,y)}$ of $\mathcal{W}$ as follows. It is $\mathbb{C}^N$ as a $\mathbb{C}$-vector space, and the action
$$
\mu_\lambda : \mathcal{W} \to \mathrm{End}_\mathbb{C}(\mathbb{C}^N)
$$
is given by
$$
\mu_\lambda(X) = x^{1/N} \cdot A, \quad \mu_\lambda(Y) = y^{1/N} \cdot B,
$$
on the generators, where $x^{1/N}$ and $y^{1/N}$ are arbitrary $N$ th roots of $x$ and $y$ respectively, and $A,B \in \mathrm{End}_\mathbb{C}(\mathbb{C}^N)$ are the unitary matrices defined by
\begin{align}
\label{eq:A_B}
A \, e_k = q^{2k} \, e_k, \qquad B \, e_k = e_{k+1}, \quad \mbox{for} \quad 0\le k \le N-1,
\end{align}
where $\{e_0,e_1,\ldots,e_{N-1}\}$ is a standard basis of $\mathbb{C}^N$, whose subscript indices are considered modulo $N$.
\end{definition}

It is easy to observe
\begin{align}
\label{eq:A_B_relations}
AB = q^2 BA, \qquad A^N = \mathrm{id} = B^N,
\end{align}
as operators; just check on basis vectors. Now, \cite[Prop.2]{Bai} translates to:
\begin{proposition}[\cite{K99}, \cite{Bai}]
\label{prop:cyclic_representations}
In the above definition, any choice of $x^{1/N}$ and $y^{1/N}$ yields isomorphic representation, which is indeed a cyclic irreducible representation, with
$$
\mu_\lambda(X^N) = x \cdot \mathrm{id}_{\mathbb{C}^N}, \qquad \mu_\lambda(Y^N) = y \cdot \mathrm{id}_{\mathbb{C}^N}.
$$
Any cyclic irreducible representation is isomorphic to some $\mu_\lambda$. Two representations $\mu_{(x,y)}$ and $\mu_{(x',y')}$ are isomorphic if and only if $x = x'$ and $y = y'$. Any cyclic representation is a direct sum of cyclic irreducible representations. \qed
\end{proposition}
One can easily understood why $X^N$ and $Y^N$ should act as scalars on an irreducible representation; note that $X^N$ and $Y^N$ commute with all elements of $\mathcal{W}$. We separately write the following fact which is related to the irreducibility of $\mu_\lambda$, as it will be used later on:
\begin{proposition}[irreducibility]
\label{prop:irreducibility}
If $U$ is any linear operator on $\mathbb{C}^N$ commuting with both $A$ and $B$, then $U$ is a scalar operator.
\end{proposition}
{\it Proof.} Write $U e_i = \sum_{j=0}^{N-1} u^j_i e_j$ for some numbers $u^j_i$. Note $AUe_i = A\sum_j u^j_i e_j = \sum_j u^j_i q^{2j} e_j$ and $UA e_i = U q^{2i} e_i = q^{2i} \sum_j u^j_i e_j$, so it must be that $u^j_i q^{2j} = u^j_i q^{2i}$, implying that $u^j_i=0$ unless $i=j$; so $U e_i = u^i_i e_i$. Note $BU e_i = B u^i_i e_i = u^i_i e_{i+1}$ and $UBe_i = U e_{i+1} = u^{i+1}_{i+1} e_{i+1}$, hence $u^i_i = u^{i+1}_{i+1}$, implying that $u^i_i = u^j_j$ for all $i,j$. Writing $u:=u^i_i$ for any $i$, we have $Ue_i = u\cdot e_i$ for all $i$, hence $U = u\cdot \mathrm{id}$. \qed

\begin{corollary}[conjugation action on $A$ and $B$ determines an operator]
\label{cor:conjugation_on_A_B_determines}
An invertible linear operator $U$ on $\mathbb{C}^N$ is completely determined up to scalar by its conjugation action on $A$ and $B$. That is, if $U,U'$ are invertible linear operators on $\mathbb{C}^N$ with $UAU^{-1} = U' A (U')^{-1}$ and $UBU^{-1} = U' B (U')^{-1}$, then $U = \alpha \cdot U'$ for some $\alpha \in \mathbb{C}^\times$. 
\end{corollary}

{\it Proof.} The condition on $U,U'$ implies $A (U^{-1} U') = (U^{-1} U') A$ and $B (U^{-1} U') = (U^{-1} U') B$, so by Prop.\ref{prop:irreducibility} the operator $U^{-1} U'$ is a scalar operator. \qed

\vs

Throughout this paper, we assume that, for each nonzero complex number $z$ we choose its $N$ th root $z^{1/N}$ arbitrarily once and fix it forever; the only condition we impose is:
\begin{align}
\label{eq:normalization_for_N_th_roots}
1^{1/N} = 1, \qquad (z^{1/N})^{-1} = (z^{-1})^{1/N} \quad\mbox{and}\quad (-z)^{1/N} = - z^{1/N}, \quad \forall z\in \mathbb{C}^\times;
\end{align}
in particular, the expression $z^{-1/N}$ makes sense. For any $z,w\in \mathbb{C}^\times$, both $z^{1/N} w^{1/N}$ and $(zw)^{1/N}$ are $N$ th roots of $zw$, hence differ by multiplication by an $N$ th root of unity. Hence there exists a unique number $m_{z,w} \in \{0,1,\ldots,N-1\}$ such that
\begin{align}
\label{eq:m_definition}
(zw)^{1/N} = q^{2m_{z,w}} \, z^{1/N} w^{1/N}.
\end{align}
Note that the assignment $(z,w) \mapsto m_{z,w}$ is symmetric, i.e. $m_{z,w}=m_{w,z}$, and satisfies the cocycle condition: $m_{z,w} + m_{zw,u} = m_{w,u} + m_{z,wu}$
.

\subsection{The Clebsch-Gordan operators for regular pairs}
\label{subsec:CG}

It is standard in the theory of Hopf algebras \cite{CP} to define the tensor product of two or more representations of $\mathcal{W}$ using the coproduct of $\mathcal{W}$, and to define a \ul{\em trivial representation} of $\mathcal{W}$ by a $\mathcal{W}$-representation where every element of $W\in \mathcal{W}$ acts by counit $\epsilon(W) \cdot \mathrm{id}$. In particular, action of $W\in \mathcal{W}$ on the tensor product $V \otimes V'$ of representations $V$ and $V'$ is given by
\begin{align}
\label{eq:action_on_tensor_product}
W.(v\otimes v') = \sum (W_{(1)}.v) \otimes (W_{(2)}.v'),
\end{align}
where we used the Sweedler notation $\Delta(W) = \sum W_{(1)} \otimes W_{(2)}$ for the coproduct.

\begin{definition}[\cite{K99}, \cite{Bai}]
A sequence $(\mu_1,\mu_2,\ldots,\mu_n)$ of $\mathcal{W}$-representations is said to be \ul{\em regular} if the representation $\mu_i \otimes \mu_{i+1} \otimes \cdots \otimes \mu_j$ is cyclic for any $1\le i\le j \le n$.
\end{definition}
\begin{definition}
A sequence $(\lambda_1,\ldots,\lambda_n)$ of weights is said to be \ul{\em regular} if the sequence $(\mu_{\lambda_1},\ldots,\mu_{\lambda_n})$ of $\mathcal{W}$-representations is regular.
\end{definition}
As pointed out in \cite{K94}, \cite{K99} and \cite{Bai}, it is easy to observe
$$
\Delta(X^N) = X^N \otimes X^N, \qquad
\Delta(Y^N) = Y^N \otimes X^N + 1 \otimes Y^N;
$$
use $\sum_{i=0}^{N-1} q^{2i} = 0$. It then follows that, for any non-singular weights $\lambda=(x,y)$ and $\lambda'=(x',y')$, the elements $X^N$ and $Y^N$ of $\mathcal{W}$ acts on the tensor product representation $V_\lambda \otimes V_{\lambda'}$ as $(xx')\cdot \mathrm{id}$ and $(yx'+y')\cdot \mathrm{id}$ respectively.  In fact, we have more:
\begin{lemma}[see \cite{Bai}, and \cite{K94}, \cite{K99}; multiplication of weights]
\label{lem:regular_pair}
Define the multiplication of any two weights $\lambda=(x,y)$ and $\lambda' = (x',y')$ as
\begin{align}
\label{eq:lambda_multiplication}
\lambda \, \lambda' := (\,x \, x', \, y \, x' + y' \,).
\end{align}
If $(\lambda,\lambda')$ is a regular pair of weights, then $\mu_\lambda \circ \mu_{\lambda'}$ is a direct sum of $N$ copies of the cyclic irreducible representation $\mu_{\lambda \, \lambda'}$. \qed
\end{lemma}

The following is an easy observation, but not found in \cite{K94}, \cite{K99}, or  \cite{Bai}.
\begin{lemma}[regular pair criterion]
A pair $(\lambda,\lambda')$ of weights is regular if and only if all three weights $\lambda,\lambda',\lambda\, \lambda'$ are non-singular.
\end{lemma}

{\it Proof.} Let $\lambda = (x,y)$ and $\lambda' = (x',y')$ be weights. Assume first that $(\lambda,\lambda')$ is a regular pair; in particular, each of $\mu_\lambda$ and $\mu_{\lambda'}$ must be cyclic, hence $\lambda,\lambda'$ are non-singular. Then, automatically, $x x' \neq 0$. Since $\Delta X = X\otimes X$, the action of $X\in \mathcal{W}$ on $V_\lambda \otimes V_{\lambda'}$ is given by $X  = (\mu_\lambda \otimes \mu_{\lambda'})(X) = x^{1/N} \, (x')^{1/N} A\otimes A$, because $X.(v\otimes v') = (X\otimes X)(v\otimes v') = (X.v)\otimes (X.v') = (x^{1/N} Av)\otimes ((x')^{1/N} Av')$; in particular, the $X$-action on $V_\lambda \otimes V_{\lambda'}$ is already invertible. Likewise, since $\Delta Y =  Y\otimes X + 1\otimes Y$, the action of $Y\in \mathcal{W}$ on $V_\lambda\otimes V_{\lambda'}$ is given by $Y = (\mu_\lambda\otimes\mu_{\lambda'})(Y)= y^{1/N} (x')^{1/N} B \otimes A + (y')^{1/N} \mathrm{id}\otimes B$. In order for this $Y$ operator to be invertible, it suffices that it has zero kernel. Suppose $\sum_{i,j} f_{i,j} e_i \otimes e_j \in \mathbb{C}^N \otimes \mathbb{C}^N = V_\lambda \otimes V_{\lambda'}$ is in the kernel of $Y$, for some numbers $f_{i,j}$. Then
\begin{align*}
0 = Y.\sum_{i,j} f_{i,j} e_i \otimes e_j & = \sum_{i,j} f_{i,j} ( y^{1/N} (x')^{1/N} (Be_i) \otimes (Ae_j) + (y')^{1/N} e_i \otimes (Be_j)) \\
& = \sum_{i,j} f_{i,j} ( y^{1/N} (x')^{1/N} q^{2j} \, e_{i+1} \otimes e_j + (y')^{1/N} \, e_i \otimes e_{j+1}) \\
& = \sum_{i,j} (f_{i-1,j} y^{1/N} (x')^{1/N} q^{2j} + f_{i,j-1} (y')^{1/N}) \, e_i \otimes e_j,
\end{align*}
which happens if and only if $0 = f_{i-1,j} y^{1/N} (x')^{1/N} q^{2j} + f_{i,j-1} (y')^{1/N}$ for all $i,j$. This condition is in turn equivalent to
\begin{align}
\label{eq:f_kernel_condition}
f_{i,j}  = - f_{i-1,j+1} \cdot \frac{ y^{1/N} (x')^{1/N} }{(y')^{1/N}} \cdot q^{2(j+1)}, \quad \forall i,j.
\end{align}
Using this recursively $N$ times, with $f_{i-N,j+N} = f_{i,j}$, one obtains $f_{i,j} = f_{i,j} \cdot (-1)^N \frac{y x'}{y'} \, q^{2(1+2+\cdots+(N-1)+N)}  = - f_{i,j} \cdot \frac{yx'}{y'} q^{N(N+1)} = - f_{i,j}\cdot yx'/y'$ for each $i,j$, where we used the fact \eqref{eq:N_odd} that $N$ is odd. Hence, if $-yx'/y' \neq 1$, then $\sum_{i,j} f_{i,j} e_i \otimes e_j \in \ker Y$ implies $f_{i,j}=0$ for all $i,j$, therefore $\ker Y =0$ as desired. On the other hand, if $-yx' / y' = 1$, then one can find numbers $f_{i,j}$ so that $\sum_{i,j} f_{i,j} e_i \otimes e_j$ is a nonzero element of $\ker Y$. For example, let $f_{0,0} := 1$, $f_{i,-i} := (-\frac{y^{1/N}(x')^{1/N}}{(y')^{1/N}})^i \, q^{- i(i-1)}$ for all $i=1,2,\ldots,N-1$, and $f_{i,j}=0$ otherwise; one can easily verify that \eqref{eq:f_kernel_condition} is satisfied, hence $\sum_{i,j} f_{i,j} e_i \otimes e_j$ is indeed a nonzero element in $\ker Y$. 

\vs

Summarizing, the regularity assumption of $(\lambda,\lambda')$ implies $\ker Y = 0$, hence $yx' + y'\neq 0$, so $\lambda \, \lambda' = (xx',yx'+y')$ is non-singular. Conversely, if $\lambda=(x,y)$, $\lambda'=(x',y')$, and $\lambda \, \lambda' =(xx',yx'+y')$ are all non-singular, then $\mu_\lambda$ and $\mu_{\lambda'}$ are cyclic, and as seen above, the operators for $X$ and $Y$ on $V_\lambda \otimes V_{\lambda'}$ are invertible (we have $\ker Y=0$), hence $\mu_\lambda \otimes \mu_{\lambda'}$ is cyclic. \qed

\vs

For a regular triple $(\lambda,\lambda',\lambda'')$ of weights representations, one has $(\lambda \, \lambda') \, \lambda'' = \lambda \, (\lambda' \, \lambda'')$, from the coassociativity of the coproduct. In fact, for this associativity of products of weights, we do not need regularity, which is easy to check. So we omit parentheses in products of weights. We collect a couple more observations:
\begin{lemma}[properties of multiplication of weights]
The multiplication of weights defined in \eqref{eq:lambda_multiplication} is associative, but not commutative. The weight $\mathbf{1} := (1,0)$ is the multiplicative identity, i.e. $\mathbf{1} \, \lambda = \lambda \, \mathbf{1}$ for each weight $\lambda$.
\end{lemma}

\begin{corollary}[regular sequence criterion]
A sequence $(\lambda_1,\ldots,\lambda_n)$ of weights is regular if and only if the product $\lambda_i \, \lambda_{i+1} \, \cdots \, \lambda_j$ is non-singular for any $1\le i \le j \le n$. \qed
\end{corollary}

The notation in \cite{K94}, \cite{K99} and \cite{Bai} is:
$$
\mbox{when $\mu = \mu_\lambda$ and $\nu = \mu_{\lambda'}$, write $\mu \nu := \mu_{\lambda \, \lambda'}$};
$$
notice that Bai \cite{Bai} used $\lambda \, \lambda' = (xx',x^{-1} y' + y)$, due to the different choice of coproduct.
\begin{definition}[\cite{K94}, \cite{K99}, \cite{Bai}]
For a regular pair $(\mu,\nu)$ of irreducible representations, define
$$
V_{(\mu,\nu)} := \mathrm{Hom}_\mathcal{W}(V_{\mu\nu}, V_\mu\otimes V_\nu) = \{ f : V_{\mu\nu} \to V_\mu\otimes V_\nu \, | \, \mbox{$f$ intertwines the $\mathcal{W}$ actions} \}.
$$
The elements of $V_{(\mu,\nu)}$ are called the \ul{\em Clebsch-Gordan operators}. The map
$$
\Omega = \Omega_{(\mu,\nu)} : V_{(\mu,\nu)} \otimes V_{\mu\nu} \to V_\mu \otimes V_\nu,
$$
defined by
$$
\Omega(f,v) := f(v), \quad \forall f\in V_{(\mu,\nu)}, \quad \forall v\in V_{\mu\nu},
$$
is called the \ul{\em canonical map}.
\end{definition}
From Lem.\ref{lem:regular_pair} it follows that
\begin{align}
\label{eq:dimension_of_intertwiner_space}
\dim_\mathbb{C}V_{(\mu,\nu)} = N.
\end{align}
As defined, $V_{(\mu,\nu)}$ can be viewed as a `space of intertwiners'.

\vs

Here, in the style of \cite{FK12}, we would like to understand $V_{(\mu,\nu)}$ as a `multiplicity space' which is in particular a trivial $\mathcal{W}$-module, and choose an explicit realization of $V_{(\mu,\nu)}$ as $\mathbb{C}^N$ as a vector space. To avoid confusion, we come up with a separate notation.

\begin{definition}[the multiplicity space]
Let $(\lambda,\lambda')$ be a regular pair of weights. Let $M_{\lambda,\lambda'}^{\lambda\lambda'}$ be $\mathbb{C}^N$ as a vector space with the standard basis $e_0,\ldots,e_{N-1}$, and let it be a trivial representation of the Hopf algebra $\mathcal{W}$.
\end{definition}

One could explicitly compute a basis of $\mathrm{Hom}_\mathcal{W}(V_{\lambda\, \lambda'}, V_\lambda\otimes V_{\lambda'})$, and thus obtain an identification of $V_{(\mu_\lambda,\mu_{\lambda'})}$ and $\mathbb{C}^N \equiv M^{\lambda\lambda'}_{\lambda,\lambda'}$. Instead, we shall construct a decomposition map of $V_\lambda \otimes V_{\lambda'}$ into the direct sum of $N$ copies of $V_{\lambda\lambda'}$, realized as the following isomorphism of $\mathcal{W}$-modules:
$$
V_\lambda \otimes V_{\lambda'} \to M^{\lambda\lambda'}_{\lambda,\lambda'} \otimes V_{\lambda  \lambda'},
$$
replacing the canonical map $\Omega$ above. Thus one understands $M^{\lambda\lambda'}_{\lambda,\lambda'}$ as the space of multiplicities (of $V_{\lambda \lambda'}$ in $V_\lambda \otimes V_{\lambda'}$). Such isomorphism is not unique, and Kashaev constructs one answer in \cite{K94} \cite{K99}. Shortly, using Faddeev-Kashaev's `cyclic quantum dilogarithm' \cite{K94} \cite{FK94}, we shall construct another answer which is written more neatly and has a favorable property.

\subsection{Cyclic quantum dilogarithm}
\label{subsec:Phi_C}

We review the notion of cyclic quantum dilogarithm of Faddeev and Kashaev \cite{FK94}. In quantum Teichm\"uller theory and related subjects, a special function named the quantum dilogarithm, established in \cite{FK94}, is used; a key intertwining operator is expressed by applying the functional calculus for a certain self-adjoint operator on this function. However, in case when $q$ is a root of unity, it is best to understand the `cyclic quantum dilogarithm' as a certain functional-calculus-type construction, instead of as a function.

\begin{definition}[\cite{FK94}, \cite{BB}; modified to two-parameter]
\label{def:w}
For $n \in \{0,1,\ldots,N-1\}$ and two complex numbers ${\bf a}, {\bf c}$ satisfying the equation
$$
{\bf c}^N - {\bf a}^N = 1,
$$
define the function $w({\bf a}, {\bf c}|n)$ as
$$
w({\bf a}, {\bf c}|0)=1, \qquad
w({\bf a}, {\bf c}|n) = \prod_{j=1}^n ({\bf c} - {\bf a} \, q^{2j})^{-1}, \quad \forall n\ge 1.
$$
\end{definition}

Using the equation ${\bf c}^N - {\bf a}^N=1$, one observes $w({\bf a},{\bf c}|n+N) = w({\bf a}, {\bf c}|n)$, so $n$ can be considered as an integer modulo $N$. The expression $({\bf c} - {\bf a} \, q^{2j})^{-1}$ always makes sense, because ${\bf c} - {\bf a} \, q^{2j} = 0$ implies ${\bf c} = {\bf a} \,  q^{2j}$, hence ${\bf c}^N = {\bf a}^N$, violating the equation ${\bf c}^N - {\bf a}^N = 1$. 

\vs

We note that in the original papers \cite{FK94} \cite{BB}, a three parameter function is used $w({\bf a}, {\bf c}|n) = \prod_{j=1}^n \frac{b}{c-aq^{2j}}$ for ${\bf a}, {\bf c}$ satisfying $a^N + b^N = c^N$; I changed it here because everything depends only on the ratios $\frac{a}{b}$ and $\frac{c}{b}$. For their notation, we have $w({\bf a}, {\bf c}|n) = w({\bf a}, 1, {\bf c}|n)$. Another advantage of this new definition is that $w({\bf a}, {\bf c}|n)$ is never zero.

\vs

For the next definition, first observe by a simple linear algebra that any linear operator $C$ on a finite dimensional complex vector space satisfying $C^N=\mathrm{id}$ is diagonalizable with eigenvalues in $\{q^{2i} : i=0,1,\ldots,N-1\}$.
\begin{definition}[modified from \cite{FK94}: the cyclic quantum dilogarithm]
\label{def:Phi_C}
Let ${\bf a}, {\bf c}$ be as in Def.\ref{def:w}. Suppose $C$ is a linear operator on a complex vector space, say $V$, with $C^N=\mathrm{id}_V$. The linear operator $\Phi(C) = \Phi^q_{{\bf a}, {\bf c}}(C)$ on $V$ is defined as the operator acting as $\zeta_{{\bf a}, {\bf c}} \cdot w({\bf a}, {\bf c}|i) \cdot \mathrm{id}$ on the $q^{2i}$-eigenspace of $C$, where the nonzero complex number $\zeta_{{\bf a}, {\bf c}}$ is to be determined later. This operator $\Phi(C)$ is called a \ul{\em cyclic quantum dilogarithm operator}.
\end{definition}

One can view the definition of $\Phi(C)$ as being an application of a finite version of `functional calculus' for $C$. For the purposes of the present paper, it is a bit modified from \cite{FK94}; there, they construct $\Psi(C)$ for operator $C$ satisfying $C^N = - \mathrm{id}_N$.  The following lemma shows that $\Phi^q_{{\bf a}, {\bf c}}(\sim)$ indeed behaves like functional calculus.
\begin{lemma}[commuting with conjugation]
\label{lem:commuting_property}
If ${\bf a}, {\bf c}, C$ are as in  Def.\ref{def:Phi_C} and $D,C'$ are linear operators on $V$ such that $DC = C' D$ and $(C')^N=\mathrm{id}$, then
$$
D \,\Phi^q_{{\bf a}, {\bf c}}(C) = \Phi^q_{{\bf a}, {\bf c}}(C')\, D.
$$
When $D$ is invertible, one has $C' = DCD^{-1}$ and
$$
D\, \Phi^q_{{\bf a}, {\bf c}}(C) \, D^{-1}  = \Phi^q_{{\bf a}, {\bf c}}(D\,C\,D^{-1}).
$$
\end{lemma}

{\it Proof.} It suffices to check the equality on each eigenvector of $C$, as $C$ is diagonalizable. Let $v_i$ be a $q^{2i}$-eigenvector of $C$; then $Dv_i$ is a $q^{2i}$-eigenvector of $C'$, because $C'Dv_i = DCv_i = Dq^{2i} v_i = q^{2i} Dv_i$. So $D \Phi^q_{{\bf a}, {\bf c}}(C) v_i = D \zeta_{{\bf a}, {\bf c}} w({\bf a}, {\bf c}|i) v_i = \zeta_{{\bf a}, {\bf c}}w({\bf a}, {\bf c}|i) Dv_i$, while $\Phi^q_{{\bf a}, {\bf c}}(C')D v_i= \zeta_{{\bf a}, {\bf c}}w({\bf a}, {\bf c}|i) Dv_i$. Hence $D\Phi^q_{{\bf a}, {\bf c}}(C) v_i = \Phi^q_{{\bf a}, {\bf c}}(C')D v_i$, as desired. \qed

\vs

The following is the `defining' functional relation of $\Phi$.
\begin{lemma}[see \cite{FK94}]
\label{lem:functional_relation}
Let ${\bf a}, {\bf c},C$ be as in Def.\ref{def:Phi_C}. Then, $\Phi(C) = \Phi^q_{{\bf a}, {\bf c}}(C)$ is invertible and satisfies
$$
\Phi^q_{{\bf a}, {\bf c}}(q^{-2}C) \, \Phi^q_{{\bf a}, {\bf c}}(C)^{-1} = {\bf c} \cdot \mathrm{id}- {\bf a} \cdot C.
$$
\end{lemma}

{\it Proof.} Each assertion can be checked on each $q^{2i}$-eigenvector of $C$, say $v_i$. The invertibility holds because $\zeta_{{\bf a}, {\bf c}} \cdot w({\bf a}, {\bf c}|i)$ is a nonzero number, hence invertible. For the last assertion, note first that $v_i$ is then a $q^{2(i-1)}$-eigenvector of $q^{-2} C$. So, observe
\begin{align*}
  \Phi^q_{{\bf a}, {\bf c}}(q^{-2} C) \, \Phi^q_{{\bf a}, {\bf c}}(C)^{-1} \, v_i & = \Phi^q_{{\bf a}, {\bf c}}(q^{-2} C) ( \zeta_{{\bf a}, {\bf c}}^{-1} \cdot w({\bf a}, {\bf c}|i)^{-1} \cdot v_i) \\
& = w({\bf a}, {\bf c}|i-1) \cdot w({\bf a}, {\bf c}|i)^{-1} \cdot v_i
= ({\bf c} - {\bf a} q^{2i}) \cdot v_i  \\ 
& \textstyle = ({\bf c} \cdot v_i - {\bf a} \cdot C v_i) =  ({\bf c}\cdot \mathrm{id} - {\bf a} \cdot C) \, v_i. \qed
\end{align*}

Faddeev and Kashaev \cite{FK94} assert without proof that this functional equation uniquely determines $\Phi^q_{{\bf a}, {\bf c}}(C)$ up to multiplicative constant; we do not attempt to prove it here.

\begin{lemma}[conjugation by cyclic quantum dilogarithm]
\label{lem:conjugation_by_Phi_C}
If ${\bf a}, {\bf c},C$ are as in Def.\ref{def:Phi_C} and $D$ is a linear operator on $V$, then
\begin{align*}
& D \circ \Phi^q_{{\bf a}, {\bf c}}(C) = \Phi^q_{{\bf a}, {\bf c}} (C)\circ ( {\bf c}\, D - {\bf a} \, CD), \qquad\quad \mbox{if}\quad  CD=q^2DC, \\
& D \circ \Phi^q_{{\bf a}, {\bf c}}(C)^{-1} = \Phi^q_{{\bf a}, {\bf c}} (C)^{-1} \circ ( {\bf c} \, D - {\bf a} \, DC),\quad \mbox{if}\quad  CD=q^{-2}DC.
\end{align*}
\end{lemma}

{\it Proof.} Suppose $CD = q^2 DC$. As $DC = (q^{-2} C)D$ and $(q^{-2} C)^N=\mathrm{id}$, note
\begin{align*}
 D\, \Phi^q_{{\bf a}, {\bf c}}(C) \, \stackrel{{\rm Lem}.\ref{lem:commuting_property}}{=} \, \Phi^q_{{\bf a}, {\bf c}}(q^{-2} C) \, D \, \stackrel{{\rm Lem}.\ref{lem:functional_relation}}{=} \, \textstyle ({\bf c}\cdot\mathrm{id} - {\bf a}\,C) \, \Phi^q_{{\bf a}, {\bf c}}(C) \, D
\, \stackrel{{\rm Lem}.\ref{lem:commuting_property}}{=} \, \Phi^q_{{\bf a}, {\bf c}}(C) \, ({\bf c}\cdot \mathrm{id} - {\bf a}\, C) \, D.
\end{align*}
Suppose $CD = q^{-2} DC$. As $DC = (q^2 C)D$ and $(q^2 C)^N = \mathrm{id}$, note
\begin{align*}
\textstyle D \, \Phi^q_{{\bf a}, {\bf c}}(C)^{-1} \,\stackrel{{\rm Lem}.\ref{lem:commuting_property}}{=}\, \Phi^q_{{\bf a}, {\bf c}}(q^2C)^{-1} D
\,\stackrel{{\rm Lem}.\ref{lem:functional_relation}}{=}\, \Phi^q_{{\bf a}, {\bf c}}(C)^{-1} \, ({\bf c}\cdot \mathrm{id}-{\bf a}\, q^2C) D. \qed
\end{align*}




\vs

One can write the above as the `conjugation action', written as
\begin{align}
\label{eq:conjugation_by_Phi_C}
\left\{
{\renewcommand{\arraystretch}{1.2} \begin{array}{ll}
\Phi^q_{{\bf a}, {\bf c}}(C)^{-1} \, D \, \Phi^q_{{\bf a}, {\bf c}}(C) = \textstyle ({\bf c}\cdot\mathrm{id} - {\bf a} \, C) \circ D, & \mbox{if $CD = q^2DC$}, \\
\Phi^q_{{\bf a}, {\bf c}}(C) \, D \, \Phi^q_{{\bf a}, {\bf c}}(C)^{-1} = \textstyle D \circ  ({\bf c} \cdot \mathrm{id} - {\bf a} \, C), & \mbox{if $CD = q^{-2} DC$}.
\end{array}}
\right.
\end{align}

\vs

The cyclic quantum dilogarithm $\Phi(C)$ satisfies a very important 5-factor identity called the quantum pentagon identity; we postpone the discussion of this until \S\ref{subsec:pentagon_identity_of_Phi_C}.

\subsection{Explicit decomposition of a tensor product}

Using the results of the previous subsection, we establish the explicit decomposition map of the tensor product representation $V_\lambda \otimes V_{\lambda'}$ into direct sum of $N$ copies of $V_{\lambda \, \lambda'}$

\begin{proposition}[explicit decomposition map of a tensor product]
\label{prop:F}
Let $(\lambda,\lambda')$ be a regular pair of weights; write $\lambda=(x,y)$, $\lambda'=(x',y')$. Define the map $\mathbf{F}_{\lambda,\lambda'} : \mathbb{C}^N \otimes \mathbb{C}^N \to \mathbb{C}^N \otimes \mathbb{C}^N$ as
\begin{align}
\label{eq:F}
\mathbf{F}_{\lambda,\lambda'} := (B_2)^{-m_{x,x'}} \circ \mathbf{S} \circ \Phi^q_{\lambda,\lambda'}(B_1 A_2 B_2^{-1}),
\end{align}
where $\mathbf{S} : \mathbb{C}^N \otimes \mathbb{C}^N \to \mathbb{C}^N \otimes \mathbb{C}^N$ is the unique linear map given on the basis by
\begin{align}
\label{eq:S}
\mathbf{S} : e_i \otimes e_j \mapsto e_i \otimes e_{i+j}, \quad \forall i,j \in \{0,\ldots,N-1\},
\end{align}
while
$$
A_1 = A\otimes \mathrm{id}, \quad B_1 = B \otimes \mathrm{id}, \quad A_2 = \mathrm{id}\otimes A, \quad \mbox{and}\quad B_2 = \mathrm{id} \otimes B,
$$
with the operators $A$ and $B$ on $\mathbb{C}^N$ are as in \eqref{eq:A_B}, $m_{x,x'}$ is as in \eqref{eq:m_definition}, and
\begin{align}
\label{eq:our_ac}
\Phi^q_{\lambda,\lambda'}: = \Phi^q_{{\bf a}, {\bf c}} \quad\mbox{with}\quad {\bf a} = - \frac{ y^{1/N} \, (x')^{1/N} }{(yx'+y')^{1/N}}, \quad {\bf c}=\frac{(y')^{1/N}}{(yx'+y')^{1/N}},
\end{align}
where $\Phi^q_{{\bf a}, {\bf c}}$ is the cyclic quantum dilogarithm (Def.\ref{def:Phi_C}). Then, $\mathbf{F}_{\lambda,\lambda'}$ is invertible and provides a $\mathcal{W}$-intertwining map
\begin{align}
\label{eq:F_conceptual}
\mathbf{F}_{\lambda,\lambda'} : V_\lambda \otimes V_{\lambda'} \to M^{\lambda\lambda'}_{\lambda,\lambda'} \otimes V_{\lambda \, \lambda'},
\end{align}
i.e. the intertwining equations hold:
\begin{align}
\label{eq:intertwining_property_F}
\mathbf{F}_{\lambda,\lambda'}( W.(v\otimes v')) = (\mathrm{id}\otimes W).\mathbf{F}_{\lambda,\lambda'}(v\otimes v'),
\end{align}
for all $W\in \mathcal{W}$, $v\in V_\lambda$, and $v' \in V_{\lambda'}$.
\end{proposition}

{\it Proof.} First, the expression $\Phi^q_{\lambda,\lambda'}(B_1A_2B_2^{-1}) = \Phi^q_{{\bf a}, {\bf c}}(B_1A_2B_2^{-1})$ makes sense, because ${\bf c}^N - {\bf a}^N = \frac{ y' }{yx' + y'} + \frac{ yx' }{yx'+y'} = 1$ and $(B_1 A_2 B_2^{-1})^N = q^{-2(1+2+\cdots+(N-1))} B_1^N A_2^N B_2^{-N} = q^{-(N-1)N}\cdot \mathrm{id}_{\mathbb{C}^N \otimes \mathbb{C}^N} = \mathrm{id}_{\mathbb{C}^N \otimes \mathbb{C}^N}$; we used the relations \eqref{eq:A_B_relations} of the operators $A$ and $B$, the fact that operators with different subscripts commute (this is because they act on different tensor factors), as well as the fact \eqref{eq:N_odd} that $N$ is odd.

\vs

Note that $\mathbf{S}$ is invertible because the linear map 
\begin{align}
\label{eq:S_inverse}
\mathbb{C}^N \otimes \mathbb{C}^N \to \mathbb{C}^N \otimes \mathbb{C}^N : e_i \otimes e_j \mapsto e_i \otimes e_{j-i}, \quad \forall i,j
\end{align}
is the inverse of $\mathbf{S}$, while $(B_2)^{-m_{x,x'}}$ is invertible with inverse $(B_2)^{m_{x,x'}}$, and $\Phi^q_{{\bf a}, {\bf c}}(-B_1A_2B_2^{-1})$ is invertible because ${\bf a}, {\bf c}$ are nonzero (Lem.\ref{lem:functional_relation}).

\vs

Let us check the intertwining property \eqref{eq:intertwining_property_F}. Note
\begin{align}
\label{eq:S_conjugation_1}
A_2 \, \mathbf{S} = \mathbf{S} \, A_1 A_2, \qquad
B_2 \, \mathbf{S} = \mathbf{S} \, B_2.
\end{align}
Indeed, on a basis vector $e_i \otimes e_j$, note $(\mathrm{id}\otimes A) \mathbf{S} \, e_i \otimes e_j = (\mathrm{id}\otimes A) \, e_i \otimes e_{i+j} = q^{2(i+j)} e_i \otimes e_{i+j}$ and $\mathbf{S} \, (A\otimes A) \, e_i \otimes e_j = \mathbf{S} \, q^{2(i+j)} e_i \otimes e_j = q^{2(i+j)} e_i \otimes e_{i+j}$ hence the first assertion. Note $(\mathrm{id}\otimes B)\mathbf{S} \, e_i \otimes e_j = (\mathrm{id}\otimes B) \, e_i \otimes e_{i+j} = e_i \otimes e_{i+j+1}$ and $\mathbf{S} \, (\mathrm{id}\otimes B) \, e_i \otimes e_j = \mathbf{S} \, e_i \otimes e_{j+1} = e_i \otimes e_{i+j+1}$, hence the second assertion.

\vs

Note now from $AB = q^2 BA$ that
$$
(A_1 A_2)(B_1 A_2 B_2^{-1}) = (B_1 A_2 B_2^{-1})(A_1A_2), \quad
B_2 (B_1A_2 B_2^{-1}) = q^{-2} (B_1A_2B_2^{-1})B_2.
$$
So $A_1 A_2$ commutes with $\Phi^q_{{\bf a}, {\bf c}}(B_1A_2B_2^{-1})$ by Lem.\ref{lem:commuting_property}, while from Lem.\ref{lem:conjugation_by_Phi_C} we have
\begin{align*}
B_2 \, \Phi^q_{{\bf a}, {\bf c}} (B_1 A_2 B_2^{-1}) & = \Phi^q_{{\bf a}, {\bf c}}(B_1A_2B_2^{-1}) \, ({\bf c}\, B_2 - {\bf a} \, B_1 A_2 B_2^{-1} B_2) \\
& = \Phi^q_{{\bf a}, {\bf c}}(B_1A_2B_2^{-1}) \, ({\bf c} \, B_2 - {\bf a} \, B_1 A_2).
\end{align*}
Also, from $AB = q^2 BA$ we have
\begin{align}
\label{eq:B_m_commuting}
AB^m = q^{2m} B^m A, \quad \mbox{for any integer $m$}.
\end{align}

\vs

Combining the results so far, we have
\begin{align*}
( (xx')^{1/N} A_2 ) \, \mathbf{F}_{\lambda,\lambda'} & = (xx')^{1/N} \, A_2 \, B_2^{-m_{x,x'}} \, \mathbf{S} \, \Phi(B_1A_2B_2^{-1}) \\
& = q^{-2m_{x,x'}} (xx')^{1/N} B_2^{-m_{x,x'}} \, A_2 \, \mathbf{S}\, \Phi(B_1A_2B_2^{-1}) \\
& = x^{1/N} (x')^{1/N} B_2^{-m_{x,x'}} \, \mathbf{S} \, A_1 A_2 \, \Phi(B_1A_2B_2^{-1}) \\
& = x^{1/N} (x')^{1/N} \, B_2^{-m_{x,x'}} \, \mathbf{S}\, \Phi(B_1A_2B_2^{-1}) \, A_1 A_2 \\
& = \mathbf{F}_{\lambda,\lambda'} \, (x^{1/N} (x')^{1/N} A_1 A_2),
\end{align*}
while
\begin{align*}
( (y x' + y')^{1/N} B_2) \, \mathbf{F}_{\lambda,\lambda'} 
& = (y x' + y')^{1/N} B_2 \, B_2^{-m_{x,x'}} \, \mathbf{S} \, \Phi(B_1A_2B_2^{-1}) \\
& = (y x' + y')^{1/N} B_2^{-m_{x,x'}} \, B_2 \, \mathbf{S} \, \Phi(B_1A_2B_2^{-1}) \\
& = (y x' + y')^{1/N} B_2^{-m_{x,x'}} \, \mathbf{S} \, B_2 \, \Phi(B_1A_2B_2^{-1}) \\
& = (y x' + y')^{1/N} B_2^{-m_{x,x'}} \, \mathbf{S} \, \Phi(B_1A_2B_2^{-1}) \, ({\bf c} \, B_2 - {\bf a} \, B_1 A_2) \\
& = \mathbf{F}_{\lambda,\lambda'} \, ( y^{1/N} B_1 \, (x')^{1/N} A_2 + (y')^{1/N} B_2),
\end{align*}
where we put in the values \eqref{eq:our_ac}. For $v\in V_\lambda$, $v' \in V_{\lambda'}$, note from Def.\ref{def:mu_lambda}, \eqref{eq:action_on_tensor_product}, and \eqref{eq:lambda_multiplication} that
\begin{align*}
& X.(v\otimes v') = (X.v)\otimes(X.v') = x^{1/N} (x')^{1/N} \, (A\otimes A)(v\otimes v') = x^{1/N} (x')^{1/N}\, A_1A_2\,(v\otimes v'), \\
& Y.(v\otimes v') = (Y.v)\otimes (X.v') + v\otimes (Y.v') = (y^{1/N} (x')^{1/N} \, B_1 \, A_2 + (y')^{1/N} \, B_2 ) \, (v\otimes v'), \\
& (\mathrm{id}\otimes X).\mathbf{F}_{\lambda,\lambda'}(v\otimes v') = (xx')^{1/N} \, A_2 \, \mathbf{F}_{\lambda,\lambda'}(v\otimes v'), \\
& (\mathrm{id}\otimes Y).\mathbf{F}_{\lambda,\lambda'}(v\otimes v') = (yx' + y')^{1/N} \, B_2 \, \mathbf{F}_{\lambda,\lambda'}(v\otimes v').
\end{align*}
Hence we proved \eqref{eq:intertwining_property_F} for $W = X, Y$, which suffices. \qed

\vs

For the purposes of the present paper, the answer \eqref{eq:F} we just found is more convenient to deal with than Kashaev's answer in eq.(1.13) of \cite{K94}, and in eq.(22)--(24) of \cite{K99}. The inverse of our $\mathbf{F}_{\lambda,\lambda'}$ can be written neatly too, just as
\begin{align}
\label{eq:F_inverse}
\mathbf{F}_{\lambda,\lambda'}^{-1} = \Phi^q_{\lambda,\lambda'}(B_1A_2B_2^{-1})^{-1} \, \mathbf{S}^{-1} \, (B_2)^{m_{x,x'}}. 
\end{align}

\vs

The non-uniqueness of the decomposition isomorphism \eqref{eq:F_conceptual} is studied as follows.
\begin{proposition}
\label{prop:non-uniqueness_of_F}
Let $(\lambda,\lambda')$ be a regular pair of weights. For any $\mathbf{U} \in \mathrm{End}_\mathbb{C}(\mathbb{C}^N)$, the map $(\mathbf{U}\otimes \mathrm{id}) \circ \mathbf{F}_{\lambda,\lambda'}$ is a $\mathcal{W}$-module isomorphism map $V_\lambda \otimes V_{\lambda'} \to M^{\lambda\lambda'}_{\lambda,\lambda'} \otimes V_{\lambda \, \lambda'}$, and any $\mathcal{W}$-module isomorphism $V_\lambda \otimes V_{\lambda'} \to M^{\lambda\lambda'}_{\lambda,\lambda'} \otimes V_{\lambda \, \lambda'}$ is of this form.
\end{proposition}

{\it Proof.} The first assertion is easy; since $M^{\lambda\lambda'}_{\lambda,\lambda'}$ is a trivial $\mathcal{W}$-module, any invertible linear map $M^{\lambda\lambda'}_{\lambda,\lambda'} \to M^{\lambda\lambda'}_{\lambda,\lambda'}$ is a $\mathcal{W}$-module isomorphism. For the second assertion, let $\mathbf{G}_{\lambda,\lambda'} : V_\lambda \otimes V_{\lambda'} \to M^{\lambda\lambda'}_{\lambda,\lambda'} \otimes V_{\lambda \, \lambda'}$ be a $\mathcal{W}$-module isomorphism. Then $\mathbf{H}:=\mathbf{G}_{\lambda,\lambda'} \circ \mathbf{F}_{\lambda,\lambda'}^{-1} : M^{\lambda\lambda'}_{\lambda,\lambda'} \otimes V_{\lambda \, \lambda'} \to M^{\lambda\lambda'}_{\lambda,\lambda'} \otimes V_{\lambda \, \lambda'}$ is a $\mathcal{W}$-module isomorphism. Let us extend to a more general case for later use; let $M,M'$ be any trivial $\mathcal{W}$-representations (of any dimensions), $V_{\lambda_0}$ be any irreducible $\mathcal{W}$-representation where $\lambda_0=(x_0,y_0)$ is a non-singular weight, and let $\mathbf{H} : M \otimes V_{\lambda_0} \to M' \otimes V_{\lambda_0}$ be any $\mathcal{W}$-module homomorphism. We shall show that $\mathbf{H}$ acts only on the first tensor factor. The intertwining equation $\mathbf{H}\, (\mathrm{id}\otimes W) \, (f\otimes v) = (\mathrm{id}\otimes W) \, \mathbf{H} \, (f\otimes v)$ must hold for each $W\in \mathcal{W}$, $f\in M$, and $v\in V_{\lambda_0}$. Choose a basis $f_1,f_2,\ldots,f_n$ of the vector space $M$ and a basis $f'_1,\ldots,f'_{n'}$ of $M'$, and write $\mathbf{H}$ as $\mathbf{H}(f_i \otimes e_j) = \sum_{k,l} H^{k,l}_{i,j} \, f_k' \otimes e_l$ for numbers $H^{k,l}_{i,j}$. The intertwining equation for $W=X$ on $f_i \otimes e_j$ says $\sum_{k,l} x_0^{1/N} q^{2j} H^{k,l}_{i,j} f_k' \otimes e_l = \sum_{k,l} x_0^{1/N} q^{2l} H^{k,l}_{i,j} f_k' \otimes e_l$, leading to $q^{2j} H^{k,l}_{i,j} = q^{2l} H^{k,l}_{i,j}$ for all $i,j,k,l$, so we must have $H^{k,l}_{i,j} = 0$ unless $j=l$. The intertwining equation for $W=Y$ on $f_i \otimes e_j$ says $\sum_{k,l} y_0^{1/N} H^{k,l}_{i,j+1} f_k' \otimes e_l = \sum_{k,l} y_0^{1/N} H^{k,l}_{i,j} f_k' \otimes e_{l+1}$, leading to $H^{k,l}_{i,j+1} = H^{k,l-1}_{i,j}$ for all $i,j,k,l$. Hence $H^{k,j}_{i,j} = H^{k,j'}_{i,j'}$ for all $i,j,k,j'$; let us denote $H^k_i := H^{k,j}_{i,j}$. Thus, in the end, we have $\mathbf{H}(f_i \otimes e_j) = \sum_k H^k_i f_k' \otimes e_j = (\sum_k H^k_i f_k') \otimes e_j$, for each $i,j$. Define a linear map $\mathbf{U} : M \to M'$ by $U(f_i) := \sum_{k=1}^{n'} H^k_i f_k'$, $\forall i=1,\ldots,n$; then we showed $\mathbf{H} = \mathbf{U} \otimes \mathrm{id}$, as asserted. For our particular case when $\mathbf{H} = \mathbf{G}_{\lambda,\lambda'} \circ \mathbf{F}_{\lambda,\lambda'}^{-1}$, we see that $U$ is invertible and we have $\mathbf{G}_{\lambda,\lambda'} = (\mathbf{U}\otimes \mathrm{id})\circ \mathbf{F}_{\lambda,\lambda'}$ as desired. \qed

\subsection{The $6j$-symbol ${\bf T}$}
\label{subsec:T}

Let $(\lambda,\lambda',\lambda'')$ be a regular triple of weights. Then one can make sense of the $\mathcal{W}$ action on the tensor product of the corresponding representations either by writing $(V_\lambda \otimes V_{\lambda'})\otimes V_{\lambda''}$, or by $V_\lambda \otimes (V_{\lambda'}\otimes V_{\lambda''})$, using the coproduct of $\mathcal{W}$ (recall that we know how to define $\mathcal{W}$-action on a tensor product of any two $\mathcal{W}$-representations using the coproduct). These two actions are the same, because of the coassociativity of the coproduct; these two parenthesizings lead to two ways of decomposing the representation $V_\lambda \otimes V_{\lambda'} \otimes V_{\lambda''}$ into the direct sum of $V_{\lambda \lambda'\lambda''}$; we shall realize this as a decomposition map into $M\otimes M \otimes V_{\lambda \,\lambda'\,\lambda''}$. We make use of our decomposition map $\mathbf{F}$. Consider the diagram
\begin{align*}
\xymatrix@C+6mm{
(V_\lambda \otimes V_{\lambda'}) \otimes V_{\lambda''} \ar[r]^-{(\mathbf{F}_{\lambda,\lambda'})_{12}} \ar[d]^{\mathrm{id}}& M^{\lambda\lambda'}_{\lambda,\lambda'} \otimes V_{\lambda \lambda'} \otimes V_{\lambda''} \ar[r]^-{(\mathbf{F}_{\lambda\lambda', \lambda''})_{23}} & M^{\lambda\lambda'}_{\lambda,\lambda'} \otimes M^{\lambda\lambda'\lambda''}_{\lambda\lambda', \lambda''} \otimes V_{\lambda  \lambda'  \lambda''} \ar@{.>}[d] \\
V_\lambda \otimes (V_{\lambda'} \otimes V_{\lambda''}) \ar[r]^-{(\mathbf{F}_{\lambda',\lambda''})_{23}} & V_\lambda \otimes M^{\lambda'\lambda''}_{\lambda',\lambda''} \otimes V_{\lambda'\lambda''} \ar[r]^-{(\mathbf{F}_{\lambda,\lambda'\lambda''})_{13}} & M^{\lambda\lambda'\lambda''}_{\lambda,\lambda'\lambda''} \otimes M^{\lambda'\lambda''}_{\lambda',\lambda''} \otimes V_{\lambda\lambda'\lambda''},
}
\end{align*}
where $\mathbf{F}_{ij}$ means $\mathbf{F}$ is being applied to $i$-th and $j$-th tensor factors; for example, $(\mathbf{F}_{\lambda,\lambda'})_{12} = \mathbf{F}_{\lambda,\lambda'} \otimes \mathrm{id}$. Each of the five solid arrows above are $\mathcal{W}$-module isomorphisms, hence so is their composition, which is the dotted arrow:
$$
(\mathbf{F}_{\lambda,\lambda'\lambda''})_{13} (\mathbf{F}_{\lambda',\lambda''})_{23} (\mathbf{F}_{\lambda,\lambda'})_{12}^{-1} \, (\mathbf{F}_{\lambda\lambda',\lambda''})_{23}^{-1} : M^{\lambda\lambda'}_{\lambda,\lambda'} \otimes M^{\lambda\lambda'\lambda''}_{\lambda\lambda',\lambda''} \otimes V_{\lambda\lambda'\lambda''} \to M^{\lambda\lambda'\lambda''}_{\lambda,\lambda'\lambda''} \otimes M^{\lambda'\lambda''}_{\lambda',\lambda''} \otimes V_{\lambda\lambda'\lambda''}.
$$
Since this is a $\mathcal{W}$-module isomorphism, where $V_{\lambda\lambda'\lambda''}$ is irreducible and both $M^{\lambda\lambda'}_{\lambda,\lambda'} \otimes M^{\lambda\lambda'\lambda''}_{\lambda\lambda',\lambda''}$ and $M^{\lambda\lambda'\lambda''}_{\lambda,\lambda'\lambda''}\otimes M^{\lambda'\lambda''}_{\lambda',\lambda''}$ are trivial $\mathcal{W}$-representations, by a claim shown in the proof of Prop.\ref{prop:non-uniqueness_of_F} we deduce that this map acts only on the first two tensor factors as an invertible map, which we denote by ${\bf T} : M\otimes M \to M\otimes M$; for our later purposes, we define it after switching the two $M$ factors:
\begin{proposition}
For any regular triple $(\lambda,\lambda',\lambda'')$ of weights, there exists a unique invertible linear operator
$$
\mathbf{T}_{\lambda,\lambda',\lambda''} : M^{\lambda\lambda'\lambda''}_{\lambda\lambda',\lambda''} \otimes M^{\lambda\lambda'}_{\lambda,\lambda'} \to M^{\lambda'\lambda''}_{\lambda',\lambda''} \otimes M^{\lambda\lambda'\lambda''}_{\lambda,\lambda'\lambda''}
$$
such that the dotted arrow in the above diagram is
\begin{align}
\label{eq:F_and_T}
(\mathbf{F}_{\lambda,\lambda'\lambda''})_{13} (\mathbf{F}_{\lambda',\lambda''})_{23} (\mathbf{F}_{\lambda,\lambda'})_{12}^{-1} \, (\mathbf{F}_{\lambda\lambda',\lambda''})_{23}^{-1} = (\mathbf{T}_{\lambda,\lambda',\lambda''})_{21}. \qed
\end{align}
\end{proposition}
Note that $(\mathbf{T}_{\lambda,\lambda',\lambda''})_{21}$ can be understood as $\mathbf{P}_{(12)} \, (\mathbf{T}_{\lambda,\lambda',\lambda''})_{12} \, \mathbf{P}_{(12)}$, where $\mathbf{P}_{(ij)} = \mathbf{P}_{(ij)}^{-1}$ is the permutation map of the $i$-th and the $j$-th tensor factors. This map $\mathbf{T}_{\lambda,\lambda',\lambda''}$, which identifies the two realizations of the multiplicity space of $V_{\lambda\lambda'\lambda''}$ inside $V_\lambda\otimes V_{\lambda'}\otimes V_{\lambda''}$, is called the \ul{\em $6j$-symbol}. One can think of this $6j$-symbol map $\mathbf{T}_{\lambda,\lambda',\lambda''}$ as encoding the parenthesize-change $(V_\lambda\otimes V_{\lambda'})\otimes V_{\lambda''} \leadsto V_\lambda \otimes (V_{\lambda'} \otimes V_{\lambda''})$.

\vs

Now consider a regular quadruple of weights $(\lambda,\lambda',\lambda'',\lambda''')$, and the corresponding tensor product representation $V_\lambda \otimes V_{\lambda'} \otimes V_{\lambda''} \otimes V_{\lambda'''}$. There are five different ways of putting the parenthesizes, each leading to a particular way of decomposing this representation into $M\otimes M\otimes M\otimes V_{\lambda\lambda'\lambda''\lambda'''}$. In the diagram
\begin{align*}
\begin{array}{c}
((V_\lambda\otimes V_{\lambda'}) \otimes V_{\lambda''}) \otimes V_{\lambda'''} \to (V_\lambda \otimes (V_{\lambda'} \otimes V_{\lambda''})) \otimes V_{\lambda'''} \to V_\lambda \otimes (( V_{\lambda'} \otimes V_{\lambda''}) \otimes V_{\lambda'''}) \\
\searrow \qquad\qquad\qquad\qquad\qquad\qquad\qquad\qquad \swarrow \\
(V_\lambda \otimes V_{\lambda'}) \otimes (V_{\lambda''} \otimes V_{\lambda'''}) \to 
V_\lambda \otimes (V_{\lambda'} \otimes (V_{\lambda''} \otimes V_{\lambda'''})),
\end{array}
\end{align*}
all arrows are identity maps and are $\mathcal{W}$-module isomorphisms, due to coassociativity of the coproduct. For each of the five parenthesizing, the parentheses tell us how to use $\mathbf{F}$ maps to decompose into $M\otimes M \otimes M \otimes V_{\lambda\lambda'\lambda''\lambda'''}$, and one can rewrite the above diagram as five arrows among such spaces; each map is again a $\mathcal{W}$-module isomorphism, and hence by a claim in the proof of Prop.\ref{prop:non-uniqueness_of_F}, it acts only on the first three $M$ factors. So we obtain the following commutative diagram
$$
\hspace{-23mm}\begin{array}{c}
M^{\lambda\lambda'\lambda''\lambda'''}_{\lambda\lambda'\lambda'', \lambda'''} \otimes M^{\lambda\lambda'\lambda''}_{\lambda\lambda',\lambda''} \otimes M^{\lambda\lambda'}_{\lambda,\lambda'} \stackrel{(\mathbf{T}_{\lambda,\lambda',\lambda''})_{23}}{\to}
M^{\lambda\lambda'\lambda''\lambda'''}_{\lambda\lambda'\lambda'',\lambda'''}\otimes M^{\lambda'\lambda''}_{\lambda',\lambda''} \otimes M^{\lambda\lambda'\lambda''}_{\lambda,\lambda'\lambda''} \stackrel{(\mathbf{T}_{\lambda,\lambda'\lambda'',\lambda'''})_{13}}{\to}
M^{\lambda'\lambda''\lambda'''}_{\lambda'\lambda'',\lambda'''} \otimes M^{\lambda'\lambda''}_{\lambda',\lambda''} \otimes M^{\lambda\lambda'\lambda''\lambda'''}_{\lambda,\lambda'\lambda''\lambda'''} \\
(\mathbf{T}_{\lambda\lambda', \lambda'', \lambda'''})_{12} ~ \searrow \qquad\qquad\qquad\qquad\qquad\qquad\qquad\qquad\qquad\qquad\qquad\qquad \swarrow ~ (\mathbf{T}_{\lambda',\lambda'',\lambda'''})_{12} \\
M^{\lambda''\lambda'''}_{\lambda'',\lambda'''} \otimes M^{\lambda\lambda'\lambda''\lambda'''}_{\lambda\lambda',\lambda''\lambda'''} \otimes M^{\lambda\lambda'}_{\lambda,\lambda'} \stackrel{(\mathbf{T}_{\lambda,\lambda',\lambda''\lambda'''})_{23}}{\to}
M^{\lambda''\lambda'''}_{\lambda'',\lambda'''} \otimes M^{\lambda'\lambda''\lambda'''}_{\lambda',\lambda''\lambda'''} \otimes M^{\lambda\lambda'\lambda''\lambda'''}_{\lambda,\lambda'\lambda''\lambda'''}.
\end{array}
$$
\begin{proposition}[the pentagon relation of the $6j$-symbol]
\label{prop:TT_TTT}
For any regular quadruple of weights $(\lambda,\lambda',\lambda'',\lambda''')$, one has
$$
(\mathbf{T}_{\lambda,\lambda',\lambda''\lambda'''})_{23} \, (\mathbf{T}_{\lambda\lambda',\lambda'',\lambda'''})_{12}
= (\mathbf{T}_{\lambda',\lambda'',\lambda'''})_{12} \, (\mathbf{T}_{\lambda,\lambda'\lambda'',\lambda'''})_{13} \, (\mathbf{T}_{\lambda,\lambda',\lambda''})_{23}. \qed
$$
\end{proposition}
Notice that, whatever constants $\zeta_{{\bf a}, {\bf c}}$ we choose for the cyclic quantum dilogarithms, the above pentagon relation of the ${\bf T}$ maps hold without any constant.

\vs

As the construction of the $6j$-symbol map $\mathbf{T}_{\lambda,\lambda',\lambda''}$ depends on a non-unique choice of the decomposition map $\mathbf{F}_{\lambda,\lambda'}$, one might worry about the uniqueness of $\mathbf{T}_{\lambda,\lambda',\lambda''}$. Note from Prop.\ref{prop:non-uniqueness_of_F} that we could have used a different choice of the decomposition map $V_\lambda \otimes V_{\lambda'} \to M^{\lambda\lambda'}_{\lambda,\lambda'} \otimes V_{\lambda,\lambda'}$ for a regular pair $(\lambda,\lambda')$ of weights other than our choice $\mathbf{F}_{\lambda,\lambda'}$, say $\til{\mathbf{F}}_{\lambda,\lambda'}$, which would be of the form $(\mathbf{U}\otimes \mathrm{id}) \circ \mathbf{F}_{\lambda,\lambda'}$, or $(\til{\mathbf{F}}_{\lambda,\lambda'})_{12} = \mathbf{U}_1 \, (\mathbf{F}_{\lambda,\lambda'})_{12}$, for some invertible linear map $\mathbf{U}$ on $M^{\lambda\lambda'}_{\lambda,\lambda'} \cong \mathbb{C}^N$. Then the corresponding new $6j$-symbol $\til{\mathbf{T}}_{\lambda,\lambda',\lambda''}$ would be
\begin{align*}
(\til{\mathbf{T}}_{\lambda,\lambda',\lambda''})_{21} & = (\til{\mathbf{F}}_{\lambda,\lambda'\lambda''})_{13} \, (\til{\mathbf{F}}_{\lambda',\lambda''})_{23} \, (\til{\mathbf{F}}_{\lambda,\lambda'})^{-1}_{12} \, (\til{\mathbf{F}}_{\lambda\lambda',\lambda''})^{-1}_{23} \\
& = \mathbf{U}_1 \, (\mathbf{F}_{\lambda,\lambda'\lambda''})_{13} \, \mathbf{U}_2 \, (\mathbf{F}_{\lambda',\lambda''})_{23} \, (\mathbf{F}_{\lambda,\lambda'})_{12}^{-1} \, \mathbf{U}_1^{-1} \, (\mathbf{F}_{\lambda\lambda',\lambda''})_{23}^{-1} \, \mathbf{U}_2^{-1} \\
& = \mathbf{U}_1 \mathbf{U}_2 \, (\mathbf{F}_{\lambda,\lambda'\lambda''})_{13} \, (\mathbf{F}_{\lambda',\lambda''})_{23} \, (\mathbf{F}_{\lambda,\lambda'})_{12}^{-1} \, (\mathbf{F}_{\lambda\lambda',\lambda''})_{23}^{-1} \, (\mathbf{U}_1\mathbf{U}_2)^{-1} \\
& = (\mathbf{U}_1\mathbf{U}_2) \, (\mathbf{T}_{\lambda,\lambda',\lambda''})_{12} \, (\mathbf{U}_1\mathbf{U}_2)^{-1},
\end{align*}
hence $\til{\mathbf{T}}_{\lambda,\lambda',\lambda''} = (\mathbf{U}\otimes \mathbf{U})\circ \mathbf{T}_{\lambda,\lambda',\lambda''} \circ (\mathbf{U}\otimes \mathbf{U})^{-1}$. So we have a control over the (non-)uniqueness of $\mathbf{T}_{\lambda,\lambda',\lambda''}$.

\section{The order three operator}

Results thus far are stated clearly already in the literature, e.g. in \cite{K94} \cite{K99}, and \cite{Bai}, or references therein. The aspects of the representation theory of the quantum torus $\mathcal{W}$ studied in the present section are not as clearly stated in the literature in the style of the present paper, although most of the ingredients appear in \cite{GKT} in a  different language. I find the formulation in \cite{GKT} somewhat enigmatic and difficult to decode. The formulation of the present paper uses only plain representation theory language, thus is more accessible.

\subsection{Dual and Hom representations}

\begin{definition}[Hom representations; modified from \cite{CP}]
\label{def:Hom_representations}
For any $\mathcal{W}$-representations $V,U$, define the \ul{\em left-Hom} and \ul{\em right-Hom} representations $\mathrm{Hom}^\mathrm{L}_\mathbb{C}(V,U)$ and $\mathrm{Hom}^\mathrm{R}_\mathbb{C}(V,U)$ to be the vector space $\mathrm{Hom}_\mathbb{C}(V,U)$ with the $\mathcal{W}$-actions defined respectively by
\begin{align*}
(W.f)(v) & := \sum W_{(1)}.(f(S(W_{(2)}).v)), \qquad \mbox{for all}\quad W\in \mathcal{W}, ~ f\in \mathrm{Hom}^\mathrm{L}_\mathbb{C}(V,U), ~ v\in V, \\
(W.f)(v) & := \sum W_{(2)}.(f(S^{-1}(W_{(1)}).v)), \quad \mbox{for all}\quad W\in \mathcal{W}, ~ f\in \mathrm{Hom}^\mathrm{R}_\mathbb{C}(V,U), ~ v\in V,
\end{align*}
where $\Delta W = \sum W_{(1)}\otimes W_{(2)}$.
\end{definition}
\begin{definition}[dual representations; \cite{CP}]
In the above definition, when $U$ is the one-dimensional trivial $\mathcal{W}$-representation $\mathbb{C}$, we call the corresponding Hom-representations $\mathrm{Hom}^\mathrm{L}_\mathbb{C}(V,\mathbb{C})$ and $\mathrm{Hom}^\mathrm{R}_\mathbb{C}(V,\mathbb{C})$ are denoted by $V^*$ and ${}^*V$ and called the \ul{\em left dual} and the \ul{\em right dual} representations of $V$, respectively.
\end{definition}
\begin{remark}
The terms `left-Hom' and `right-Hom' are coined here as such, inspired by the names `left dual' and `right dual' which are used e.g. in \cite{CP}.
\end{remark}
Let us check that the above defined Hom and dual representations are well-defined. First, the actions on the dual representations can be simplified to:
\begin{align}
\label{eq:explicit_dual_actions}
\left\{ {\renewcommand{\arraystretch}{1.2} \begin{array}{ll}
(W.f)(v) = f(S(W).v), & \forall W\in \mathcal{W}, ~ f\in V^*, ~ v\in V,  \\
(W.f)(v) = f(S^{-1}(W).v), & \forall W\in \mathcal{W}, ~ f\in {}^*V, ~ v\in V,
\end{array}} \right.
\end{align}
using the facts $\sum \epsilon(W_{(1)}) \, S(W_{(2)}) = S(W)$ and $\sum S^{-1}(W_{(1)}) \, \epsilon(W_{(2)}) = S^{-1}(W)$, which follows from Hopf algebra axioms, or can just be checked in our case for $W=X$ and $Y$. It is easy to see that each of these two indeed defines a representation of the algebra $\mathcal{W}$ on the vector space $\mathrm{Hom}_\mathbb{C}(V,\mathbb{C})$, using the fact that $S:\mathcal{W} \to \mathcal{W}$ is a anti-homomorphism of rings; for the left dual, $(W_1.(W_2.f))(v) = (W_2.f)(S(W_1).v) = f(S(W_2).(S(W_1).v)) = f((S(W_2)S(W_1)).v)=f(S(W_1W_2).v) = ((W_1W_2).f)(v)$, and similarly for the right dual.

\vs

Now, consider the {\em canonical} vector space isomorphisms
\begin{align}
\label{eq:Hom_as_tensor_products}
J_{V,U}^\mathrm{L} : \mathrm{Hom}^\mathrm{L}_\mathbb{C}(V,U) \stackrel{\sim}{\to} U \otimes V^* \qquad \mbox{and} \qquad J_{V,U}^\mathrm{R}: \mathrm{Hom}^\mathrm{R}_\mathbb{C}(V,U) \stackrel{\sim}{\to} {}^* V \otimes U.
\end{align}
Of course, as vector spaces, the superscripts L and R do not matter, $V^*$ and ${}^*V$ are the same, and tensor products are commutative. The isomorphism sends each element $\sum_i u_i \otimes f_i$ (a finite sum) of $U\otimes V^*$, equivalently an element $\sum_i f_i \otimes u_i$ of ${}^* V\otimes U$, to the element of $\mathrm{Hom}_\mathbb{C}(V,U)$ that sends each $v\in V$ to $\sum_i u_i \cdot f_i(v)$. 
\begin{lemma}[composition of canonical maps]
\label{lem:composition_of_canonical_maps}
The map $J^\mathrm{R}_{V,U} \circ (J^\mathrm{L}_{V,U})^{-1} : U \otimes V^* \to {}^* V \otimes U$ equals the factor permuting map $\mathbf{P}_{(12)}$ sending each $u\otimes \xi$ to $\xi \otimes u$, if $V^*$ and ${}^*V$ are identified naturally as vector spaces. \qed
\end{lemma}
Now, one can view $U\otimes V^*$ and ${}^*V \otimes U$ as $\mathcal{W}$-representations using the coproduct of $\mathcal{W}$; in particular, the order of the tensor product does matter. Then, in view of Def.\ref{def:Hom_representations}, one observes that the isomorphisms \eqref{eq:Hom_as_tensor_products} are isomorphisms of $\mathcal{W}$-representations. This shows that Def.\ref{def:Hom_representations} indeed defines well-defined $\mathcal{W}$-representations.

\begin{definition}[\cite{CP}]
For a $\mathcal{W}$-representation $V$, the \ul{\em invariant subspace} of $V$ is defined to be the set of all elements of $V$ on which the Hopf algebra $\mathcal{W}$ acts trivially, i.e. by the counit:
$$
\mathrm{Inv}(V) = \{ v\in V : W.v = \epsilon(W)\cdot v, ~ \forall W\in\mathcal{W} \}.
$$
\end{definition}

\begin{lemma}[space of intertwiners as invariant subspace]
\label{lem:space_of_intertwiners_as_invariant_subspace}
For $\mathcal{W}$-representations $V,U$, both spaces $\mathrm{Inv}(\mathrm{Hom}^\mathrm{L}_\mathbb{C}(V,U))$ and $\mathrm{Inv}(\mathrm{Hom}^\mathrm{R}_\mathbb{C}(V,U))$ coincide with the subspace $\mathrm{Hom}_\mathcal{W}(V,U) := \{ f\in \mathrm{Hom}_\mathbb{C}(V,U) : f(W.v) = W.(f(v)), ~ \forall W\in \mathcal{W}, ~ \forall v\in V \}$ of $\mathrm{Hom}_\mathbb{C}(V,U)$.
\end{lemma}
{\it Proof.} We first show that $\mathrm{Hom}_\mathcal{W}(V,U)$ is included in both invariant subspaces. Let $f\in \mathrm{Hom}_\mathcal{W}(V,U) \subset \mathrm{Hom}_\mathbb{C}(V,U)$. Then, for each $W\in \mathcal{W}$ and $v\in V$ one has $\sum W_{(1)}.(f(S(W_{(2)}).v)) = \sum W_{(1)}.(S(W_{(2)}).f(v)) = (\sum W_{(1)} \, S(W_{(2)})). (f(v)) = \epsilon(W) \cdot f(v)$, thus the left-Hom action of $W$ on $f$ is trivial, hence $f\in \mathrm{Inv}(\mathrm{Hom}^\mathrm{L}_\mathbb{C}(V,U))$. Likewise, observe for each $W\in \mathcal{W}$ and $v\in V$ that $\sum W_{(2)}.(f(S^{-1}(W_{(1)}.v))) = \sum W_{(2)}.(S^{-1}(W_{(1)}).(f(v))) = (\sum W_{(2)} \, S^{-1}(W_{(1)})).(f(v)) = S^{-1}(\sum W_{(1)} S(W_{(2)})).(f(v)) = S^{-1}(\epsilon(W)).(f(v)) = \epsilon(W)\cdot f(v)$, thus the right-Hom action of $W$ on $f$ is trivial, so $f\in \mathrm{Inv}(\mathrm{Hom}^\mathrm{R}_\mathbb{C}(V,U))$ as asserted.

\vs

The reverse inclusion $\mathrm{Inv}(\mathrm{Hom}^\mathrm{L}_\mathbb{C}(V,U)) \subset \mathrm{Hom}_\mathcal{W}(V,U)$ is a standard exercise, using Hopf algebra axioms, and one can adapt such a proof to show $\mathrm{Inv}(\mathrm{Hom}^\mathrm{R}_\mathbb{C}(V,U)) \subset \mathrm{Hom}_\mathcal{W}(V,U)$. For $f\in \mathrm{Inv}(\mathrm{Hom}^\mathrm{R}_\mathbb{C}(V,U))$, $v\in V$, and $W\in \mathcal{W}$ with the Sweedler notations $\Delta W = \sum W_{(1)} \otimes W_{(2)}$ and $((\Delta\otimes 1) \Delta) W = \sum W_{(1)}\otimes W_{(2)}\otimes W_{(3)}$, note
\begin{align*}
f(W.v) & \textstyle = f((\sum W_{(1)} \, \epsilon(W_{(2)})). v) = \sum \epsilon(W_{(2)}) \cdot f(W_{(1)}.v)
= \sum (W_{(2)}.f)(W_{(1)}.v) \\
& \textstyle = \sum W_{(3)}.(f(S^{-1}(W_{(2)}).(W_{(1)}.v))) 
= \sum W_{(3)}.(f( S^{-1}( S(W_{(1)}) \, W_{(2)}). v)) \\
& \textstyle = \sum W_{(2)}.(f(S^{-1}(\epsilon(W_{(1)})).v)) = \sum W_{(2)}.(f(\epsilon(W_{(1)}) \cdot v)) \\
& \textstyle = \sum (\epsilon(W_{(1)}) \cdot W_{(2)}).(f(v)) = W.(f(v)).
\end{align*}
For those who are not familiar with Sweedler notations, we give explicit computation here for our particular Hopf algebra $\mathcal{W}$. Let $f\in \mathrm{Inv}(\mathrm{Hom}^\mathrm{L}_\mathbb{C}(V,U))$. The triviality of the left-Hom action says $\sum W_{(1)}.(f(S(W_{(2)}).v)) = \epsilon(W)\cdot f(v)$ for all $W\in \mathcal{W}$ and $v\in V$. Putting $W=X^{\pm 1}$ gives $X.(f(X^{-1}.v)) = f(v) = X^{-1}.(f(X.v))$, thus $f(X^{-1}.v) = X^{-1}.(f(v))$ and $f(X.v) = X.(f(v))$. Putting $W=Y$ gives $0 = Y.(f(X^{-1}.v)) + f((-YX^{-1}).v) = Y.(f(X^{-1}.v)) - f (Y.(X^{-1}.v))$; put $v = X.v'$ for arbitrary $v'\in V$ to get $f(Y.v')=Y.(f(v'))$. We thus proved $f\in \mathrm{Hom}_\mathcal{W}(V,U)$, so $\mathrm{Inv}(\mathrm{Hom}^\mathrm{L}_\mathbb{C}(V,U)) \subseteq \mathrm{Hom}_\mathcal{W}(V,U)$, leading to the equality.

\vs

Likewise, let $f\in \mathrm{Inv}(\mathrm{Hom}^\mathrm{R}_\mathbb{C}(V,U))$. Similarly as above, the triviality of the right-Hom action of $W=X^{\pm 1}$ gives $f(X^{\pm 1}.v) = X^{\pm 1}.(f(v))$ for all $v\in V$, while that for $W=Y$ yields $0 = X.(f(S^{-1}(Y).v)) + Y.(f(S^{-1}(1).v)) = X.(f( (-X^{-1} Y).v)) + Y.(f(v)) = - X.(f(X^{-1}.(Y.v))) + Y.(f(v)) = -X.(X^{-1}.(f(Y.v))) + Y.(f(v)) = - (X X^{-1}).(f(Y.v)) + Y.(f(v)) = - f(Y.v) + Y.(f(v))$, so $f(Y.v) = Y.(f(v))$ for all $v\in V$. Here we used $S^{-1}(X) = X^{-1}$ and $S^{-1}(Y) = -X^{-1}Y$ which can easily be verified. Thus $f\in \mathrm{Hom}_\mathcal{W}(V,U)$, so $\mathrm{Inv}(\mathrm{Hom}^\mathrm{R}_\mathbb{C}(V,U)) \subset \mathrm{Hom}_\mathcal{W}(V,U)$, leading to the equality. \qed

\vs

The left-Hom representation is standard. Note that the right-Hom representation is different from what is mentioned in the remark in \S4.1.C of \cite{CP}; what is there is the representation isomorphic to $U \otimes {}^*V$. The above Lem.\ref{lem:space_of_intertwiners_as_invariant_subspace} somewhat justifies why it is reasonable to consider the left-Hom and right-Hom representations as defined in Def.\ref{def:Hom_representations}. I claim without proof that for a general Hopf algebra, say $\mathcal{A}$, if we define an $\mathcal{A}$-representation on the vector space $\mathrm{Hom}_\mathbb{C}(V,U)$ using the realization $U\otimes {}^*V$ or $V^*\otimes U$ instead of $U\otimes V^*$ or ${}^* V\otimes U$ as we did \eqref{eq:Hom_as_tensor_products}, then the corresponding invariant subspace does not equal the space of intertwiners $\mathrm{Hom}_\mathcal{A}(V,U)$.

\subsection{Isomorphisms to dual representations}
\label{subsec:isomorphisms_to_dual_representations}

For a non-singular weight $\lambda = (x,y)$, let us realize the left and the right dual representations $V_\lambda^*$ and ${}^*V_\lambda$ of the cyclic irreducible representation $V_\lambda$ more explcitly on the vector space $\mathbb{C}^N$ with the basis $e_0,e_1,\ldots,e_{N-1}$, as we did for the cyclic irreducibles, instead of the dual vector space $\mathrm{Hom}_\mathbb{C}(V_\lambda,\mathbb{C}) \equiv \mathrm{Hom}_\mathbb{C}(\mathbb{C}^N, \mathbb{C})$. Here we view $\mathbb{C}^N$ as the dual of $\mathbb{C}^N$, using the bilinear pairing
\begin{align}
\label{eq:canonical_pairing_on_C_N}
\mathbb{C}^N \times \mathbb{C}^N \to \mathbb{C}, \quad (e_i, e_j) \mapsto \langle e_i, e_j \rangle = \delta_{i,j}, \quad \forall i,j,
\end{align}
where $\delta_{\cdot,\cdot}$ is the Kronecker delta. For $W \in \mathcal{W}$ denote by $\mu_\lambda^*(W)$ and ${}^*\mu_\lambda(W)$ the $\mathcal{W}$-actions on $V^*_\lambda \equiv \mathbb{C}$ and ${}^*V_\lambda \equiv \mathbb{C}$ respectively. Then
\begin{align*}
\langle \mu^*_\lambda(X).e_i, e_j\rangle & \stackrel{\eqref{eq:explicit_dual_actions},\eqref{eq:W_Hopf_algebra_structure}}{=} \langle e_i, \mu_\lambda(X^{-1}).e_j\rangle = \langle e_i, x^{-1/N} q^{-2j} e_j\rangle  = x^{-1/N} q^{-2j} \delta_{i,j} \\
& \,\,\,\,\,\, = x^{-1/N} q^{-2i} \delta_{i,j} = \langle x^{-1/N} q^{-2i} e_i, e_j\rangle \stackrel{\eqref{eq:A_B}}{=} \langle x^{-1/N} A^{-1} e_i, e_j\rangle, \\
\langle \mu^*_\lambda(Y).e_i,e_j\rangle & \stackrel{\eqref{eq:explicit_dual_actions},\eqref{eq:W_Hopf_algebra_structure}}{=} \langle e_i, \mu_\lambda (-YX^{-1}).e_j \rangle = \langle e_i, - y^{1/N} x^{-1/N} q^{-2j} e_{j+1} \rangle  = - y^{1/N} x^{-1/N} q^{-2j} \delta_{i,j+1} \\
& \,\,\,\,\,\, = - y^{1/N} x^{-1/N} q^{-2(i-1)} \delta_{i-1,j}
= \langle - y^{1/N}x^{-1/N} q^{-2(i-1)} e_{i-1}, e_j\rangle \\
& \,\,\,\,\,\, \stackrel{\eqref{eq:A_B}}{=} \langle - y^{1/N}x^{-1/N} A^{-1} B^{-1} e_i, e_j\rangle,
\end{align*}
hence the left dual action is given by
\begin{align}
\label{eq:left_dual_of_cyclic_irreducible}
\mu^*_\lambda(X) = x^{-1/N} \, A^{-1}, \qquad
\mu^*_\lambda(Y) = - y^{1/N} x^{-1/N} \, A^{-1} B^{-1}.
\end{align}
Similar computation shows that the right dual action is given by
\begin{align}
\label{eq:right_dual_of_cyclic_irreducible}
{}^*\mu_\lambda(X) = x^{-1/N} A^{-1}, \qquad
{}^*\mu_\lambda(Y) = - y^{1/N} x^{-1/N} \, B^{-1} A^{-1}.
\end{align}
In particular, $\mu^*_\lambda$ and ${}^*\mu_\lambda$ are cyclic, and hence also irreducible due to their dimensions and Prop.\ref{prop:cyclic_representations}. To find out their corresponding weights, observe
$$
(\mu^*_\lambda(X))^N = x^{-1} \cdot \mathrm{id} = ({}^*\mu_\lambda(X))^N, \qquad
(\mu^*_\lambda(Y))^N = -y x^{-1} \cdot \mathrm{id} = ({}^*\mu_\lambda(Y))^N.
$$
So, both $V_\lambda^*$ and ${}^*V_\lambda$ are isomorphic to the cyclic irreducible representation $V_{(x^{-1}, - yx^{-1})}$.

\begin{definition}
\label{def:dual_weight}
For a non-singular weight $\lambda = (x,y)$, define its \ul{\em dual weight} as
$$
\lambda^* := (x^{-1}, -yx^{-1}).
$$
\end{definition}

Easy observation:
\begin{lemma}
\label{lem:dual_weight}
For any non-singular weight $\lambda$, its dual weight $\lambda^*$ is non-singular, and one has
$$
(\lambda^*)^* = \lambda, \qquad \lambda \, \lambda^* = \mathbf{1} = \lambda^* \, \lambda,
$$
where $\mathbf{1} = (1,0)$. \qed
\end{lemma}
So $\lambda^*$ is the (unique) inverse of $\lambda$, with respect to the multiplication of weights.
\begin{corollary}
\label{cor:dual_weight_of_product}
For any non-singular weights $\lambda,\lambda'$, one has $(\lambda\lambda')^* = (\lambda')^* \lambda^*$. \qed
\end{corollary}

Let us construct explicit isomorphisms to duals:
\begin{proposition}[isomorphisms to duals]
\label{prop:isomorphisms_to_duals}
For a non-singular weight $\lambda=(x,y)$, the linear maps
\begin{align*}
& C_\lambda : V_\lambda^* \to V_{\lambda^*}, \qquad\quad e_i \mapsto q^{-2i\cdot m_{y,x^{-1}}}\, q^{-i(i-1)} \, e_{-i}, \quad \forall i, \\
& D_\lambda : V_{\lambda^*} \to {}^*V_\lambda, \qquad \,\,\,\,
e_i \mapsto q^{-2i\cdot m_{y,x^{-1}}} \, q^{i(i-1)} \, e_{-i}, \quad \forall i,
\end{align*}
are isomorphisms of $\mathcal{W}$-representations.  

\vs

In particular, the left dual and the right dual of a cyclic irreducible $\mathcal{W}$-representation are also cyclic irreducible and isomorphic to each other. \qed
\end{proposition}

Due to the irreduciblity, the maps $C_\lambda$ and $D_\lambda$ are unique $\mathcal{W}$-module isomorphisms up to scalar (see Cor.\ref{cor:conjugation_on_A_B_determines}). For a proof of the above proposition and for later use, we compute:
\begin{lemma}[conjugation action of $C_\lambda,D_\lambda$ on $A,B$]
\label{lem:conjugation_of_C_D}
For a non-singular weight $\lambda = (x,y)$, one has the following equalities of operators on $\mathbb{C}^N$:
\begin{align*}
& C_\lambda^{-1} \, A \, C_\lambda = A^{-1}, \qquad
C_\lambda^{-1} \, B \, C_\lambda = q^{-2m_{y,x^{-1}}} A^{-1} B^{-1}, \\
& D_\lambda^{-1} \, A \, D_\lambda = A^{-1}, \qquad
D_\lambda^{-1} \, B \, D_\lambda = q^{-2m_{y,x^{-1}} } A B^{-1}.
\end{align*}
\end{lemma}
{\it Proof.} For convenience, write $m = m_{y,x^{-1}}$. It is easy to see that $C_\lambda$ and $D_\lambda$ are invertible, with inverses given by $C_\lambda^{-1}(e_i) = q^{-2im} q^{i(i+1)} e_{-i}$ and $D_\lambda^{-1}(e_i) = q^{-2im} q^{-i(i+1)} e_{-i}$. On the basis vector $e_i$,
\begin{align*}
C_\lambda^{-1} A C_\lambda e_i & = C_\lambda^{-1} A q^{-2im} q^{-i(i-1)} e_{-i}
= C_\lambda^{-1} q^{-2im} q^{-i(i-1)} q^{-2i} e_{-i}
= q^{-2i} C_\lambda^{-1} C_\lambda e_i = q^{-2i} e_i = A^{-1} e_i, \\
C_\lambda^{-1} B C_\lambda e_i & = C_\lambda^{-1} B q^{-2im} q^{-i(i-1)} e_{-i}
= C_\lambda^{-1} q^{-2im} q^{-i(i-1)} e_{-i+1} \\
& = q^{-2im} q^{-i(i-1)} q^{-2(-i+1)m} q^{(-i+1)(-i+2)} e_{i-1}
= q^{-2m} q^{-2(i-1)} e_{i-1}
= q^{-2m} A^{-1} B^{-1} e_i, \\
D_\lambda^{-1} A D_\lambda e_i & = D_\lambda^{-1} A q^{-2im} q^{i(i-1)} e_{-i} = D_\lambda^{-1} q^{-2im} q^{i(i-1)} q^{-2i} e_{-i} = q^{-2i} D_\lambda^{-1} D_\lambda e_i = A^{-1} e_i, \\
D_\lambda^{-1} B D_\lambda e_i & = D_\lambda^{-1} B q^{-2im} q^{i(i-1)} e_{-i} = D_\lambda^{-1} q^{-2im} q^{i(i-1)} e_{-i+1} \\
& = q^{-2im} q^{i(i-1)} q^{-2(-i+1)m} q^{-(-i+1)(-i+2)} e_{i-1} = q^{-2m} q^{2(i-1)} e_{i-1} = q^{-2m} A B^{-1} e_i. \qed
\end{align*}

\vs

{\it Proof of Prop.\ref{prop:isomorphisms_to_duals}.} Note first from \eqref{eq:normalization_for_N_th_roots} that $(-yx^{-1})^{1/N} = - (yx^{-1})^{1/N} = - q^{2m_{y,x^{-1}}} \, y^{1/N} x^{-1/N}$. For convenience, denote $m = m_{y,x^{-1}}$. Observe
\begin{align*}
& C_\lambda^{-1} \, \mu_{\lambda^*}(X) = C_\lambda^{-1} x^{-1/N} A \,\stackrel{{\rm Lem}.\ref{lem:conjugation_of_C_D}}{=}\, x^{-1/N} A^{-1} C_\lambda^{-1} \, \stackrel{\eqref{eq:left_dual_of_cyclic_irreducible}}{=} \, \mu_\lambda^*(X) \, C_\lambda^{-1}, \\
& C_\lambda^{-1} \, \mu_{\lambda^*}(Y) = C_\lambda^{-1} (-yx^{-1})^{1/N} B \,\stackrel{{\rm Lem}.\ref{lem:conjugation_of_C_D}}{=}\, - y^{1/N} x^{-1/N} A^{-1} B^{-1} C_\lambda^{-1} \, \stackrel{\eqref{eq:left_dual_of_cyclic_irreducible}}{=} \, \mu^*_\lambda(Y) \, C_\lambda^{-1}, \\
& D_\lambda^{-1} \, {}^* \mu_\lambda(X) \, \stackrel{\eqref{eq:right_dual_of_cyclic_irreducible}}{=} \, D_\lambda^{-1} x^{-1/N} A^{-1} \, \stackrel{{\rm Lem}.\ref{lem:conjugation_of_C_D}}{=} \, x^{-1/N} A D_\lambda^{-1} = \mu_{\lambda^*}(X) \, D_\lambda^{-1}, \\
& D_\lambda^{-1} \, {}^* \mu_\lambda(Y) \, \stackrel{\eqref{eq:right_dual_of_cyclic_irreducible}}{=} \, D_\lambda^{-1} (-y^{1/N} x^{-1/N} B^{-1} A^{-1}) \, \stackrel{{\rm Lem}.\ref{lem:conjugation_of_C_D}}{=} \, -y^{1/N} x^{-1/N} q^{2m} B A^{-1} A D_\lambda^{-1} \\
& \qquad \qquad \quad  = (-yx^{-1})^{1/N} B D_\lambda^{-1} = \mu_{\lambda^*}(Y) D_\lambda^{-1},
\end{align*}
so $C_\lambda^{-1} \, \mu_{\lambda^*}(W) = \mu^*_\lambda(W) C_\lambda^{-1}$ and $D_\lambda^{-1} \, {}^*\mu_\lambda(W) = \mu_{\lambda^*}(W) D_\lambda^{-1}$ for $W=X,Y$, so $C_\lambda^{-1} : V_{\lambda^*} \to V_\lambda^*$ and $D_\lambda^{-1} : {}^* V_\lambda \to V_{\lambda^*}$ are indeed $\mathcal{W}$-module maps. \qed

\begin{corollary}
For a non-singular weight $\lambda$, we have isomorphisms of $\mathcal{W}$-representations
$$
(V_\lambda^*)^* \cong {}^*(V_\lambda^*) \cong ({}^*V_\lambda)^* \cong {}^*({}^*V_\lambda) \cong V_\lambda. \qed
$$
\end{corollary}

\subsection{Some lemmas on intertwiner spaces}

\begin{lemma}[invariant subspace and trivial representations]
\label{lem:inv_and_triv}
Let $V$ be a $\mathcal{W}$-representation and $M$ be a trivial $\mathcal{W}$-representation. We have the equality of vector subspaces
$$
M \otimes \mathrm{Inv}(V) = \mathrm{Inv}(M\otimes V)
$$
of $M\otimes V$.
\end{lemma}
{\it Proof.} Each element $W\in \mathcal{W}$ acts on $M\otimes V$ as $(\mu_M \otimes \mu_V)(\sum W_{(1)} \otimes W_{(2)}) = \sum \epsilon(W_{(1)}) \otimes \mu_V(W_{(2)}) = \mathrm{id} \otimes \mu_V( \sum \epsilon(W_{(1)}) \, W_{(2)}) = \mathrm{id} \otimes \mu_V(W)$, where $\Delta W = \sum W_{(1)} \otimes W_{(2)}$. Let $f_1,\ldots,f_n$ be a basis of $M$. Then each element of $M\otimes V$ can be written uniquely as $\sum_{i=1}^n f_i \otimes v_i$ for some vectors $v_i \in V$. As just seen, this element lies in $\mathrm{Inv}(M\otimes V)$ if and only if $\sum_i f_i \otimes (W.v_i) = \epsilon(W)\cdot \sum_i f_i \otimes v_i$, which holds if and only if $W.v_i = \epsilon(W)\cdot v_i$ for each $i$, which is equivalent to the condition $\sum_i f_i \otimes v_i \in M\otimes \mathrm{Inv}(V)$. \qed

\begin{lemma}[taking out trivial representation from Hom]
\label{lem:taking_out_triv}
The linear map
$$
M \otimes \mathrm{Hom}_\mathcal{W}(V,U) \to \mathrm{Hom}_\mathcal{W}(V, M \otimes U)
$$
sending $f \otimes g$, for each $f\in M$ and $g\in \mathrm{Hom}_\mathcal{W}(V,U)$, to the element of $V \to M\otimes U$ defined by $v\mapsto f\otimes (g(v))$, $\forall v\in V$, is an isomorphism.
\end{lemma}

{\it Proof.} We have a following chain of identifications and isomorphisms:
\begin{align*}
M\otimes \mathrm{Hom}_\mathcal{W}(V,U) & \stackrel{\tiny \circled{1}}{=} M\otimes \mathrm{Inv}(\mathrm{Hom}_\mathbb{C}^\mathrm{L}(V,U)) \stackrel{\tiny \circled{2}}{\to} M \otimes \mathrm{Inv}(U \otimes V^*) \stackrel{\tiny \circled{3}}{=} \mathrm{Inv}(M\otimes U \otimes V^*) \\
& \stackrel{\tiny \circled{4}}{\to} \mathrm{Inv}( \mathrm{Hom}^\mathrm{L}_\mathbb{C}(V, M\otimes U) ) \stackrel{\tiny \circled{5}}{=} \mathrm{Hom}_\mathcal{W}(V,M\otimes U).
\end{align*}
The equalities \circled{1} and \circled{5} are from Lem.\ref{lem:space_of_intertwiners_as_invariant_subspace}, the maps \circled{2} and \circled{4} are obtained by the restricting to the invariant subspaces of the (left) canonical map in \eqref{eq:Hom_as_tensor_products}, and the equality \circled{3} is from Lem.\ref{lem:inv_and_triv}. From the defintion of the (left) canonical map \eqref{eq:Hom_as_tensor_products}, it is straightforward to verify that the composition of the above five equalities and arrows is the asserted map. \qed

\begin{corollary}[multiplicity space as Hom]
\label{cor:M_as_Hom}
If $V_1,V_2$ are $\mathcal{W}$-representations, $V_1$ is irreducible, $M$ is a trivial $\mathcal{W}$-representation, and $F : V_2 \cong M\otimes V_1$ is an isomorphism of $\mathcal{W}$-representations, then the linear map
\begin{align}
\label{eq:M_as_Hom}
M \to \mathrm{Hom}_\mathcal{W}(V_1, V_2)
\end{align}
sending each $f\in M$ to the element $V_1 \to V_2$ defined by $v \mapsto F^{-1}(f\otimes v)$ is an isomorphism.
\end{corollary}

{\it Proof.} Consider the chain of isomorphisms
$$
\mathrm{Hom}_\mathcal{W}(V_1,V_2) \stackrel{\tiny \circled{1}}{\to} \mathrm{Hom}_\mathcal{W}(V_1, M\otimes V_1) \stackrel{\tiny \circled{2}}{\to} M \otimes \mathrm{Hom}_\mathcal{W}(V_1,V_1) \stackrel{\tiny \circled{3}}{\to} M.
$$
The map \circled{1} is induced by $F$, \circled{2} is by Lem.\ref{lem:taking_out_triv}, and \circled{3} sends each $f\otimes \mathrm{id}$ to $f$; we used irreducibility of $V_1$ to deduce that $\mathrm{Hom}_\mathcal{W}(V_1,V_1)$ is one-dimensional, generated by $\mathrm{id}$. It is straightforward to verify that the inverse of the composition of these three isomosphisms is the asserted map. \qed

\begin{lemma}[composition lemma]
\label{lem:composition_lemma}
Let $V_1,V_2,V_3,V_4$ be $\mathcal{W}$-representations. 

1) Suppose that $V_1,V_2$ are irreducible, and that there exist trivial $\mathcal{W}$-representations $M_1$ and $M_2$ such that $V_4 \cong M_1 \otimes V_2$ and $V_2\otimes V_3 \cong M_2 \otimes V_1$ hold as isomorphisms of $\mathcal{W}$-representations. Then the linear map
$$
\mathrm{Hom}_\mathcal{W}(V_1, V_2\otimes V_3) \otimes \mathrm{Hom}_\mathcal{W}(V_2,V_4) \to \mathrm{Hom}_\mathcal{W}(V_1, V_4 \otimes V_3),
$$
sending $f \otimes g$ to $(g\otimes \mathrm{id}_{V_3})\circ f$ for each $f\in \mathrm{Hom}_\mathcal{W}(V_1,V_2\otimes V_3)$ and $g\in \mathrm{Hom}_\mathcal{W}(V_2,V_4)$, is an isomorphism.

\vs

2) Suppose that $V_1,V_3$ are irreducible, and that there exist trivial $\mathcal{W}$-representations $M_1$ and $M_2$ such that $V_4 \cong M_1 \otimes V_3$ and $V_1 \otimes V_3 \cong M_2 \otimes V_1$ hold as isomorphisms of $\mathcal{W}$-representations. Then the linear map
$$
\mathrm{Hom}_\mathcal{W}(V_3,V_4) \otimes \mathrm{Hom}_\mathcal{W}(V_1,V_2\otimes V_3) \to \mathrm{Hom}_\mathcal{W}(V_1, V_2 \otimes V_4),
$$
sending $f\otimes g$ to $(\mathrm{id}_{V_2} \otimes f)\circ g$ for each $f\in \mathrm{Hom}_\mathcal{W}(V_3,V_4)$ and $g\in \mathrm{Hom}_\mathcal{W}(V_1,V_2\otimes V_3)$, is an isomorphism.
\end{lemma}

{\it Proof.} 1) Consider the chain of isomorphisms
\begin{align*}
\mathrm{Hom}_\mathcal{W}(V_1, V_4 \otimes V_3) & \stackrel{\tiny \circled{1}}{\to} \mathrm{Hom}_\mathcal{W}(V_1, M_1 \otimes V_2 \otimes V_3) \stackrel{\tiny \circled{2}}{\to} \mathrm{Hom}_\mathcal{W}(V_1, M_1 \otimes M_2 \otimes V_1) \\
& \stackrel{\tiny \circled{3}}{\to} M_1 \otimes M_2 \otimes \mathrm{Hom}_\mathcal{W}(V_1, V_1) \stackrel{\tiny \circled{4}}{\to} M_1 \otimes M_2 \stackrel{\tiny \circled{5}}{\to} M_2 \otimes M_1 \\
& \stackrel{\tiny \circled{6}}{\to} \mathrm{Hom}_\mathcal{W}(V_1,V_2\otimes V_3) \otimes \mathrm{Hom}_\mathcal{W}(V_2, V_4).
\end{align*}
\circled{1} is induced by $V_4 \cong M_1 \otimes V_2$, \circled{2} by $V_2 \otimes V_3 \cong M_2 \otimes V_1$, \circled{3} follows from Lem.\ref{lem:taking_out_triv}, \circled{4} sends $f_1 \otimes f_2 \otimes \mathrm{id}$ to $f_1\otimes f_2$ (we used $\mathrm{Hom}_\mathcal{W}(V_1,V_1) = \mathbb{C} \cdot \mathrm{id}_{V_1}$), \circled{5} is just factor permuting map $\mathbf{P}_{(12)}$, and \circled{6} is from Cor.\ref{cor:M_as_Hom}. 

\vs

To check that this is the asserted map, pick any $f_1\in M_1, f_2\in M_2$, and denote the given isomorphisms by $F_1 : V_4 \to M_1 \otimes V_2$ and $F_2 : V_2\otimes V_3 \to M_2 \otimes V_1$. Then, by \circled{5},\circled{6} the element $f_1\otimes f_2 \in M_1\otimes M_2$ is sent to $f\otimes g$, with $f\in \mathrm{Hom}_\mathcal{W}(V_1,V_2\otimes V_3)$ given by $f(v_1) = F_2^{-1}(f_2 \otimes v_1)$, $\forall v_1 \in V_1$, and $g\in \mathrm{Hom}_\mathcal{W}(V_2,V_4)$ given by $g(v_2) = F_1^{-1}(f_1 \otimes v_2)$, $\forall v_2 \in V_2$. Meanwhile, by the inverses of \circled{4}, \circled{3}, \circled{2}, \circled{1}, the element $f_1 \otimes f_2$ is sent to the element of $\mathrm{Hom}_\mathcal{W}(V_1,V_4\otimes V_3)$ sending each $v_1\in V_1$ to $(F_1^{-1} \otimes \mathrm{id})(\mathrm{id}\otimes F_2^{-1}) (f_1 \otimes f_2\otimes v_1)$, which equals $(g\otimes \mathrm{id}_{V_3})f(v_1)$, as desired. We leave 2) as an exercise to the readers. \qed

\subsection{The factor-permuting operator ${\bf A}$}

As promised in \S\ref{subsec:CG}, we first give an identification of our multiplicity space $M_{\lambda,\lambda'}$ with the space of intertwiners $\mathrm{Hom}_\mathcal{W}(V_{\lambda\lambda'}, V_\lambda\otimes V_{\lambda'})$, using Cor.\ref{cor:M_as_Hom}.
\begin{definition}[identification of the multiplicity space and the space of intertwiners]
\label{def:I}
For a regular pair $(\lambda,\lambda')$ of weights, denote by
\begin{align}
\label{eq:I}
I^{\lambda\lambda'}_{\lambda,\lambda'} : M^{\lambda\lambda'}_{\lambda,\lambda'} \to \mathrm{Hom}_\mathcal{W}(V_{\lambda\lambda'}, V_\lambda\otimes V_{\lambda'})
\end{align}
the vector space isomorphism given in Cor.\ref{cor:M_as_Hom}, which is given by
$$
(I^{\lambda\lambda'}_{\lambda,\lambda'} f)(v) = \mathbf{F}^{-1}_{\lambda,\lambda'}(f \otimes v), \quad \forall f\in M^{\lambda\lambda'}_{\lambda,\lambda'}, \quad \forall v\in V_{\lambda\lambda'}.
$$
\end{definition}

\vs

Now, consider the maps
\begin{align*}
\mathrm{Hom}^\mathrm{L}_\mathbb{C}(V_{\lambda\lambda'}, V_\lambda \otimes V_{\lambda'})
\to V_\lambda \otimes V_{\lambda'} \otimes V_{\lambda\lambda'}^*
\to {}^*V_{\lambda^*} \otimes V_{\lambda'} \otimes V_{(\lambda\lambda')^*}
\to \mathrm{Hom}^\mathrm{R}_\mathbb{C}(V_{\lambda^*}, V_{\lambda'} \otimes V_{(\lambda\lambda')^*}),
\end{align*}
where the first and the last arrows are the canonical maps $J^\mathrm{L}$ and $J^\mathrm{R}$ in \eqref{eq:Hom_as_tensor_products}, and the middle one is $D_{\lambda^*} \otimes \mathrm{id} \otimes C_{\lambda\lambda'}$ (see Prop.\ref{prop:isomorphisms_to_duals}); so all three arrows are $\mathcal{W}$-module isomorphisms. So these maps restrict to isomorphisms among the corresponding invariant subspaces. In view of Lem.\ref{lem:space_of_intertwiners_as_invariant_subspace}, the restriction of the composition of these maps to the invariant subspaces yield an isomorphism
\begin{align}
\label{eq:A_as_Hom}
{\bf A}^\mathrm{Hom} = (\mathbf{A}^{\lambda\lambda'}_{\lambda,\lambda'})^\mathrm{Hom} : \mathrm{Hom}_\mathcal{W}(V_{\lambda\lambda'}, V_\lambda\otimes V_{\lambda'}) \to \mathrm{Hom}_\mathcal{W}(V_{\lambda^*}, V_{\lambda'} \otimes V_{(\lambda\lambda')^*}),
\end{align}
which via the identification maps defined in Def.\ref{def:I} translates to the invertible linear map
$$
\mathbf{A} = \mathbf{A}_{\lambda,\lambda'}^{\lambda\lambda'} : M^{\lambda\lambda'}_{\lambda,\lambda'} \to M^{\lambda^*}_{\lambda', (\lambda\lambda')^*}.
$$
This can be thought of as cyclically permuting the roles of the three weights $\lambda,\lambda',\lambda \lambda'$ with certain dualizing; note that this initial triple $(\lambda,\lambda',\lambda\lambda')$ of non-singular weights is characterized by the fact that the third is the product of the first and the second. We first cyclically shift to $(\lambda',\lambda\lambda',\lambda)$ and dualize the new second and third, to get the triple $( \lambda',(\lambda\lambda')^*,\lambda^* )$. I claim that this is also a triple of non-singular weights whose third is the product of the first and the second. Non-singularity follows from Lem.\ref{lem:dual_weight}, and one just needs to verify 
$
\lambda' \, (\lambda\lambda')^* = \lambda^*,
$
which follows easily from Cor.\ref{cor:dual_weight_of_product} and Lem.\ref{lem:dual_weight}.

\vs

Apply this transformation of triples of weights
$$
(\lambda,\lambda',\lambda\lambda') \leadsto ( \lambda', (\lambda\lambda')^*, \lambda^*)
$$
three times and we are back into the same situation as the first triple
$$
(\lambda,\lambda',\lambda\lambda') \leadsto ( \lambda', (\lambda\lambda')^*, \lambda^*) \leadsto ((\lambda\lambda')^*, \lambda, (\lambda')^*) \leadsto (\lambda, \lambda', \lambda\lambda');
$$
so it is reasonable to expect:
\begin{proposition}[order three relation]
\label{prop:AAA}
The operator $\mathbf{A}^3$ is a scalar operator. More precisely, the composition
$$
\xymatrix@C+8mm{
M_{\lambda,\lambda'}^{\lambda\lambda'} \ar[r]^-{\mathbf{A}_{\lambda,\lambda'}^{\lambda\lambda'}} & M_{\lambda', (\lambda\lambda')^*}^{\lambda^*} \ar[r]^-{\mathbf{A}_{\lambda',(\lambda\lambda')^*}^{\lambda^*} } & M_{(\lambda\lambda')^*, \lambda}^{(\lambda')^*}  \ar[r]^-{\mathbf{A}_{(\lambda\lambda')^*, \lambda}^{(\lambda')^*} } & M_{\lambda,\lambda'}^{\lambda\lambda'}.
}
$$
is $q^{-2m_{x,x'}}$ times the identity operator.
\end{proposition}

{\it Proof.} We shall prove that the composition of three maps \eqref{eq:A_as_Hom} on the $\mathrm{Hom}_\mathcal{W}$'s is the identity. This composition can be obtained first by composing all arrows in the following diagram and then restrict the resulting map $\mathrm{Hom}^\mathrm{L}_\mathbb{C}(V_{\lambda\lambda'}, V_\lambda \otimes V_{\lambda'}) \to \mathrm{Hom}^\mathrm{R}_\mathbb{C}(V_{\lambda\lambda'}, V_\lambda \otimes V_{\lambda'})$ to the invariant subspaces:
\begin{align*}
\hspace{-10mm} \xymatrix@C-4mm@R-3mm{
\mathrm{Hom}^\mathrm{L}_\mathbb{C}(V_{\lambda\lambda'},V_\lambda\otimes V_{\lambda'})
\ar[r]^-{\tiny \circled{1}} &  V_\lambda \otimes V_{\lambda'} \otimes V^*_{\lambda\lambda'}
\ar[r]^-{\tiny \circled{2}} & {}^*V_{\lambda^*} \otimes V_{\lambda'} \otimes V_{(\lambda\lambda')^*}
\ar[r]^-{\tiny \circled{3}} & \mathrm{Hom}^\mathrm{R}_\mathbb{C}(V_{\lambda^*}, V_{\lambda'} \otimes V_{(\lambda\lambda')^*}) \ar[d]^-{\tiny \circled{4}} \\
\mathrm{Hom}^\mathrm{R}_\mathbb{C} (V_{(\lambda')^*}, V_{(\lambda\lambda')^*} \otimes V_\lambda) \ar[d]^-{\tiny \circled{8}} & {}^* V_{(\lambda')^*} \otimes V_{(\lambda\lambda')^*} \otimes V_\lambda \ar[l]_-{\tiny \circled{7}} & V_{\lambda'} \otimes V_{(\lambda\lambda')^*} \otimes V_{\lambda^*}^* \ar[l]_-{\tiny \circled{6}} &  \mathrm{Hom}^\mathrm{L}_\mathbb{C}(V_{\lambda^*}, V_{\lambda'} \otimes V_{(\lambda\lambda')^*}) \ar[l]_-{\tiny \circled{5}} \\
\mathrm{Hom}^\mathrm{L}_\mathbb{C}(V_{(\lambda')^*}, V_{(\lambda\lambda')^*} \otimes V_\lambda) \ar[r]^-{\tiny \circled{9}} & V_{(\lambda\lambda')^*} \otimes V_\lambda \otimes V^*_{(\lambda')^*} \ar[r]^-{\tiny \circled{10}} & {}^* V_{\lambda\lambda'} \otimes V_\lambda \otimes V_{\lambda'} \ar[r]^-{\tiny \circled{11}} & \mathrm{Hom}^\mathrm{R}_\mathbb{C}(V_{\lambda\lambda'}, V_\lambda \otimes V_{\lambda'}),
}
\end{align*}
where we enumerated the arrows for convenience. The maps \circled{4} and \circled{9} are identity maps on the corresponding $\mathrm{Hom}_\mathbb{C}(\sim, \sim)$, the maps \circled{1},\circled{3},\circled{5},\circled{7},\circled{9},\circled{11} are canonical maps in \eqref{eq:Hom_as_tensor_products}, and the remaining three maps are
$$
\circled{2} = (D_{\lambda^*})_1 (C_{\lambda\lambda'})_3, \quad \circled{6} = (D_{(\lambda')^*})_1 (C_{\lambda^*})_3, \quad \circled{10} = (D_{\lambda\lambda'})_1 (C_{(\lambda')^*})_3.
$$
Observe by Lem.\ref{lem:composition_of_canonical_maps} that the composition $\circled{5} \circ \circled{4} \circ \circled{3}$ is just the factor permuting map $\mathbf{P}_{(132)}$, sending each $f\otimes v_1 \otimes v_2$ to $v_1 \otimes v_2 \otimes f$. By the way, on any $n$-fold tensor product vector space and for any permutation $\gamma$ of $1,2,\ldots,n$, one defines the permutation map $\mathbf{P}_\gamma$ as the one sending each $i$-th tensor factor to $\gamma(i)$-th tensor factor, as seen in \S\ref{subsec:T} for the case when $\gamma$ is a transposition and just now for the case when $\gamma = (132)$. The properties we use are
\begin{align*}
& \mathbf{P}_{\mathrm{id}} = \mathrm{id}, \qquad \mathbf{P}_\gamma \circ \mathbf{P}_\eta = \mathbf{P}_{\gamma\circ \eta}, \qquad (\mbox{hence} ~ \mathbf{P}_\gamma^{-1} = \mathbf{P}_{\gamma^{-1}}) \\
& \mathbf{P}_\gamma \, (\mbox{expression on $i_1,i_2,\ldots$-th factors}) \, \mathbf{P}_{\gamma^{-1}} = (\mbox{same expression on $\gamma(i_1),\gamma(i_2),\ldots$-th factors}).
\end{align*}
Likewise, $\circled{9}\circ \circled{8} \circ \circled{7}$ and $\circled{1} \circ \mathrm{id}\circ \circled{11}$ are also $\mathbf{P}_{(132)}$. So the composition \circled{2},\circled{3},\ldots,\circled{10} is given by
\begin{align*}
& (D_{\lambda\lambda'})_1 (C_{(\lambda')^*})_3 \, \mathbf{P}_{(132)} \, (D_{(\lambda')^*})_1 (C_{\lambda^*})_3 \, \mathbf{P}_{(132)} \, (D_{\lambda^*})_1 (C_{\lambda\lambda'})_3 \\
& = \mathbf{P}_{(132)} \, (D_{\lambda\lambda'})_2 (C_{(\lambda')^*})_1 (D_{(\lambda')^*})_1 (C_{\lambda^*})_3 \, \mathbf{P}_{(132)} \, (D_{\lambda^*})_1 (C_{\lambda\lambda'})_3 \\
& = \mathbf{P}_{(132)} \, \mathbf{P}_{(132)} \, (D_{\lambda\lambda'})_3 (C_{(\lambda')^*})_2 (D_{(\lambda')^*})_2 (C_{\lambda^*})_1 \, (D_{\lambda^*})_1 (C_{\lambda\lambda'})_3 \\
& = \mathbf{P}_{(123)} \, (C_{\lambda^*} \, D_{\lambda^*})_1 \, (C_{(\lambda')^*} \, D_{(\lambda')^*})_2 \, (D_{\lambda\lambda'} \, C_{\lambda\lambda'})_3.
\end{align*}
For any non-singular weight $\lambda = (x,y)$, from the definition of the linear maps $D_\lambda$ and $C_\lambda$ in Prop.\ref{prop:isomorphisms_to_duals}, observe
\begin{align*}
& D_\lambda \, C_\lambda : e_i \mapsto q^{-2im_{y,x^{-1}}}q^{-i(i-1)} q^{2im_{y,x^{-1}}} q^{-i(-i-1)}e_i = q^{2i}e_i = Ae_i, \\
& C_\lambda \, D_\lambda : e_i \mapsto q^{-2im_{y,x^{-1}}} q^{i(i-1)} q^{2im_{y,x^{-1}}} e^{i(-i-1)} e_i = q^{-2i} e_i = A^{-1} e_i,
\end{align*}
so
\begin{align}
\label{eq:DC_and_CD}
D_\lambda \, C_\lambda = A, \qquad C_\lambda \, D_\lambda = A^{-1},
\end{align}
as linear maps on $\mathbb{C}^N$, for any non-singular weight $\lambda$. Hence,
$$
\circled{10} \circ \cdots \circ \circled{2} = \mathbf{P}_{(123)} \, A_1^{-1} A_2^{-1} A_3.
$$
Consider applying this to an element of $\mathrm{Inv}( V_\lambda \otimes V_{\lambda'} \otimes V_{\lambda\lambda'}^* )$. This element is invariant under the action of $\mathcal{W}$; in particular, the $X$-action is $(\Delta \circ \mathrm{id})(\Delta X) = X\otimes X\otimes X = x^{1/N} (x')^{1/N}(xx')^{-1/N} A_1 A_2A_3^{-1} = q^{-2m_{x,x'}} A_1A_2A_3^{-1}$ in view of \eqref{eq:right_dual_of_cyclic_irreducible}, $\lambda\lambda' = (xx', yx' + y')$, and \eqref{eq:m_definition}, while $\epsilon(X)=1$. Hence the application of $A_1^{-1} A_2^{-1} A_3$ to this element results in multiplcation by $q^{-2m_{x,x'}}$. Thus, for any element $\xi \in \mathrm{Inv}(\mathrm{Hom}^\mathrm{L}_\mathbb{C}(V_{\lambda\lambda'},V_\lambda\otimes V_{\lambda'}))$, the application of the composed map $\circled{11} \circ \cdots \circ \circled{1}$ results in
$$
\circled{11} \circ \cdots \circ \circled{1} \, \xi = \circled{11} \, \mathbf{P}_{(123)} \, A_1^{-1} A_2^{-1} A_3 \, (\circled{1} \, \xi) = \circled{11} \, \mathbf{P}_{(123)} \, q^{-2m_{x,x'}} \, \circled{1} \, \xi
= q^{-2m_{x,x'}} \, \xi,
$$
where we used $\circled{1} \, \xi \in \mathrm{Inv}( V_\lambda \otimes V_{\lambda'} \otimes V^*_{\lambda\lambda'}) \otimes V_\lambda \otimes V_{\lambda'})$ and $\circled{1} \circ \circled{11} = \mathbf{P}_{(132)} = \mathbf{P}_{(123)}^{-1}$. So the restriction to the invariant subspaces of the composition of the maps \circled{1},\ldots,\circled{11} is just the scalar multiplication by $q^{-2m_{x,x'}}$, on $\mathrm{Hom}_\mathcal{W}(V_{\lambda\lambda'}, V_\lambda\otimes V_{\lambda'})$. Translating this to the map $M_{\lambda,\lambda'}^{\lambda\lambda'} \to M_{\lambda,\lambda'}^{\lambda\lambda'}$ via the identification $I_{\lambda,\lambda'}^{\lambda\lambda'} : M_{\lambda,\lambda'}^{\lambda\lambda'} \to \mathrm{Hom}_\mathcal{W}(V_{\lambda\lambda'}, V_\lambda\otimes V_{\lambda'})$ in Def.\ref{def:I}, we get $(I_{\lambda,\lambda'}^{\lambda\lambda'})^{-1} \, (q^{-2m_{x,x'}} \cdot \mathrm{id}) \, I_{\lambda,\lambda'}^{\lambda\lambda'} = q^{-2m_{x,x'}} \cdot \mathrm{id}$, as asserted. \qed

\vs

The above proposition says
$$
\mathbf{A}_{(\lambda\lambda')^*,\lambda}^{(\lambda')^*} \circ \mathbf{A}_{\lambda', (\lambda\lambda')^*}^{\lambda^*} \circ \mathbf{A}_{\lambda,\lambda'}^{\lambda\lambda'} = q^{-2m_{x,x'}}\cdot \mathrm{id}.
$$
If we apply the same proposition after replacing the initial triple $(\lambda,\lambda',\lambda\lambda')$ by $(\lambda',(\lambda\lambda')^*,\lambda^*)$, we see that we must have $q^{-2m_{x,x'}} = q^{-2m_{ x', (xx')^{-1} }} = q^{ - 2m_{(xx')^{-1}, x}}$, for consistency. Indeed, it is a straightforward exercise to show $m_{x,x'} = m_{x',(xx')^{-1}} = m_{(xx')^{-1}, x}$.

\section{Computation of formulas of the operators ${\bf T}$ and ${\bf A}$}
\label{sec:computation_of_formulas}

In the present section we compute explicit formulas of the operators ${\bf T}$ and ${\bf A}$ constructed in the previous sections. If one is not so interested in these formulas, one may skip this section and go to the next one.

\subsection{Quantum pentagon identity of cyclic quantum dilogarithm}
\label{subsec:pentagon_identity_of_Phi_C}

The following looks like what was established in \cite{FK94} \cite{BB}, but a little different because the settting and definition for the cyclic quantum dilogarithm are modified in the present paper (see \S\ref{subsec:Phi_C}). So we give an explicit proof.
\begin{proposition}[quantum pentagon identity; primitive form]
\label{prop:pentagon_primitive}
For complex numbers ${\bf a}_1$, \ldots, ${\bf a}_5$, ${\bf c}_1$, \ldots, ${\bf c}_5$ and an integer $m$ satisfying
\begin{align}
\label{eq:condition_of_numbers}
{\bf c}_1 = {\bf c}_4 {\bf c}_5, \quad
{\bf a}_1 {\bf c}_2 = {\bf c}_4 {\bf a}_5, \quad
{\bf a}_1 {\bf a}_2 = - q^{2m} {\bf a}_4, \quad
{\bf c}_2 = {\bf c}_3 {\bf c}_4, \quad
{\bf a}_2 {\bf c}_5 = {\bf a}_3,
\end{align}
one has the identity of operators
$$
\Phi^q_{{\bf a}_2, {\bf c}_2}(B) \, \Phi^q_{{\bf a}_1, {\bf c}_1}(A)
= \alpha \, \Phi^q_{{\bf a}_5, {\bf c}_5}(A) \, \Phi^q_{ {\bf a}_4, {\bf c}_4 }(q^{2m} AB) \, \Phi^q_{{\bf a}_3, {\bf c}_3}(B),
$$
for some constant $\alpha$ depending on the above eleven numbers ${\bf a}_1,\ldots,{\bf c}_1,\ldots,m$.
\end{proposition}

{\it Proof.} Each of the five factors is invertible. The expression $\Phi^q_{{\bf a}_4, {\bf c}_4}(q^{2m}AB)$ makes sense because $(q^{2m}AB)^N = \mathrm{id}$. Let's compare the conjugation actions of both sides on the operators $A$ and $B$. Observe
\begin{align*}
& \Phi^q_{{\bf a}_2, {\bf c}_2}(B) \, \Phi^q_{{\bf a}_1, {\bf c}_1}(A) \, B^{-1} \, \Phi^q_{{\bf a}_1, {\bf c}_1}(A)^{-1} \, \Phi^q_{{\bf a}_2, {\bf c}_2}(B)^{-1}
= \Phi^q_{{\bf a}_2, {\bf c}_2}(B) \, ( {\bf c}_1 B^{-1} - {\bf a}_1 B^{-1} A) \, \Phi^q_{{\bf a}_2, {\bf c}_2}(B)^{-1} \\
& = {\bf c}_1 B^{-1} - {\bf a}_1 ( {\bf c}_2 B^{-1} A - {\bf a}_2 B^{-1} A B) = {\bf c}_1 B^{-1} - {\bf a}_1 {\bf c}_2 B^{-1} A + {\bf a}_1 {\bf a}_2 q^2 A,
\end{align*}
which would equal
\begin{align*}
& \Phi^q_{{\bf a}_5, {\bf c}_5}(A) \, \Phi^q_{ {\bf a}_4, {\bf c}_4 }(q^{2m} AB) \, \Phi^q_{{\bf a}_3, {\bf c}_3}(B) \, B^{-1} \Phi^q_{{\bf a}_3, {\bf c}_3}(B)^{-1} \, \Phi^q_{{\bf a}_4, {\bf c}_4}(q^{2m} AB)^{-1} \, \Phi^q_{{\bf a}_5, {\bf c}_5}(A)^{-1} \\
& = \Phi^q_{{\bf a}_5, {\bf c}_5}(A) \, \Phi^q_{ {\bf a}_4, {\bf c}_4 }(q^{2m} AB) \, B^{-1} \, \Phi^q_{{\bf a}_4, {\bf c}_4}(q^{2m} AB)^{-1} \, \Phi^q_{{\bf a}_5, {\bf c}_5}(A)^{-1} \\ 
& = \Phi^q_{{\bf a}_5, {\bf c}_5}(A) \, ( {\bf c}_4 B^{-1} - {\bf a}_4 B^{-1} q^{2m} AB) \, \Phi^q_{{\bf a}_5, {\bf c}_5}(A)^{-1}
=  \Phi^q_{{\bf a}_5, {\bf c}_5}(A) \, ( {\bf c}_4 B^{-1} - {\bf a}_4 q^{2m} q^2 A ) \, \Phi^q_{{\bf a}_5, {\bf c}_5}(A)^{-1} \\
& = {\bf c}_4 ({\bf c}_5 B^{-1} - {\bf a}_5 B^{-1} A) - {\bf a}_4 q^{2m} q^2 A = {\bf c}_4 {\bf c}_5 B^{-1} - {\bf c}_4 {\bf a}_5 B^{-1} A - {\bf a}_4 q^{2m} q^2 A
\end{align*}
because of the stipulated condition of the numbers. Similarly, observe
\begin{align*}
& \Phi^q_{{\bf a}_5, {\bf c}_5}(A)^{-1} \, \Phi^q_{{\bf a}_2, {\bf c}_2}(B) \, \Phi^q_{{\bf a}_1, {\bf c}_1}(A) \, A \, \Phi^q_{{\bf a}_1, {\bf c}_1}(A)^{-1} \, \Phi^q_{{\bf a}_2, {\bf c}_2}(B)^{-1} \, \Phi^q_{{\bf a}_5, {\bf c}_5}(A) \\
& = \Phi^q_{{\bf a}_5, {\bf c}_5}(A)^{-1} \, \Phi^q_{{\bf a}_2, {\bf c}_2}(B) \, A \, \Phi^q_{{\bf a}_2, {\bf c}_2}(B)^{-1} \, \Phi^q_{{\bf a}_5, {\bf c}_5}(A)
 = \Phi^q_{{\bf a}_5, {\bf c}_5}(A)^{-1} ( {\bf c}_2 A - {\bf a}_2 AB) \Phi^q_{{\bf a}_5, {\bf c}_5}(A) \\
& = {\bf c}_2 A - {\bf a}_2 ( {\bf c}_5 AB - {\bf a}_5 A AB)
= {\bf c}_2 A - {\bf a}_2 {\bf c}_5 AB + {\bf a}_2 {\bf a}_5 A^2 B,
\end{align*}
which would equal
\begin{align*}
& \Phi^q_{{\bf a}_4, {\bf c}_4}(q^{2m} AB) \, \Phi^q_{{\bf a}_3, {\bf c}_3}(B) \, A \, \Phi^q_{{\bf a}_3, {\bf c}_3}(B)^{-1} \, \Phi^q_{{\bf a}_4, {\bf c}_4}(q^{2m} AB)^{-1} \\
& = \Phi^q_{{\bf a}_4, {\bf c}_4}(q^{2m} AB) \, ( {\bf c}_3 A - {\bf a}_3 AB) \, \Phi^q_{{\bf a}_4, {\bf c}_4}(q^{2m} AB)^{-1} \\
& = {\bf c}_3 ({\bf c}_4 A - {\bf a}_4 A q^{2m} AB) - {\bf a}_3 AB
= {\bf c}_3 {\bf c}_4 A - {\bf c}_3 {\bf a}_4 q^{2m} A^2 B - {\bf a}_3 AB,
\end{align*}
because of the stipulated condition of the numbers; we used ${\bf a}_2 {\bf a}_5 = - q^{2m} {\bf c}_3 {\bf a}_4$, which follows from the conditions: $\frac{{\bf a}_2 {\bf a}_5}{-q^{2m} {\bf c}_3 {\bf a}_4} = \frac{{\bf a}_1 {\bf a}_2}{-q^{2m}{\bf a}_4} \cdot \frac{{\bf c}_4 {\bf a}_5}{{\bf a}_1 {\bf c}_2} \cdot \frac{{\bf c}_2}{{\bf c}_3{\bf c}_4} = 1$. So we verified that the conjugation actions of the LHS and the RHS on $B^{-1}$ (hence on $B$) are equal, and those on $A$ are equal. So the assertion follows from Cor.\ref{cor:conjugation_on_A_B_determines}.  \qed

\begin{theorem}[quantum pentagon identity; see \cite{FK94} \cite{BB} for a different version]
\label{thm:pentagon}
For complex numbers ${\bf a}_1,\ldots,{\bf a}_5, {\bf c}_1, \ldots, {\bf c}_5$, $\alpha$ as in Prop.\ref{prop:pentagon_primitive}, and for linear operators $C,D$ on a complex vector space satisfying $C^N = D^N = \mathrm{id}$ and $CD = q^2 DC$, one has
$$
\Phi^q_{{\bf a}_2, {\bf c}_2}(D) \, \Phi^q_{{\bf a}_1, {\bf c}_1}(C)
= \alpha \, \Phi^q_{{\bf a}_5, {\bf c}_5}(C) \, \Phi^q_{ {\bf a}_4, {\bf c}_4 }(q^{2m} CD) \, \Phi^q_{{\bf a}_3, {\bf c}_3}(D).
$$
\end{theorem}

{\it Proof.} The assignment of operators $X\mapsto C$ and $Y\mapsto D$ provides a cyclic representation of $\mathcal{W}$, hence by Prop.\ref{prop:cyclic_representations} decomposes into the direct sum of cyclic irreducibles; as $C^N = D^N = \mathrm{id}$, each of these cyclic irreducibles is labeled by the weight $\lambda = (1,1)$ hence is isomorphic to the one given by the assignment of operators $X \mapsto A$ and $Y \mapsto B$ on $\mathbb{C}^N$. For the operators $A$ and $B$ we proved the identity in Prop.\ref{prop:pentagon_primitive}. \qed

\subsection{The ${\bf T}$ operator}

How does the $6j$-symbol map $\mathbf{T}_{\lambda,\lambda',\lambda''}$ look like? In fact it almost looks like the decomposition map $\mathbf{F}$ itself:
\begin{proposition}[computation of the $6j$-symbol ${\bf T}$]
\label{prop:computation_of_T}
For any regular triple $(\lambda,\lambda',\lambda'')$ of weights, one has
$$
\mathbf{T}_{\lambda,\lambda',\lambda''} = \alpha \, \Phi^q_{{\bf a}, {\bf c}}(B_2 A_1 B_1^{-1})^{-1} \, {\bf S}_{21}^{-1} \, B_1^{m_{x,x'}},
$$
where $\alpha$ is some nonzero complex constant depending on $\lambda,\lambda',\lambda''$, and
\begin{align}
\label{eq:T_ac}
{\bf a} = - \frac{y^{1/N} (x')^{1/N} (y'')^{1/N}}{(yx'+y')^{1/N}(y'x''+y'')^{1/N}}, \qquad
{\bf c} = \frac{ (y')^{1/N} (yx'x'' + y'x'' + y'')^{1/N} }{(yx'+y')^{1/N} (y'x'' + y'')^{1/N} },
\end{align}
where $\lambda = (x,y)$, $\lambda'=(x',y')$, $\lambda''=(x'',y'')$.
\end{proposition}

{\it Proof.} We shall prove
\begin{align}
\label{eq:F_T_computation_objective}
(\mathbf{F}_{\lambda\lambda',\lambda''})_{23} \, (\mathbf{F}_{\lambda,\lambda'})_{12} = (\mathbf{F}_{\lambda''',\lambda''''})_{12} \, (\mathbf{F}_{\lambda,\lambda'\lambda''})_{13} \, (\mathbf{F}_{\lambda',\lambda''})_{23},
\end{align}
for some to-be-found regular pair $(\lambda''',\lambda'''')$ of weights. Then, from \eqref{eq:F_and_T} we would have $(\mathbf{T}_{\lambda,\lambda',\lambda''})_{21} = (\mathbf{F}_{\lambda''',\lambda''''})_{12}^{-1}$, yielding the desired result.  I first claim 
$$
\mathbf{S}_{23} \, \mathbf{S}_{12} = \mathbf{S}_{12} \, \mathbf{S}_{13} \, \mathbf{S}_{23},
$$
which can be directly verified on each basis vector $e_i \otimes e_j \otimes e_k$. In view of the definition \eqref{eq:S}, one has $\mathbf{S}_{23} \, \mathbf{S}_{12} \, e_i \otimes e_j \otimes e_k = \mathbf{S}_{23} \,e_i \otimes e_{i+j} \otimes e_k = e_i \otimes e_{i+j} \otimes e_{i+j+k}$, while $\mathbf{S}_{12}\, \mathbf{S}_{13}\, \mathbf{S}_{23}\, e_i \otimes e_j \otimes e_k = \mathbf{S}_{12} \, \mathbf{S}_{13}\, e_i \otimes e_j \otimes e_{j+k} = \mathbf{S}_{12} \, e_i \otimes e_j \otimes e_{i+j+k} = e_i \otimes e_{i+j} \otimes e_{i+j+k}$. Like in \eqref{eq:S_conjugation_1}, note $A_1 \mathbf{S} \, e_i \otimes e_j = A_1 \, e_i \otimes e_{i+j} = q^{2i} \, e_i \otimes e_{i+j} = \mathbf{S} (q^{2i} e_i \otimes e_j) = \mathbf{S}\, A_1 \, e_i \otimes e_j$ and $B_1 \mathbf{S} \, e_i \otimes e_j = B_1 e_i \otimes e_{i+j} = e_{i+1} \otimes e_{i+j} = \mathbf{S} \, e_{i+1} \otimes e_{j-1} = \mathbf{S} \, B_1 B_2^{-1} \, e_i \otimes e_{j-1}$. Summarizing, we have
\begin{align}
\label{eq:S_commuting}
A_2 \, \mathbf{S} = \mathbf{S} \, A_1A_2, \quad
B_2 \, \mathbf{S} = \mathbf{S} \, B_2, \quad
A_1 \, \mathbf{S} = \mathbf{S} \, A_1, \quad
B_1 \, \mathbf{S} = \mathbf{S} \, B_1 B_2^{-1},
\end{align}
which we will use when moving the factor $\mathbf{S}$ `to the left'. Recall also $AB^m = q^{2m}B^mA$ \eqref{eq:B_m_commuting} for each integer $m$, and that any two factors whose subscript indices do not intersect commute. 
\vs

Denote $m_{\lambda,\lambda'} = m_{x,x'}$ for any regular pair $(\lambda,\lambda')$ of non-singular weights; write $m_1 = m_{\lambda,\lambda'}$, $m_2 = m_{\lambda,\lambda'\lambda''}$, $m_3 = m_{\lambda',\lambda''}$, $m_4 = m_{\lambda\lambda',\lambda''}$, $m_5 = m_{\lambda''',\lambda''''}$. From \eqref{eq:F}, $\mathbf{F}_{\lambda,\lambda'} = (B_2)^{-m_{x,x'}} \, \mathbf{S} \, \Phi^q_{\lambda,\lambda'}(B_1A_2B_2^{-1})$, which equals $\mathbf{S} \, B_2^{-m_{\lambda,\lambda'}} \, \Phi^q_{\lambda,\lambda'}(B_1A_2B_2^{-1})$.  Note (at each stage, I underlined the part that is being moved to the left or being changed)
\begin{align*}
(\mathbf{F}_{\lambda''',\lambda''''})_{12} \, (\mathbf{F}_{\lambda,\lambda'\lambda''})_{13} \, (\mathbf{F}_{\lambda',\lambda''})_{23}
& = (\mathbf{S}_{12} \, B_2^{-m_5} \,  \Phi^q_{\lambda''',\lambda''''}(B_1A_2B_2^{-1})) \, (\ul{ \mathbf{S}_{13} } \, B_3^{-m_2}\, \Phi^q_{\lambda,\lambda'\lambda''}(B_1A_3B_3^{-1}) ) \\
& \quad \cdot(\mathbf{S}_{23} \, B_3^{-m_3} \, \Phi^q_{\lambda',\lambda''}(B_2A_3B_3^{-1})) \\
& = \mathbf{S}_{12} \, \mathbf{S}_{13}\, B_2^{-m_5} \,  \Phi^q_{\lambda''',\lambda''''}(B_1B_3^{-1} A_2B_2^{-1}) \, B_3^{-m_2} 
 \, \Phi^q_{\lambda,\lambda'\lambda''}(B_1A_3B_3^{-1}) \\
& \quad \cdot \ul{ \mathbf{S}_{23}  }\, B_3^{-m_3} \, \Phi^q_{\lambda',\lambda''}(B_2A_3B_3^{-1}) \quad (\because \eqref{eq:S_commuting}) \\
& = \mathbf{S}_{12} \, \mathbf{S}_{13}\, \mathbf{S}_{23}\, 
B_2^{-m_5} \, B_3^{m_1} \,  \Phi^q_{\lambda''',\lambda''''}(B_1\cancel{B_3^{-1}} A_2 B_2^{-1}\cancel{B_3}) \, \ul{ B_3^{-m_2} } \\
& \quad \cdot \Phi^q_{\lambda,\lambda'\lambda''}(B_1A_2A_3B_3^{-1})\, B_3^{-m_3} \, 
 \, \Phi^q_{\lambda',\lambda''}(B_2A_3B_3^{-1}) \quad (\because \eqref{eq:S_commuting}) \\
& = \mathbf{S}_{12} \, \mathbf{S}_{13}\, \mathbf{S}_{23}\, 
B_2^{-m_5} \, B_3^{m_1-m_2} \,  \Phi^q_{\lambda''',\lambda''''}(B_1 A_2 B_2^{-1}) \\
& \quad \cdot \Phi^q_{\lambda,\lambda'\lambda''}(B_1A_2A_3B_3^{-1})\, \ul{ B_3^{-m_3}} \, 
 \, \Phi^q_{\lambda',\lambda''}(B_2A_3B_3^{-1}) \quad (\because \eqref{eq:B_m_commuting}) \\
& = \mathbf{S}_{12} \, \mathbf{S}_{13}\, \mathbf{S}_{23}\, 
B_2^{-m_5} \, B_3^{m_1-m_2-m_3} \,  \Phi^q_{\lambda''',\lambda''''}(B_1 A_2 B_2^{-1}) \\
& \quad \cdot \Phi^q_{\lambda,\lambda'\lambda''}(q^{-2m_3}B_1A_2A_3B_3^{-1})\, 
 \, \Phi^q_{\lambda',\lambda''}(B_2A_3B_3^{-1}) \quad (\because \eqref{eq:B_m_commuting}),
\end{align*}
while
\begin{align*}
(\mathbf{F}_{\lambda\lambda',\lambda''})_{23} \, (\mathbf{F}_{\lambda,\lambda'})_{12} 
& = (\mathbf{S}_{23}\, B_3^{-m_4} \, \Phi^q_{\lambda\lambda',\lambda''}(B_2A_3B_3^{-1})) \, (\ul{ \mathbf{S}_{12} } \, B_2^{-m_1} \, \Phi^q_{\lambda,\lambda'}(B_1A_2B_2^{-1})) \\
& = \mathbf{S}_{23}\, \mathbf{S}_{12} \, B_3^{-m_4} \, \Phi^q_{\lambda\lambda',\lambda''}(B_2A_3B_3^{-1})  \, \ul{ B_2^{-m_1} } \, \Phi^q_{\lambda,\lambda'}(B_1A_2B_2^{-1}) \quad (\because \eqref{eq:S_commuting}) \\
& = \mathbf{S}_{23}\, \mathbf{S}_{12} \, B_2^{-m_1} \, B_3^{-m_4} \, \Phi^q_{\lambda\lambda',\lambda''}(B_2A_3B_3^{-1})   \, \Phi^q_{\lambda,\lambda'}(B_1A_2B_2^{-1}).
\end{align*}
We already proved ${\bf S}_{12} {\bf S}_{13} {\bf S}_{23} = {\bf S}_{23} {\bf S}_{12}$. Note $m_1-m_2-m_3 = -m_4 \Leftrightarrow m_1 + m_4 = m_2 +m_3$, which is equivalent to $m_{\lambda,\lambda'} + m_{\lambda\lambda',\lambda''} = m_{\lambda,\lambda'\lambda''} + m_{\lambda',\lambda''}$, which is true. We shall see that $m_1 = m_5$ can be attained. What remains to establish is (up to constant)
$$
\hspace{-8mm} \Phi^q_{\lambda''',\lambda''''}(B_1A_2B_2^{-1}) \Phi^q_{\lambda,\lambda'\lambda''}(q^{-2m_3}B_1A_2A_3B_3^{-1}) \Phi^q_{\lambda',\lambda''}(B_2A_3B_3^{-1})
= \Phi^q_{\lambda\lambda',\lambda''}(B_2A_3B_3^{-1}) \Phi^q_{\lambda,\lambda'}(B_1A_2B_2^{-1}).
$$
Put $C = B_1 A_2 B_2^{-1}$ and $D = B_2A_3 B_3^{-1}$, which are operators on $\mathbb{C}^N \otimes \mathbb{C}^N \otimes \mathbb{C}^N$. They satisfy $C^N = D^N = \mathrm{id}$ and $CD = q^2DC$. We now apply Thm.\ref{thm:pentagon}. For $q^{2m} CD = q^{-2m_3} B_1A_2A_3B_3^{-1} = q^{-2m_3} CD$, we put $m=-m_3= - m_{\lambda',\lambda''} = -m_{x',x''}$ into the theorem. For the other parameters, we match up factor by factor
\begin{align*}
& \Phi^q_{\mathbf{a}_1,\mathbf{c}_1} = \Phi^q_{\lambda,\lambda'} = \Phi^q_{(x,y),(x',y')}: \quad {\bf a}_1 = - \frac{ y^{1/N} (x')^{1/N} }{ (yx' + y')^{1/N} }, \quad {\bf c}_1 = \frac{ (y')^{1/N} }{ (yx' + y')^{1/N} } \\
& \Phi^q_{\mathbf{a}_2, \mathbf{c}_2} = \Phi^q_{\lambda\lambda',\lambda''} = \Phi^q_{(xx',yx'+y'), (x'',y'')} : \quad {\bf a}_2 = - \frac{ (yx'+y')^{1/N} (x'')^{1/N} }{ ( (yx'+y')x'' + y'')^{1/N} }, \\
& \qquad\qquad\qquad\qquad\qquad\qquad\qquad\qquad\qquad {\bf c}_2 = \frac{ (y'')^{1/N} }{( (yx'+y')x'' + y'')^{1/N}} \\
& \Phi^q_{{\bf a}_3, {\bf c}_3} = \Phi^q_{\lambda',\lambda''} = \Phi^q_{(x',y'),(x'',y'')} : \quad {\bf a}_3 = - \frac{ (y')^{1/N} (x'')^{1/N} }{(y'x'' + y'')^{1/N}}, \quad {\bf c}_3 = \frac{(y'')^{1/N}}{(y'x'' + y'')^{1/N}} \\
& \Phi^q_{{\bf a}_4, {\bf c}_4} = \Phi^q_{\lambda,\lambda'\lambda''} = \Phi^q_{(x,y),(x'x'',y'x''+y'')} : \quad {\bf a}_4 = - \frac{y^{1/N} (x'x'')^{1/N}}{(yx'x'' + y'x'' + y'')^{1/N}}, \\
& \qquad\qquad\qquad\qquad\qquad\qquad\qquad\qquad\qquad {\bf c}_4 = \frac{(y'x''+y'')^{1/N}}{(yx'x'' + y'x'' + y'')^{1/N}}, \\
& \Phi^q_{{\bf a}_5, {\bf c}_5} = \Phi^q_{\lambda''',\lambda''''}.
\end{align*}
In fact, we specify ${\bf a}_5$ and ${\bf c}_5$ as the numbers ${\bf a}$ and ${\bf c}$ in \eqref{eq:T_ac}, instead of weights $\lambda''',\lambda''''$. We leave it as an exercise to check ${\bf c}_5^N - {\bf a}_5^N =1$ and that the condition \eqref{eq:condition_of_numbers} is satisfied. So, in order for \eqref{eq:F_T_computation_objective} to hold, it suffices to have
$$
\mathbf{F}_{\lambda''',\lambda''''} = \alpha\, B_2^{-m_1} \, \mathbf{S} \, \Phi^q_{ {\bf a}_5, {\bf c}_5 }(B_1 A_2 B_2^{-1});
$$
that is, we use the RHS of this equation as the definition of the symbol on the LHS, instead of trying to find actual weights $\lambda''',\lambda''''$. \qed

\subsection{The ${\bf A}$ operator}
\label{subsec:computation_of_A}

As in the infinite dimensional case, the formula for the ${\bf A}$ operator is an analog of a `Fourier transform' or `exponential of quadratic expressions in Schr\"odinger operators' (see \S\ref{subsec:discussion_on_Phi_C} for more discussion on such operators). Instead of computing ${\bf A}$ directly using its definition which would be quite a complicated job, we compute its conjugation action on $A$ and $B$ to determine it indirectly, using Cor.\ref{cor:conjugation_on_A_B_determines}.

\begin{proposition}[computation of the $\mathbf{A}$ operator]
\label{prop:computation_of_A}
Let $(\lambda,\lambda')$ be a regular pair of weights, with $\lambda = (x,y)$, $\lambda'=(x',y')$. Then the map
$$
{\bf A}_{\lambda,\lambda'}^{\lambda\lambda'} : M_{\lambda,\lambda'}^{\lambda\lambda'} \to M_{\lambda',(\lambda\lambda')^*}^{\lambda^*}
$$
coincides up to multiplicative constant with the following linear map
\begin{align}
\label{eq:A_formula}
\mathbb{C}^N \to \mathbb{C}^N : e_i \mapsto \sum_{j=0}^{N-1} q^{-2ij - j^2 +(2m_1+1)j + 2(m_1-m_2) i} e_j, \qquad \forall i=0,1,\ldots,N-1,
\end{align}
where
\begin{align}
\label{eq:m1_and_m2}
m_1 := m_{x,x'} - m_{yx'+y', (xx')^{-1}}, \quad
m_2 := m_{-yx^{-1},x}.
\end{align}

\end{proposition}

{\it Proof.} Let $f \in M_{\lambda,\lambda'}^{\lambda\lambda'} \equiv \mathbb{C}^N$. Apply $I_{\lambda,\lambda'}^{\lambda\lambda'}$ in Def.\ref{def:I} to get $I_{\lambda,\lambda'}^{\lambda\lambda'} f \in \mathrm{Hom}_\mathcal{W}(V_{\lambda\lambda'}, V_\lambda \otimes V_{\lambda'})$, defined by $(I_{\lambda,\lambda'}^{\lambda\lambda'} f)(e_j) = \mathbf{F}_{\lambda,\lambda'}^{-1}(f \otimes e_j)$ for each $j$. Apply the canonical map $\mathrm{Hom}^\mathrm{L}_\mathbb{C}(V_{\lambda\lambda'},V_\lambda\otimes V_{\lambda'}) \to (V_\lambda\otimes V_{\lambda'}) \otimes V_{\lambda\lambda'}^*$ in \eqref{eq:Hom_as_tensor_products}; in view of the realization of $V_{\lambda\lambda'}^*$ as $\mathbb{C}^N$ as we did in \S\ref{subsec:isomorphisms_to_dual_representations} using the pairing \eqref{eq:canonical_pairing_on_C_N}, this map sends $I_{\lambda,\lambda'}^{\lambda\lambda'} f$ to $\sum_j \mathbf{F}^{-1}_{\lambda,\lambda'}(f\otimes e_j)\otimes e_j$. Then apply $D_{\lambda^*} \otimes \mathrm{id}\otimes C_{\lambda\lambda'}$ to land in the invariant subspace of ${}^*V_{\lambda^*} \otimes V_{\lambda'} \otimes V_{(\lambda\lambda')^*}$.

\vs

On the other hand, let us now start from the other end, i.e. from the image $\mathbf{A} f = \mathbf{A}_{\lambda,\lambda'}^{\lambda\lambda'} f \in M_{\lambda',(\lambda\lambda')^*}^{\lambda^*}$. Apply $I_{\lambda',(\lambda\lambda')^*}^{\lambda^*}$ in Def.\ref{def:I} to get $I_{\lambda',(\lambda\lambda')^*}^{\lambda^*}( \mathbf{A} f) \in \mathrm{Hom}_\mathcal{W} (V_{\lambda^*}, V_{\lambda'} \otimes V_{(\lambda\lambda')^*})$, defined by $( I_{\lambda',(\lambda\lambda')^*}^{\lambda^*} \mathbf{A}f)(e_j) = \mathbf{F}_{\lambda', (\lambda\lambda')^*}^{-1} ( (\mathbf{A}f) \otimes e_j)$ for each $j$. Apply the canonical map $\mathrm{Hom}_\mathbb{C}^\mathrm{R}(V_{\lambda^*}, V_{\lambda'} \otimes V_{(\lambda\lambda')^*}) \to {}^*V_{\lambda^*} \otimes (V_{\lambda'} \otimes V_{(\lambda\lambda')^*})$ in \eqref{eq:Hom_as_tensor_products}; using the realization ${}^*V_{\lambda^*}$ as $\mathbb{C}^N$, this map sends $I_{\lambda',(\lambda\lambda')^*}^{\lambda^*} \mathbf{A}f$ to $\sum_j e_j\otimes \mathbf{F}_{\lambda',(\lambda\lambda')^*}^{-1}( (\mathbf{A}f) \otimes e_j)$. Hence we are to solve the equation
\begin{align*}
\textstyle (D_{\lambda^*} \otimes \mathrm{id} \otimes C_{\lambda\lambda'})(\sum_j \mathbf{F}_{\lambda,\lambda'}^{-1} (f\otimes e_j)) \otimes e_j = \sum_j e_j \otimes \mathbf{F}_{\lambda', (\lambda\lambda')^*}^{-1} ( (\mathbf{A} f) \otimes e_j).
\end{align*}
By applying $D_{\lambda^*}^{-1} \otimes \mathrm{id} \otimes \mathrm{id}$ from the left on both sides, we get an equivalent equation
$$
\textstyle
\sum_j \mathbf{F}_{\lambda,\lambda'}^{-1} (f\otimes e_j ) \otimes (C_{\lambda\lambda'} e_j)
= \sum_j (D_{\lambda^*}^{-1} e_j) \otimes \mathbf{F}_{\lambda',(\lambda\lambda')^*}^{-1} ( (\mathbf{A}f)\otimes e_j)
$$
of elements in $\mathrm{Inv}(V_\lambda \otimes V_{\lambda'} \otimes V_{(\lambda\lambda')^*})$. Instead of putting in the formulas of $D_{\lambda^*}$, $C_{\lambda\lambda'}$, $\mathbf{F}_{\lambda,\lambda'}^{-1}$, and $\mathbf{F}_{\lambda',(\lambda\lambda')^*}^{-1}$ and compute what $\mathbf{A}$ is, we shall compute the conjugation action of $\mathbf{A}$ on $A$ and $B$, which will determine $\mathbf{A}$ up to scalar, in view of Cor.\ref{cor:conjugation_on_A_B_determines}.

\vs

To ease computation, we first collect conjugation actions of $\mathbf{F}_{\lambda,\lambda'}^{-1} = \Phi^q_{\lambda,\lambda'}(B_1A_2B_2^{-1})^{-1} \mathbf{S}^{-1} B_2^{m_{x,x'}}$. In Prop.\ref{prop:F} we already computed the conjugation on $A_2$ and $B_2$:
$$
\mathbf{F}_{\lambda,\lambda'}^{-1} \, A_2 \, \mathbf{F}_{\lambda,\lambda'}
= q^{-2m_{x,x'}} A_1 A_2, \qquad
\mathbf{F}_{\lambda,\lambda'}^{-1} \, B_2 \, \mathbf{F}_{\lambda,\lambda'} = \textstyle (\mathbf{c} \cdot \mathrm{id} - \mathbf{a} \, B_1 A_2 B_2^{-1} ) B_2,
$$
where ${\bf a} = - \frac{ y^{1/N} (x')^{1/N} }{(yx'+y')^{1/N}}$ and ${\bf c} = \frac{(y')^{1/N}}{(yx'+y')^{1/N}}$ as in \eqref{eq:our_ac}. Let's now compute conjugation on $A_1$ and $B_1$. First, note 
\begin{align*}
\mathbf{S}^{-1} A_1^{-1} = A_1^{-1} \mathbf{S}^{-1}, \qquad \mathbf{S}^{-1} B_1^{-1} = B_1^{-1} B_2 \mathbf{S}^{-1}.
\end{align*}
Indeed, on a basis vector $e_i\otimes e_j$, note $\mathbf{S}^{-1} (A^{-1}\otimes \mathrm{id})(e_i \otimes e_j) = \mathbf{S}^{-1}(q^{-2i} e_i \otimes e_j) = q^{-2i} e_i \otimes e_{j-i}$, while $(A^{-1} \otimes \mathrm{id}) \mathbf{S}^{-1} e_i \otimes e_j = (A^{-1} \otimes\mathrm{id}) e_i \otimes e_{j-i} = q^{-2i} e_i \otimes e_{j-i}$, and $\mathbf{S}^{-1}(B^{-1}\otimes\mathrm{id})(e_i\otimes e_j) = \mathbf{S}^{-1} e_{i-1} \otimes e_j = e_{i-1} \otimes e_{j-i+1}$ while $(B^{-1} \otimes B)\mathbf{S}^{-1} e_i \otimes e_j = (B^{-1}\otimes B) e_i \otimes e_{j-i} = e_{i-1} \otimes e_{j-i+1}$.

\vs

As $(B_1A_2B_2^{-1}) A_1^{-1} = q^2 A_1^{-1} (B_1A_2B_2^{-1})$, from Lem.\ref{lem:conjugation_by_Phi_C} or \eqref{eq:conjugation_by_Phi_C}
one has
$$
\textstyle \Phi^q_{\lambda,\lambda'}(B_1A_2B_2^{-1})^{-1} \, A_1^{-1}  = ({\bf c} \cdot \mathrm{id} - {\bf a} \, B_1A_2B_2^{-1}) \, A_1^{-1} \, \Phi^q_{\lambda,\lambda'}(B_1A_2B_2^{-1})^{-1},
$$
where $\Phi^q_{\lambda,\lambda'} = \Phi^q_{{\bf a},{\bf c}}$. Observe
\begin{align*}
\mathbf{F}_{\lambda,\lambda'}^{-1} A_1^{-1}
& = \Phi^q_{\lambda,\lambda'}(B_1A_2B_2^{-1})^{-1} \mathbf{S}^{-1} B_2^{m_{x,x'}} A_1^{-1} \\
& = \Phi^q_{\lambda,\lambda'}(B_1A_2B_2^{-1})^{-1} \mathbf{S}^{-1} A_1^{-1} B_2^{m_{x,x'}} \\
& = \Phi^q_{\lambda,\lambda'}(B_1A_2B_2^{-1})^{-1} A_1^{-1} \mathbf{S}^{-1} B_2^{m_{x,x'}} \\
& \textstyle = ({\bf c} \cdot \mathrm{id} - {\bf a} \, B_1 A_2 B_2^{-1} ) \, A_1^{-1}  \, \mathbf{F}_{\lambda,\lambda'}^{-1}
\end{align*}
Now, as $(B_1A_2B_2^{-1}) B_1^{-1} B_2 = q^2 B_1^{-1} B_2 (B_1A_2B_2^{-1})$, from Lem.\ref{lem:conjugation_by_Phi_C} or \eqref{eq:conjugation_by_Phi_C} one has
$$
\textstyle \Phi^q_{\lambda,\lambda'}(B_1A_2B_2^{-1})^{-1} \, B_1^{-1} B_2 =  ({\bf c} \cdot \mathrm{id} - {\bf a}\, B_1 A_2 B_2^{-1}) \, B_1^{-1} B_2 \, \Phi^q_{\lambda,\lambda'}(B_1A_2B_2^{-1})^{-1},
$$
so
\begin{align*}
\mathbf{F}_{\lambda,\lambda'}^{-1} \, B_1^{-1} & = \Phi^q_{\lambda,\lambda'}(B_1A_2B_2^{-1})^{-1} \mathbf{S}^{-1} B_2^{m_{x,x'}} B_1^{-1} \\
& = \Phi^q_{\lambda,\lambda'}(B_1A_2 B_2^{-1})^{-1} B_1^{-1} B_2 \, \mathbf{S}^{-1} B_2^{m_{x,x,'}} \\
& \textstyle = ({\bf c}\cdot \mathrm{id} - {\bf a} \, B_1A_2B_2^{-1}) \, B_1^{-1} B_2\, \mathbf{F}_{\lambda,\lambda'}^{-1}.
\end{align*}
Summarizing, we obtained
\begin{align*}
\textstyle
\mathbf{F}_{\lambda,\lambda'}^{-1} \, A_1^{-1} \, \mathbf{F}_{\lambda,\lambda'} = ({\bf c} \cdot \mathrm{id} - {\bf a}\, B_1A_2B_2^{-1}) \,A_1^{-1}, \qquad
\mathbf{F}_{\lambda,\lambda'}^{-1} \, B_1^{-1} \, \mathbf{F}_{\lambda,\lambda'} = ({\bf c} \cdot \mathrm{id} - {\bf a}\, B_1A_2 B_2^{-1}) \, B_1^{-1}B_2.
\end{align*}
Conjugation of an invertible operator by an invertible one is also invertible, so we deduce that $({\bf c}\cdot \mathrm{id} - {\bf a}\, B_1A_2B_2^{-1})$ is invertible. This allows us to make the above results a bit neater, by certain combinations:
\begin{align*}
\mathbf{F}_{\lambda,\lambda'}^{-1} \, A_1 B_2 \, \mathbf{F}_{\lambda,\lambda'} = A_1 B_2, \qquad
\mathbf{F}_{\lambda,\lambda'}^{-1} \, B_1 B_2 \, \mathbf{F}_{\lambda,\lambda'} = B_1.
\end{align*}
Then, it is easy to deduce the following conjugation results for $\mathbf{F}_{\lambda,\lambda'}$, which we shall use shortly:
\begin{align}
\label{eq:F_conjugation_on_A1B2_and_B1}
& \mathbf{F}_{\lambda,\lambda'} \, A_1 B_2 \, \mathbf{F}_{\lambda,\lambda'}^{-1} = A_1 B_2, \qquad\quad
\mathbf{F}_{\lambda,\lambda'} \, B_1 \, \mathbf{F}_{\lambda,\lambda'}^{-1} = B_1 B_2.
\end{align}
We will also need:
\begin{align*}
  \mathbf{F}_{\lambda,\lambda'} B_2^{-1} \mathbf{F}_{\lambda,\lambda'}^{-1}
& = B_2^{-m_{x,x'}} \mathbf{S} \, \Phi^q_{\lambda,\lambda'}(B_1A_2B_2^{-1}) B_2^{-1} \, \Phi^q_{\lambda,\lambda'}(B_1A_2B_2^{-1})^{-1} {\bf S}^{-1} B_2^{m_{x,x'}}\\
& = B_2^{-m_{x,x'}} \mathbf{S} B_2^{-1} ({\bf c} - {\bf a} \, B_1A_2B_2^{-1}) \, {\bf S}^{-1} B_2^{m_{x,x'}} \\
& = B_2^{-m_{x,x'}} B_2^{-1}( {\bf c} - {\bf a} \, q^{-2} B_1 A_1^{-1} A_2) B_2^{m_{x,x'}} \\
& = B_2^{-1}({\bf c} - {\bf a} \, q^{-2} B_1 A_1^{-1} A_2 q^{2m_{x,x'}}),
\end{align*}
yielding, in particular,
\begin{align}
\label{eq:F_conjugation_on_A1}
\mathbf{F}_{\lambda,\lambda'} \, A_1 \, \mathbf{F}_{\lambda,\lambda'}^{-1} 
& = {\bf c}  \, A_1  - {\bf a} \, B_1  A_2 q^{2m_{x,x'}}, \\
\label{eq:F_conjugation_on_A1B1}
\mathbf{F}_{\lambda,\lambda'} \, A_1 B_1^{-1} \, \mathbf{F}_{\lambda,\lambda'}^{-1}  
& = {\bf c} \, A_1 B_1^{-1} B_2^{-1} - {\bf a} \, q^{2m_{x,x'}} A_2 B_2^{-1}.
\end{align}







\vs

Replace $f$ by $A^{-1} f$. Then $\mathbf{F}_{\lambda,\lambda'}^{-1} (A^{-1} f\otimes e_j) = \mathbf{F}^{-1}_{\lambda,\lambda'} A^{-1}_1 (f\otimes e_j)$. So
\begin{align*}
& \textstyle \sum_j (D_{\lambda^*}^{-1}e_j) \otimes \mathbf{F}_{\lambda',(\lambda\lambda')^*}^{-1}( (\mathbf{A} A^{-1} f) \otimes e_j) \\
& \textstyle = \sum_j \mathbf{F}_{\lambda,\lambda'}^{-1}(A^{-1} f\otimes e_j) \otimes (C_{\lambda\lambda'} e_j) \\
& \textstyle = ({\bf c} \cdot \mathrm{id} - {\bf a}\, B_1A_2B_2^{-1}) \, A_1^{-1} \, (\sum_j \mathbf{F}_{\lambda,\lambda'}^{-1}(f\otimes e_j) \otimes (C_{\lambda\lambda'} e_j) ) \\
& \textstyle = ({\bf c} \cdot \mathrm{id} - {\bf a}\, B_1A_2B_2^{-1}) \, A_1^{-1} \, \sum_j (D_{\lambda^*}^{-1} e_j) \otimes \mathbf{F}_{\lambda',(\lambda\lambda')^*}^{-1} ((\mathbf{A} f) \otimes e_j)) \\
& \textstyle = ({\bf c} \, A_1^{-1} - {\bf a}\, B_1 A_1^{-1} A_2B_2^{-1}) \,  (D_{\lambda^*}^{-1})_1(\mathbf{F}_{\lambda',(\lambda\lambda')^*}^{-1} )_{23}  \, \sum_j e_j \otimes (\mathbf{A} f)\otimes e_j .
\end{align*}
We shall move the factors $(D_{\lambda^*}^{-1})_1 (\mathbf{F}_{\lambda',(\lambda\lambda')^*}^{-1} )_{23}$ to the left. From Lem.\ref{lem:conjugation_of_C_D} we have $A_1^{-1} (D_{\lambda^*}^{-1})_1 = (D_{\lambda^*}^{-1})_1 A_1$ and $B_1A_1^{-1} (D_{\lambda^*}^{-1})_1 = (D_{\lambda^*}^{-1})_1 q^{-2m_{-yx^{-1},x}} B_1^{-1}$, in view of $\lambda^* = (x^{-1}, -yx^{-1})$. So
\begin{align*}
& ({\bf c} \, A_1^{-1} - {\bf a}\, B_1 A_1^{-1} A_2B_2^{-1}) \, \ul{ (D_{\lambda^*}^{-1})_1 } \, (\mathbf{F}_{\lambda',(\lambda\lambda')^*}^{-1} )_{23} \\
& = (D_{\lambda^*}^{-1})_1 \, ({\bf c} \, A_1 - {\bf a} \, q^{-2m_{-yx^{-1},x}} \, B_1^{-1} \ul{ A_2 B_2^{-1} } ) \, \ul{ (\mathbf{F}_{\lambda',(\lambda\lambda')^*}^{-1})_{23} } \\
& = (D_{\lambda^*}^{-1})_1 \, (\mathbf{F}_{\lambda',(\lambda\lambda')^*}^{-1})_{23} \, \underbrace{ ({\bf c} \, A_1 - {\bf a} \, q^{-2m_{-yx^{-1},x}} \, B_1^{-1} ({\bf c}_2 \, A_2 B_2^{-1} B_3^{-1} - {\bf a}_2 \, q^{2m_{x',(xx')^{-1}}} \, A_3 B_3^{-1}) )},
\end{align*}
where we used \eqref{eq:F_conjugation_on_A1B1} for $({\bf F}_{\lambda',(\lambda\lambda')^*}^{-1})_{23}$; so
$$
{\bf a}_2 = - \frac{ (y')^{1/N} (xx')^{-1/N} }{(y'(xx')^{-1} - (yx'+y')(xx')^{-1})^{1/N}} = \frac{ (y')^{1/N} (xx')^{-1/N} }{(yx^{-1})^{1/N}}, \quad
{\bf c}_2 = \frac{((yx'+y')(xx')^{-1})^{1/N}}{(yx^{-1})^{1/N}},
$$
in view of $\lambda' = (x',y')$, $(\lambda\lambda')^* = (xx',yx'+y')^* = ( (xx')^{-1}, -(yx'+y')(xx')^{-1})$. Notice that the underbraced part is being applied to $\sum_j e_j \otimes ({\bf A} f) \otimes e_j$, which has `invariance' properties. Applying $B_1^{-1} B_3^{-1}$ to this element yields $\sum_j e_{j-1} \otimes ({\bf A} f)\otimes e_{j-1}$ which equals itself; likewise application of $A_1$ on this element is same as application of $A_3$. So the underbraced part can be replaced by
$$
{\bf c} \, A_3 - {\bf a} \, q^{-2m_{-yx^{-1},x}} \, {\bf c}_2 \, A_2 B_2^{-1} + {\bf a} \, q^{-2m_{-yx^{-1},x}} \, {\bf a}_2 \, q^{2m_{x',(xx')^{-1}}} \, A_3.
$$
I claim that the first and the third terms cancel. Indeed, note
\begin{align*}
& {\bf c}^{-1} \, {\bf a} \, q^{-2m_{-yx^{-1},x}} \, {\bf a}_2 \, q^{2m_{x',(xx')^{-1}}} \\
& = \frac{(yx'+y')^{1/N}}{(y')^{1/N}} \cdot \frac{-y^{1/N} (x')^{1/N}}{(yx'+y')^{1/N}} \cdot q^{-2m_{-yx^{-1},x}} \cdot \frac{(y')^{1/N} (xx')^{-1/N}}{(yx^{-1})^{1/N}} \cdot q^{2m_{x',(xx')^{-1}}} \\
& = x^{1/N} y^{1/N} \cdot  \frac{ q^{2m_{x',(xx')^{-1}}} (x')^{1/N}( (xx')^{-1} )^{1/N}   }{q^{2m_{-yx^{-1},x}} (-yx^{-1})^{1/N} x^{1/N} }
= x^{1/N} y^{1/N} \cdot \frac{ x^{-1/N} }{(-y)^{1/N}} = -1.
\end{align*}
Let's simplify the coefficient of the surviving term:
\begin{align*}
- {\bf a} \, q^{-2m_{-yx^{-1},x}} \, {\bf c}_2
& = \frac{ y^{1/N} (x')^{1/N} }{(yx' + y')^{1/N}} \cdot q^{-2m_{-yx^{-1},x}} \cdot \frac{ ( (yx'+y')(xx')^{-1} )^{1/N} }{(yx^{-1} )^{1/N}} \\
& = - \frac{ y^{1/N} (x')^{1/N} x^{1/N}}{(yx' + y')^{1/N}} \cdot \frac{ q^{2m_{yx'+y', (xx')^{-1}}} (yx'+y')^{1/N} (xx')^{-1/N} }{q^{2m_{-yx^{-1},x}} (-yx^{-1} )^{1/N} x^{1/N}} \\
& = - q^{2m_{yx'+y', (xx')^{-1}}} \frac{ y^{1/N} (x')^{1/N} x^{1/N} (xx')^{-1/N}}{-y^{1/N}}
= q^{2(-m_{x,x'} + m_{yx'+y',(xx')^{-1}})}.
\end{align*}
In the end, we just proved
$$
{\bf A} \, A^{-1} \, f = q^{2(-m_{x,x'} + m_{yx'+y',(xx')^{-1}})} \, AB^{-1} \, \mathbf{A} \, f, \qquad \forall f\in M_{\lambda,\lambda'}^{\lambda\lambda'},
$$
which can be written as the following conjugation action of ${\bf A}$ on $A^{-1}$:
$$
{\bf A} \, A^{-1} \, {\bf A}^{-1} = q^{-2m_1} \, AB^{-1},
$$
where we used $m_1 = m_{x,x'} - m_{yx'+y',(xx')^{-1}}$ from \eqref{eq:m1_and_m2}.

\vs

Similarly, we find out the conjugation action on $B^{-1}$, using the accumulated results. Note
\begin{align*}
& \textstyle \sum_j (D_{\lambda^*}^{-1} e_j) \otimes \mathbf{F}_{\lambda',(\lambda\lambda')^*}^{-1}( ({\bf A} B^{-1} f) \otimes e_j) \\
& \textstyle = \sum_j \mathbf{F}_{\lambda,\lambda'}^{-1} (B^{-1} f \otimes e_j) \otimes (C_{\lambda\lambda'} e_j) \\
& \textstyle = (\mathbf{c}\, B_1^{-1} B_2 - {\bf a} \, A_2) \, \sum_j \mathbf{F}_{\lambda,\lambda'}^{-1} (f\otimes e_j) \otimes (C_{\lambda\lambda'} e_j) \\
& \textstyle = (\mathbf{c}\, B_1^{-1} B_2 - {\bf a} \, A_2) \, \sum_j (D_{\lambda^*}^{-1} e_j) \otimes \mathbf{F}_{\lambda',(\lambda\lambda')^*}^{-1}( ({\bf A} f) \otimes e_j) \\
& \textstyle = (\mathbf{c}\, B_1^{-1} B_2 - {\bf a} \, A_2) \, \ul{ (D_{\lambda^*}^{-1})_1 } \, (\mathbf{F}_{\lambda',(\lambda\lambda')^*}^{-1})_{23} \, \sum_j e_j \otimes ({\bf A} f) \otimes e_j \\
& \textstyle \stackrel{{\rm Lem}.\ref{lem:conjugation_of_C_D}}{=} (D_{\lambda^*}^{-1})_1 \, (\mathbf{c}\, q^{2m_{-yx^{-1},x}} A_1 B_1 \ul{ B_2 } - {\bf a} \, \ul{ A_2} )  \, \ul{ (\mathbf{F}_{\lambda',(\lambda\lambda')^*}^{-1})_{23} } \, \sum_j e_j \otimes ({\bf A} f) \otimes e_j \\
& \stackrel{\eqref{eq:F_conjugation_on_A1B2_and_B1}, \eqref{eq:F_conjugation_on_A1}}{=} (D_{\lambda^*}^{-1})_1 \,(\mathbf{F}_{\lambda',(\lambda\lambda')^*}^{-1})_{23}\,  \\
& \qquad \textstyle \circ \underbrace{ (\mathbf{c}\, q^{2m_{-yx^{-1},x}} A_1 B_1 B_2 B_3 - {\bf a} \, ({\bf c}_2 \, A_2 - {\bf a}_2 \, B_2 A_3 q^{2m_{x',(xx')^{-1}}}) )}  \, \sum_j e_j \otimes ({\bf A} f) \otimes e_j.
\end{align*}
Again, we simplify the underbraced part by the same idea; application of $B_1B_3$ to $\sum_j e_j \otimes ({\bf A} f)\otimes e_j$ is identity application, and application of $A_1$ is same as that of $A_3$:
$$
{\bf c} \, q^{2m_{-yx^{-1},x}} \, B_2 A_3 - {\bf a} \, {\bf c}_2 \, A_2 + {\bf a} \, {\bf a}_2 \, q^{2m_{x',(xx')^{-1}}} \, B_2 A_3.
$$
I claim that the first and the third term cancels; for this we need ${\bf c} \, q^{2m_{-yx^{-1},x}} = - {\bf a} \, {\bf a}_2 \, q^{2m_{x',(xx')^{-1}}}$, which we already showed. For the coefficient of the surviving term, we also use the computed result $- {\bf a} {\bf c}_2 = q^{2(-m_{x,x'} + m_{yx'+y',(xx')^{-1}} + m_{-yx^{-1},x})} = q^{2(-m_1+m_2)}$ (see \eqref{eq:m1_and_m2} for $m_2$). Hence we showed 
$$
{\bf A} \, B^{-1} \, {\bf A}^{-1} = q^{2(-m_1+m_2)} \, A.
$$
Taking the inverses, we have the following conjugation actions of ${\bf A}$ on $A$ and $B$:
\begin{align}
\label{eq:A_conjugation}
 {\bf A} \, A \, {\bf A}^{-1} = q^{2m_1} B A^{-1}, \qquad
{\bf A}\, B \, {\bf A}^{-1} = q^{2(m_1- m_2)} A^{-1},
\end{align}
which, in view of Cor.\ref{cor:conjugation_on_A_B_determines}, determines ${\bf A}$ up to constant, as desired.

\vs

It remains to verify that the map defined in \eqref{eq:A_formula} satisfies the conjugation equations \eqref{eq:A_conjugation}. For convenience, denote the map defined in \eqref{eq:A_formula} by $\til{\bf A}$. Then
\begin{align*}
\til{\bf A} A e_i & \textstyle = \til{\bf A} q^{2i} e_i = q^{2i} \sum_j q^{-2ij - j^2 + (2m_1+1)j + 2(m_1-m_2)i} e_j, \\
q^{2m_1} BA^{-1} {\bf A} e_i & =\textstyle  q^{2m_1} BA^{-1} \sum_j  q^{-2ij - j^2 +(2m_1+1)j + 2(m_1-m_2) i} e_j \\
& \textstyle = q^{2m_1} \sum_j q^{-2ij - j^2 +(2m_1+1)j + 2(m_1-m_2) i} q^{-2j} e_{j+1} \\
& \textstyle = q^{2m_1} \sum_j q^{-2i(j-1) - (j-1)^2 + (2m_1+1)(j-1) + 2(m_1-m_2)i} q^{-2(j-1)} e_j,
\end{align*}
and one can verify that $2i -2ij - j^2 + (2m_1+1)j + 2(m_1-m_2)i$ coincides with $2m_1 - 2i(j-1) - (j-1)^2 + (2m_1 + 1)(j-1) + 2(m_1-m_2)i-2(j-1)$ for each $i,j$. Likewise, note
\begin{align*}
 \til{\bf A} B e_i & = \til{\bf A} e_{i+1} = \textstyle \sum_j q^{-2(i+1)j - j^2 +(2m_1+1)j + 2(m_1-m_2) (i+1)} e_j, \\
q^{2(m_1-m_2)} A^{-1} {\bf A} e_i & = q^{2(m_1-m_2)} A^{-1} \textstyle \sum_j q^{-2ij - j^2 +(2m_1+1)j + 2(m_1-m_2) i} e_j \\
& = q^{2(m_1-m_2)} \textstyle \sum_j q^{-2ij-j^2 +(2m_1+1)j + 2(m_1-m_2)i} q^{-2j} e_j,
\end{align*}
and one can check $-2(i+1)j - j^2 + (2m_1+1)j + 2(m_1-m_2)(i+1)$ coincides with $2(m_1-m_2) - 2ij -j^2 + (2m_1+1)j + 2(m_1-m_2)i-2j$ for each $i,j$. \qed

\section{Representation of Kashaev groupoid}


\subsection{Graphical encoding}
\label{subsec:graphical_encoding}

To proceed further with the representation theory with the help of the operators ${\bf T}$ and ${\bf A}$, we pause for a moment and discuss how to encode the results obtained so far by pictures. First, the multiplicity space $M_{\lambda,\lambda'}^{\lambda\lambda'}$, which is identified with the intertwiner space $\mathrm{Hom}_\mathcal{W}(V_{\lambda\lambda'}, V_\lambda \otimes V_{\lambda'})$, is encoded as the following triangle.
\begin{definition}
A \ul{\em labeled dotted triangle}, or \ul{\em LD-triangle}, is a triangle on an oriented surface with a distinguished corner depicted by a dot $\bullet$, with all three edges labeled by non-singular weights.
\end{definition}

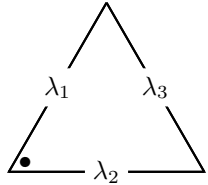
\begin{figure}[htbp!]
\centering
\begin{pspicture}[showgrid=false](0,0)(3.3,2.8)
\rput[bl](0,-0.5){
\PstTriangle[unit=1.5,PolyName=P]
\pcline(P1)(P2)\ncput*{$\lambda_1$}
\pcline(P2)(P3)\ncput*{$\lambda_2$}
\pcline(P3)(P1)\ncput*{$\lambda_3$}
\rput[l]{30}(P2){\hspace{1,7mm}$\bullet$}
}
\end{pspicture}
\caption{A triangle representing $M_{\lambda_1,\lambda_2}^{\lambda_3} \cong \mathrm{Hom}_\mathcal{W}(V_{\lambda_3}, V_{\lambda_1}\otimes V_{\lambda_2})$}
\label{fig:one_triangle}
\end{figure}

\begin{definition}
\label{def:sane}
Suppose that the edge labels of a LD-triangle are $\lambda_1,\lambda_2,\lambda_3$, with this counterclockwise order, with $\lambda_3$ being the edge at the opposite of the dot. This LD-triangle is said to be \ul{\em sane} if $\lambda_1 \lambda_2 = \lambda_3$.
\end{definition}
Such a sane LD-triangle encodes the situation $M_{\lambda_1,\lambda_2}^{\lambda_3} \cong \mathrm{Hom}_\mathcal{W}(V_{\lambda_3}, V_{\lambda_1}\otimes V_{\lambda_2})$, namely, `multiply $\lambda_1$ and $\lambda_2$ and decompose by $\lambda_3$'; see Fig.\ref{fig:one_triangle}. We allow to glue LD-triangles along the edges with same label. For several LD-triangles glued together, we assign the tensor product of the corresponding $M$ spaces, or $\mathrm{Hom}_\mathcal{W}$ spaces.

\vspace{-5mm}

\vspace{-5mm}

\begin{figure}[htbp!]
\centering
\begin{pspicture}[showgrid=false](0,0.5)(8.4,3.5)
\rput[bl](0.5,0){
\PstTriangle[unit=1.5,PolyName=P]
\pcline(P1)(P2)\ncput*{$\lambda_1$}
\pcline(P2)(P3)\ncput*{$\lambda_2$}
\pcline(P3)(P1)\ncput*{$\lambda_3$}
\rput[l]{30}(P2){\hspace{1,7mm}$\bullet$}
}
\rput[l](2.0,1.5){$t$}
\rput[bl](4.5,0){
\PstTriangle[unit=1.5,PolyName=P]
\pcline(P1)(P2)\ncput*{$\lambda_1^*$}
\pcline(P2)(P3)\ncput*{$\lambda_2$}
\pcline(P3)(P1)\ncput*{$\lambda_3^*$}
\rput[l]{150}(P3){\hspace{1,7mm}$\bullet$}
}
\rput[l](6,1.5){$t$}
\rput[l](3.7,2){\pcline{->}(0,0)(1;0)\Aput{${\bf A}_t$}}
\end{pspicture}
\caption{The move representing ${\bf A}_t = {\bf A} : M_{\lambda_1,\lambda_2}^{\lambda_3} \to M_{\lambda_2, \lambda_3^*}^{\lambda_1^*}$}
\label{fig:A}
\end{figure}
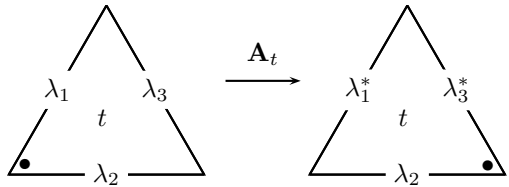

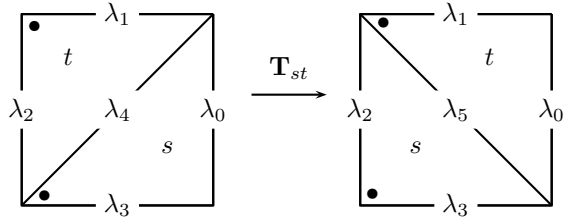
\begin{figure}[htbp!]
\centering
\begin{pspicture}[showgrid=false](0,0.5)(8.4,3.5)
\rput[bl](0;0){
\PstSquare[unit=1.8,PolyName=P]
\pcline(P1)(P2)\ncput*{$\lambda_1$}
\pcline(P2)(P3)\ncput*{$\lambda_2$}
\pcline(P3)(P4)\ncput*{$\lambda_3$}
\pcline(P4)(P1)\ncput*{$\lambda_0$}
\pcline(P1)(P3)\ncput*{$\lambda_4$}
\rput[l]{-45}(P2){\hspace{1,4mm}$\bullet$}
\rput[l]{22}(P3){\hspace{2,4mm}$\bullet$}
}
\rput[l](2.5,1.3){$s$}
\rput[l](1.2,2.5){$t$}
\rput[bl](4.5;0){
\PstSquare[unit=1.8,PolyName=P]
\pcline(P1)(P2)\ncput*{$\lambda_1$}
\pcline(P2)(P3)\ncput*{$\lambda_2$}
\pcline(P3)(P4)\ncput*{$\lambda_3$}
\pcline(P4)(P1)\ncput*{$\lambda_0$}
\pcline(P2)(P4)\ncput*{$\lambda_5$}
\rput[l]{-23}(P2){\hspace{2,4mm}$\bullet$}
\rput[l]{41}(P3){\hspace{1,4mm}$\bullet$}
}
\rput[l](5.8,1.3){$s$}
\rput[l](6.8,2.5){$t$}
\rput[l](3.7,2){\pcline{->}(0,0)(1;0)\Aput{${\bf T}_{st}$}}
\end{pspicture}
\caption{The move representing ${\bf T}_{st} = {\bf T}_{12} : M_{\lambda_4,\lambda_3}^{\lambda_0} \otimes M_{\lambda_1,\lambda_2}^{\lambda_4} \to M_{\lambda_2,\lambda_3}^{\lambda_5} \otimes M_{\lambda_1,\lambda_5}^{\lambda_0}$}
\label{fig:T}
\end{figure}

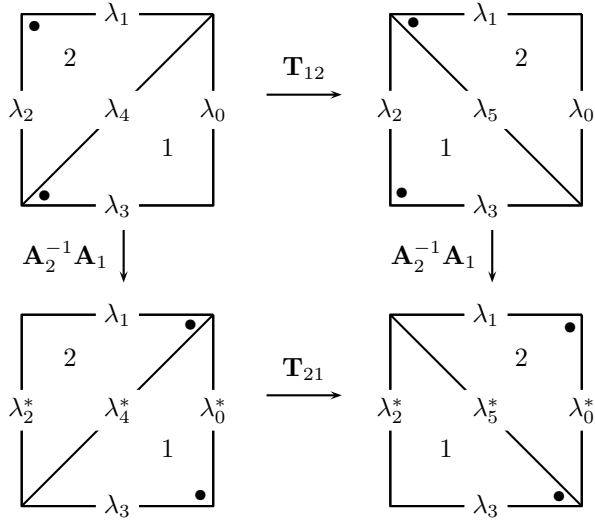
\begin{figure}[htbp!]
\centering
\begin{pspicture}[showgrid=false](0,-3.4)(8.4,3.1)
\rput[bl](0;0){
\PstSquare[unit=1.8,PolyName=P]
\pcline(P1)(P2)\ncput*{$\lambda_1$}
\pcline(P2)(P3)\ncput*{$\lambda_2$}
\pcline(P3)(P4)\ncput*{$\lambda_3$}
\pcline(P4)(P1)\ncput*{$\lambda_0$}
\pcline(P1)(P3)\ncput*{$\lambda_4$}
\rput[l]{-45}(P2){\hspace{1,4mm}$\bullet$}
\rput[l]{22}(P3){\hspace{2,4mm}$\bullet$}
}
\rput[l](2.5,1.3){$1$}
\rput[l](1.2,2.5){$2$}
\rput[bl](0,-4){
\PstSquare[unit=1.8,PolyName=P]
\pcline(P1)(P2)\ncput*{$\lambda_1$}
\pcline(P2)(P3)\ncput*{$\lambda_2^*$}
\pcline(P3)(P4)\ncput*{$\lambda_3$}
\pcline(P4)(P1)\ncput*{$\lambda_0^*$}
\pcline(P1)(P3)\ncput*{$\lambda_4^*$}
\rput[l]{202}(P1){\hspace{2,4mm}$\bullet$}
\rput[l]{135}(P4){\hspace{1,4mm}$\bullet$}
}
\rput[l](2.5,-2.7){$1$}
\rput[l](1.2,-1.5){$2$}
\rput[bl](4.9,0){
\PstSquare[unit=1.8,PolyName=P]
\pcline(P1)(P2)\ncput*{$\lambda_1$}
\pcline(P2)(P3)\ncput*{$\lambda_2$}
\pcline(P3)(P4)\ncput*{$\lambda_3$}
\pcline(P4)(P1)\ncput*{$\lambda_0$}
\pcline(P2)(P4)\ncput*{$\lambda_5$}
\rput[l]{-22}(P2){\hspace{2,4mm}$\bullet$}
\rput[l]{45}(P3){\hspace{1,4mm}$\bullet$}
}
\rput[l](6.2,1.3){$1$}
\rput[l](7.2,2.5){$2$}
\rput[bl](4.9,-4){
\PstSquare[unit=1.8,PolyName=P]
\pcline(P1)(P2)\ncput*{$\lambda_1$}
\pcline(P2)(P3)\ncput*{$\lambda_2^*$}
\pcline(P3)(P4)\ncput*{$\lambda_3$}
\pcline(P4)(P1)\ncput*{$\lambda_0^*$}
\pcline(P2)(P4)\ncput*{$\lambda_5^*$}
\rput[l]{-135}(P1){\hspace{1,4mm}$\bullet$}
\rput[l]{155}(P4){\hspace{2,4mm}$\bullet$}
}
\rput[l](6.2,-2.7){$1$}
\rput[l](7.2,-1.5){$2$}
\rput[l](3.9,2){\pcline{->}(0,0)(1;0)\Aput{${\bf T}_{12}$}}
\rput[l](2,0.2){\pcline{->}(0,0)(0,-0.7)\Bput{${\bf A}_2^{-1} {\bf A}_1$}}
\rput[l](3.9,-2){\pcline{->}(0,0)(1;0)\Aput{${\bf T}_{21}$}}
\rput[l](6.9,0.2){\pcline{->}(0,0)(0,-0.7)\Bput{${\bf A}_2^{-1} {\bf A}_1$}}
\end{pspicture}
\caption{The relation ${\bf A}_1 {\bf T}_{12} {\bf A}_2 = {\bf A}_2 {\bf T}_{21} {\bf A}_1$}
\label{fig:ATA_ATA}
\end{figure}

\begin{figure}[htbp!]
\centering
\begin{pspicture}[showgrid=false](0,-3.5)(8.4,3.1)
\rput[bl](0;0){
\PstSquare[unit=1.8,PolyName=P]
\pcline(P1)(P2)\ncput*{$\lambda_1$}
\pcline(P2)(P3)\ncput*{$\lambda_2$}
\pcline(P3)(P4)\ncput*{$\lambda_3$}
\pcline(P4)(P1)\ncput*{$\lambda_0$}
\pcline(P1)(P3)\ncput*{$\lambda_4$}
\rput[l]{-45}(P2){\hspace{1,4mm}$\bullet$}
\rput[l]{22}(P3){\hspace{2,4mm}$\bullet$}
}
\rput[l](2.5,1.3){$2$}
\rput[l](1.2,2.5){$1$}
\rput[bl](0,-4){
\PstSquare[unit=1.8,PolyName=P]
\pcline(P1)(P2)\ncput*{$\lambda_1^*$}
\pcline(P2)(P3)\ncput*{$\lambda_2$}
\pcline(P3)(P4)\ncput*{$\lambda_3$}
\pcline(P4)(P1)\ncput*{$\lambda_0^*$}
\pcline(P1)(P3)\ncput*{$\lambda_4^*$}
\rput[l]{68}(P3){\hspace{2,4mm}$\bullet$}
\rput[l]{135}(P4){\hspace{1,4mm}$\bullet$}
}
\rput[l](2.5,-2.7){$1$}
\rput[l](1.2,-1.5){$2$}
\rput[bl](4.9,0){
\PstSquare[unit=1.8,PolyName=P]
\pcline(P1)(P2)\ncput*{$\lambda_1$}
\pcline(P2)(P3)\ncput*{$\lambda_2$}
\pcline(P3)(P4)\ncput*{$\lambda_3$}
\pcline(P4)(P1)\ncput*{$\lambda_0$}
\pcline(P2)(P4)\ncput*{$\lambda_5$}
\rput[l]{-22}(P2){\hspace{2,4mm}$\bullet$}
\rput[l]{45}(P3){\hspace{1,4mm}$\bullet$}
}
\rput[l](6.2,1.3){$2$}
\rput[l](7.2,2.5){$1$}
\rput[bl](4.9,-4){
\PstSquare[unit=1.8,PolyName=P]
\pcline(P1)(P2)\ncput*{$\lambda_1^*$}
\pcline(P2)(P3)\ncput*{$\lambda_2$}
\pcline(P3)(P4)\ncput*{$\lambda_3$}
\pcline(P4)(P1)\ncput*{$\lambda_0^*$}
\pcline(P2)(P4)\ncput*{$\lambda_5$}
\rput[l]{112}(P4){\hspace{2,4mm}$\bullet$}
\rput[l]{45}(P3){\hspace{1,4mm}$\bullet$}
}
\rput[l](6.2,-2.7){$2$}
\rput[l](7.2,-1.5){$1$}
\rput[l](3.9,2){\pcline{->}(0,0)(1;0)\Aput{${\bf T}_{21}$}}
\rput[l](2,0.2){\pcline{->}(0,0)(0,-0.7)\Bput{${\bf A}_1 {\bf A}_2 {\bf P}_{(12)}$}}
\rput[l](3.9,-2){\pcline{->}(1;0)(0,0)\Aput{${\bf T}_{12}$}}
\rput[r](6.9,0.2){\pcline{->}(0,0)(0,-0.7)\Bput{${\bf A}_1$}}
\end{pspicture}
\caption{The relation ${\bf T}_{12} {\bf A}_1 {\bf T}_{21} = {\bf A}_1 {\bf A}_2 {\bf P}_{(12)}$}
\label{fig:TAT_AAP}
\end{figure}
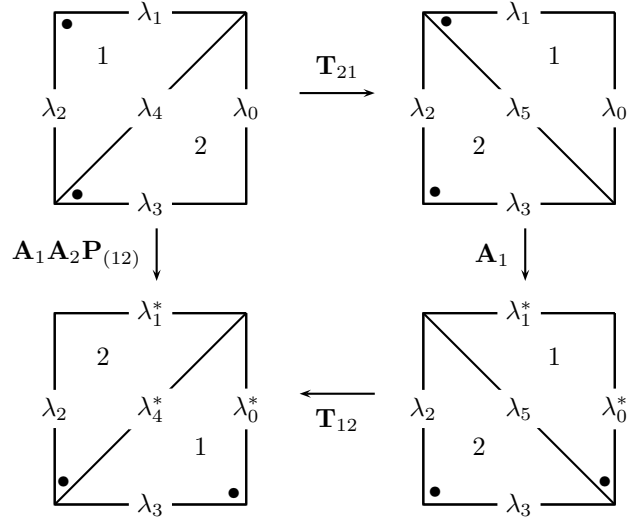

Notice that the ${\bf A}$ operator can easily be understood as just `moving the dot of an LD-triangle, with appropriate change of edge-labels', and the ${\bf T}$ operator as the `flip of a quadrilateral, with appropriate change of edge-label of the inner edge'; see Fig.\ref{fig:A} and \ref{fig:T} \footnote{I note here that the source codes of Figures \ref{fig:one_triangle}--\ref{fig:TAT_AAP} are taken from \cite{FK12}.}. These pictorial counterpart of ${\bf A}$ and ${\bf T}$ operators satisfy certain consistency relations. For example, do the picture version of ${\bf A}$ three times, then we are back to the same picture. And when three LD-triangles are glued, after a certain sequence of five flips of inner edges we get back to the original picture. These two relations are indeed satisfied by our operators ${\bf A}$ and ${\bf T}$, as we showed in Propositions \ref{prop:TT_TTT} and \ref{prop:AAA}. The two diagrams of pictures, in Fig.\ref{fig:ATA_ATA} and \ref{fig:TAT_AAP} suggest two more kinds.

Indeed, in the following subsections, we show that these relations are satisfied by our operators ${\bf A}$ and ${\bf T}$. As shall be pointed out again later, these are {\em all} possible consistency relations to be checked.

\subsection{First consistency relation involving ${\bf T}$ and ${\bf A}$}

For any regular triple $(\lambda,\lambda',\lambda'')$ of weights, recall the ${\bf T}$ map:
$$
{\bf T}_{\lambda,\lambda',\lambda''} : M_{\lambda\lambda',\lambda''}^{\lambda\lambda'\lambda''} \otimes M_{\lambda,\lambda'}^{\lambda\lambda'} \to M_{\lambda',\lambda''}^{\lambda'\lambda''} \otimes M_{\lambda,\lambda'\lambda''}^{\lambda\lambda'\lambda''}
$$
We put
\begin{align}
\label{eq:lambda_1_to_5}
\lambda_1 = \lambda, \quad
\lambda_2 = \lambda', \quad
\lambda_3 = \lambda'', \quad
\lambda_0 = \lambda\lambda'\lambda'', \quad
\lambda_4 = \lambda\lambda', \quad
\lambda_5 = \lambda'\lambda'',
\end{align}
in accordance with the upper row of the diagram in Fig.\ref{fig:ATA_ATA}; write
$$
\lambda_i = (x_i,y_i), \quad \forall i=0,1,\ldots,5.
$$ 
We assume all five $\lambda_i$'s are non-singular; from the relations among them, we notice that $\lambda_1,\lambda_2,\lambda_3$ determine $\lambda_4,\lambda_5$. In terms of $\mathrm{Hom}_\mathcal{W}$ spaces the above $\mathbf{T}$ map translates first to the following version, via the identification maps in Def.\ref{def:I}
\begin{align*}
\mathbf{T}^\mathrm{Hom}_{\lambda_1,\lambda_2,\lambda_3} & : \mathrm{Hom}_\mathcal{W}(V_{\lambda_0}, V_{\lambda_4} \otimes V_{\lambda_3}) \otimes \mathrm{Hom}_\mathcal{W}(V_{\lambda_4}, V_{\lambda_1} \otimes V_{\lambda_2}) \\
& \quad \to
\mathrm{Hom}_\mathcal{W}(V_{\lambda_5}, V_{\lambda_2} \otimes V_{\lambda_3}) \otimes \mathrm{Hom}_\mathcal{W}(V_{\lambda_0}, V_{\lambda_1} \otimes V_{\lambda_5}),
\end{align*}
which in turn translates via Lem.\ref{lem:composition_lemma} to the following identity map
$$
\mathrm{id} : \mathrm{Hom}_\mathcal{W}(V_{\lambda_0}, (V_{\lambda_1} \otimes V_{\lambda_2}) \otimes V_{\lambda_3}) \to \mathrm{Hom}_\mathcal{W}(V_{\lambda_0}, V_{\lambda_1} \otimes (V_{\lambda_2} \otimes V_{\lambda_3})).
$$

\begin{lemma}[Hom version of ${\bf T}$ map]
\label{lem:T_Hom}
Pick any elements $\sum f_1 \otimes f_2 \in \mathrm{Hom}_\mathcal{W}(V_{\lambda_0}, V_{\lambda_4} \otimes V_{\lambda_3}) \otimes \mathrm{Hom}_\mathcal{W}(V_{\lambda_4}, V_{\lambda_1} \otimes V_{\lambda_2})$ and $\sum f_3 \otimes f_4 \in \mathrm{Hom}_\mathcal{W}(V_{\lambda_5}, V_{\lambda_2} \otimes V_{\lambda_3}) \otimes \mathrm{Hom}_\mathcal{W}(V_{\lambda_0}, V_{\lambda_1} \otimes V_{\lambda_5})$. Then $\mathbf{T}^\mathrm{Hom}(\sum f_1 \otimes f_2) = \sum f_3 \otimes f_4$ if and only if $\sum (f_2\otimes \mathrm{id})\circ f_1 = \sum (\mathrm{id}\otimes f_3) \circ f_4$, as elements of $\mathrm{Hom}_\mathcal{W}(V_{\lambda_0}, V_{\lambda_1}\otimes V_{\lambda_2} \otimes V_{\lambda_3})$. \qed
\end{lemma}

The first consistency relation involving both ${\bf T}$ and ${\bf A}$ and corresponding to the diagram in Fig.\ref{fig:ATA_ATA} is formulated and shown as follows. Instead of a computational proof using the explicit formulas of ${\bf T}$ and ${\bf A}$ obtained in \S\ref{sec:computation_of_formulas}, which would be very complicated and brute force, we give a clean and more enlightening representation theoretic proof.
\begin{proposition}[${\bf ATA} = {\bf ATA}$ relation]
\label{prop:ATA_ATA}
Let $\lambda_0,\lambda_1,\lambda_2,\lambda_3,\lambda_4,\lambda_5$ be as above. Consider the map $\mathbf{T}_{\lambda,\lambda',\lambda''}$ written as
$$
\mathbf{T}_{\lambda_1,\lambda_2,\lambda_3} : M_{\lambda_4,\lambda_3}^{\lambda_0} \otimes M_{\lambda_1,\lambda_2}^{\lambda_4} \to M_{\lambda_2,\lambda_3}^{\lambda_5} \otimes M_{\lambda_1, \lambda_5}^{\lambda_0},
$$
and another ${\bf T}$ map
$$
\mathbf{T}_{\lambda_3, \lambda_0^*, \lambda_1} : M_{\lambda_4^*, \lambda_1}^{\lambda_2^*} \otimes M_{\lambda_3, \lambda_0^*}^{\lambda_4^*} \to M_{\lambda_0^*, \lambda_1}^{\lambda_5^*} \otimes M_{\lambda_3, \lambda_5^*}^{\lambda_2^*}.
$$
Then this map ${\bf T}_{\lambda_3,\lambda_0^*,\lambda_1}$ makes sense as a ${\bf T}$ map, and one has
\begin{align}
\label{eq:TAA_AAT}
({\bf T}_{\lambda_3, \lambda_0^*, \lambda_1})_{21} \, ({\bf A}_{\lambda_1, \lambda_4^*}^{\lambda_2^*})_2^{-1} \, ({\bf A}_{\lambda_4, \lambda_3}^{\lambda_0})_1 = q^{2m}\, ({\bf A}_{\lambda_1,\lambda_0^*}^{\lambda_5^*})_2^{-1} \, ({\bf A}_{\lambda_2,\lambda_3}^{\lambda_5})_1 \, (\mathbf{T}_{\lambda_1,\lambda_2,\lambda_3})_{12},
\end{align}
as maps $M_{\lambda_4,\lambda_3}^{\lambda_0} \otimes M_{\lambda_1,\lambda_2}^{\lambda_4} \to M_{\lambda_3,\lambda_5^*}^{\lambda_2^*} \otimes M_{\lambda_0^*,\lambda_1}^{\lambda_5^*}$, where
$$
m = m_{x_1,x_2} - m_{x_1,x_2x_3} = m_{x,x'} - m_{x,x'x''}.
$$
\end{proposition}
One can symbolically write down the above relation \eqref{eq:TAA_AAT} as ${\bf A}_1 {\bf T}_{12} {\bf A}_2 = q^{-2m} {\bf A}_2 {\bf T}_{21} {\bf A}_1$.

\vs

{\it Proof.} We prove the Hom version. The Hom versions of the two ${\bf T}$ maps are
\begin{align*}
(\mathbf{T}_{\lambda_1,\lambda_2,\lambda_3}^\mathrm{Hom})_{12} & : \mathrm{Hom}_\mathcal{W}(V_{\lambda_0},V_{\lambda_4}\otimes V_{\lambda_3}) \otimes \mathrm{Hom}_\mathcal{W}(V_{\lambda_4}, V_{\lambda_1} \otimes V_{\lambda_2}) \\
& \quad \to \mathrm{Hom}_\mathcal{W}(V_{\lambda_5}, V_{\lambda_2} \otimes V_{\lambda_3}) \otimes \mathrm{Hom}_\mathcal{W}(V_{\lambda_0}, V_{\lambda_1} \otimes V_{\lambda_5}), \\
(\mathbf{T}_{\lambda_3, \lambda_0^*, \lambda_1}^\mathrm{Hom})_{21} & : \mathrm{Hom}_\mathcal{W}(V_{\lambda_4^*}, V_{\lambda_3} \otimes V_{\lambda_0^*}) \otimes \mathrm{Hom}_\mathcal{W}(V_{\lambda_2^*}, V_{\lambda_4^*} \otimes V_{\lambda_1}) \\
& \quad \to \mathrm{Hom}_\mathcal{W}(V_{\lambda_2^*}, V_{\lambda_3} \otimes V_{\lambda_5^*}) \otimes \mathrm{Hom}_\mathcal{W}(V_{\lambda_5^*}, V_{\lambda_0^*} \otimes V_{\lambda_1}).
\end{align*}
Pick any element $\sum f_1 \otimes f_2 \in \mathrm{Hom}_\mathcal{W}(V_{\lambda_0}, V_{\lambda_4} \otimes V_{\lambda_3}) \otimes \mathrm{Hom}_\mathcal{W}(V_{\lambda_4}, V_{\lambda_1} \otimes V_{\lambda_2})$. Then, as in the previous lemma, the element $\mathbf{T}^\mathrm{Hom}_{12}(\sum f_1 \otimes f_2) = \sum f_3 \otimes f_4$ of $\mathrm{Hom}_\mathcal{W}(V_{\lambda_5}, V_{\lambda_2} \otimes V_{\lambda_3}) \otimes \mathrm{Hom}_\mathcal{W}(V_{\lambda_0}, V_{\lambda_1} \otimes V_{\lambda_5})$ can be described completely by the equation $\sum (f_2\otimes \mathrm{id})\circ f_1 = \sum (\mathrm{id}\otimes f_3) \circ f_4$ of elements in $\mathrm{Hom}_\mathcal{W}(V_{\lambda_0}, V_{\lambda_1}\otimes V_{\lambda_2} \otimes V_{\lambda_3})$.

\vs

We also use the Hom version \eqref{eq:A_as_Hom} of the ${\bf A}$ operators. Understood appropriately, we would like to show
$$
q^{2m} \, \sum {\bf A}^\mathrm{Hom} f_3 \otimes ({\bf A}^\mathrm{Hom})^{-1} f_4 = \mathbf{T}^\mathrm{Hom}_{21} (\sum {\bf A}^\mathrm{Hom} f_1 \otimes (\mathbf{A}^\mathrm{Hom})^{-1} f_2),
$$
which is equivalent to $\sum ({\bf A}^\mathrm{Hom} f_1 \otimes \mathrm{id}) \circ (({\bf A}^\mathrm{Hom})^{-1} f_2) = q^{2m} \, \sum (\mathrm{id}\otimes (\mathbf{A}^\mathrm{Hom})^{-1} f_4) \circ ({\bf A}^\mathrm{Hom} f_3)$, which is an equality of elements of $\mathrm{Hom}_\mathcal{W}(V_{\lambda_2^*}, V_{\lambda_3} \otimes V_{\lambda_0^*} \otimes V_{\lambda_1})$.

\vs

We now go to Inv versions. Define the following elements via the left and right canonical maps \eqref{eq:Hom_as_tensor_products}:
\begin{align*}
& \til{f}_1 := J^\mathrm{L} f_1 \in \mathrm{Inv}(V_{\lambda_4} \otimes V_{\lambda_3} \otimes V_{\lambda_0}^*), \qquad
\til{f}_2 := J^\mathrm{R} f_2 \in \mathrm{Inv}({}^* V_{\lambda_4}\otimes V_{\lambda_1} \otimes V_{\lambda_2}), \\
& \til{f}_3 := J^\mathrm{L} f_3 \in \mathrm{Inv}(V_{\lambda_2} \otimes V_{\lambda_3} \otimes V_{\lambda_5}^*), \qquad
\til{f}_4 := J^\mathrm{R} f_4 \in \mathrm{Inv}({}^* V_{\lambda_0} \otimes V_{\lambda_1} \otimes V_{\lambda_5}), \\
& \til{f}_5 := J^\mathrm{R} \mathbf{A}^\mathrm{Hom} f_1 \in \mathrm{Inv}({}^* V_{\lambda_4^*} \otimes V_{\lambda_3} \otimes V_{\lambda_0^*}), \quad
\til{f}_6 := J^\mathrm{L} (\mathbf{A}^\mathrm{Hom})^{-1} f_2 \in \mathrm{Inv}(V_{\lambda_4^*} \otimes V_{\lambda_1} \otimes V_{\lambda_2^*}^*), \\
& \til{f}_7 := J^\mathrm{R} \mathbf{A}^\mathrm{Hom}f_3 \in \mathrm{Inv}({}^* V_{\lambda_2^*} \otimes V_{\lambda_3} \otimes V_{\lambda_5^*}), \quad
\til{f}_8 := J^\mathrm{L} (\mathbf{A}^\mathrm{Hom})^{-1} f_4 \in \mathrm{Inv}(V_{\lambda_0^*} \otimes V_{\lambda_1} \otimes V_{\lambda_5^*}^*).
\end{align*}
Definition \eqref{eq:A_as_Hom} of the $\mathbf{A}^\mathrm{Hom}$ map says
\begin{align*}
\til{f}_5 = (D_{\lambda_4^*})_1 (C_{\lambda_0})_3 \til{f}_1, \quad
\til{f}_6 = (D_{\lambda_4})_1^{-1} (C_{\lambda_2^*})_3^{-1} \til{f}_2, \quad
\til{f}_7 = (D_{\lambda_2^*})_1 (C_{\lambda_5})_3 \til{f}_3, \quad
\til{f}_8 = (D_{\lambda_0})_1^{-1} (C_{\lambda_5^*})_3^{-1} \til{f}_4.
\end{align*}

\vs

Let us now make use of the condition $\sum (f_2 \otimes \mathrm{id}) \circ f_1 = \sum (\mathrm{id} \otimes f_3) \circ f_4$. Note that $\til{f}_2 \otimes \til{f}_1 \in {}^*V_{\lambda_4} \otimes V_{\lambda_1} \otimes V_{\lambda_2} \otimes V_{\lambda_4} \otimes V_{\lambda_3} \otimes V_{\lambda_0}^*$, and the symbol $\circ$ in $(f_2 \otimes \mathrm{id})\circ f_1$ means to pair out ${}^*V_{\lambda_4}$ and $V_{\lambda_4}$; these are $1$st and $4$th tensor factors, and we denote this pairing out map of these factors by $\mathrm{ev}_{14}$, so that $\sum \mathrm{ev}_{14}(\til{f}_2 \otimes \til{f}_1) \in V_{\lambda_1} \otimes V_{\lambda_2} \otimes V_{\lambda_3} \otimes V_{\lambda_0}^*$, which represents the element $\sum (f_2 \otimes \mathrm{id}) \circ f_1$ of $\mathrm{Hom}_\mathbb{C}(V_{\lambda_0}, V_{\lambda_1}\otimes V_{\lambda_2} \otimes V_{\lambda_3})$ via the canonical map $J^\mathrm{L}$ of \eqref{eq:Hom_as_tensor_products}. In general, we denote by $\mathrm{ev}_{ij}$ the map that `pairs out' the $i$-th and $j$-th factors; so the codomain has $2$ tensor factors less than the domain.  Likewise, $\til{f}_4 \otimes \til{f}_3 \in {}^*V_{\lambda_0} \otimes V_{\lambda_1} \otimes V_{\lambda_5} \otimes V_{\lambda_2} \otimes V_{\lambda_3} \otimes V^*_{\lambda_5}$, so that $\sum \mathrm{ev}_{36} (\til{f}_4 \otimes \til{f}_3) \in {}^*V_{\lambda_0} \otimes V_{\lambda_1} \otimes V_{\lambda_2} \otimes V_{\lambda_3}$ represents the element $\sum (\mathrm{id}\otimes f_3) \circ f_4$ of $\mathrm{Hom}_\mathbb{C}(V_{\lambda_0}, V_{\lambda_1}\otimes V_{\lambda_2} \otimes V_{\lambda_3})$ via the canonical map $J^\mathrm{R}$ of \eqref{eq:Hom_as_tensor_products}. In view of Lem.\ref{lem:composition_of_canonical_maps}, the equality $\sum(f_2 \otimes \mathrm{id}) \circ f_1 = \sum(\mathrm{id}\otimes f_3) \circ f_4$ translates to the following equality of elements in the vector space $V_{\lambda_1} \otimes V_{\lambda_2} \otimes V_{\lambda_3} \otimes V_{\lambda_0}^*$:
\begin{align}
\label{eq:ATA_ATA_to_assume}
 \textstyle \sum \mathrm{ev}_{14}(\til{f}_2 \otimes \til{f}_1) = \mathbf{P}_{(1432)} \sum \mathrm{ev}_{36} (\til{f}_4 \otimes \til{f}_3).
\end{align}

\vs

We now move on to the equation $\sum ({\bf A}^\mathrm{Hom} f_1 \otimes \mathrm{id}) \circ (({\bf A}^\mathrm{Hom})^{-1} f_2) = q^{2m} \, \sum (\mathrm{id}\otimes (\mathbf{A}^\mathrm{Hom})^{-1} f_4) \circ ({\bf A}^\mathrm{Hom} f_3)$ which we would like to prove. Note $\til{f}_5 \otimes \til{f}_6 \in {}^*V_{\lambda_4^*} \otimes V_{\lambda_3} \otimes V_{\lambda_0^*} \otimes V_{\lambda_4^*} \otimes V_{\lambda_1} \otimes V_{\lambda_2^*}^*$, so that $\sum \mathrm{ev}_{14} (\til{f}_5 \otimes \til{f}_6) \in V_{\lambda_3} \otimes V_{\lambda_0^*} \otimes V_{\lambda_1} \otimes V_{\lambda_2^*}^*$ represents the element $\sum ({\bf A}^\mathrm{Hom} f_1 \otimes \mathrm{id}) \circ (({\bf A}^\mathrm{Hom})^{-1} f_2)$ of $\mathrm{Hom}_\mathbb{C}(V_{\lambda_2^*}, V_{\lambda_3} \otimes V_{\lambda_0^*} \otimes V_{\lambda_1})$. Likewise, $\til{f}_7 \otimes \til{f}_8 \in {}^* V_{\lambda_2^*} \otimes V_{\lambda_3} \otimes V_{\lambda_5^*} \otimes V_{\lambda_0^*} \otimes V_{\lambda_1} \otimes V_{\lambda_5^*}^*$, so that $\sum \mathrm{ev}_{36}(\til{f}_7 \otimes \til{f}_8) \in {}^*V_{\lambda_2^*} \otimes V_{\lambda_3} \otimes V_{\lambda_0^*} \otimes V_{\lambda_1}$ represents the element $\sum (\mathrm{id}\otimes (\mathbf{A}^\mathrm{Hom})^{-1} f_4) \circ ({\bf A}^\mathrm{Hom} f_3)$ of $\mathrm{Hom}_\mathbb{C}(V_{\lambda_2^*}, V_{\lambda_3} \otimes V_{\lambda_0^*} \otimes V_{\lambda_1})$. Thus, the equation that we would like to prove is equivalent to the following equation of elements in $V_{\lambda_3} \otimes V_{\lambda_0^*} \otimes V_{\lambda_1} \otimes V_{\lambda_2^*}^*$:
\begin{align}
\label{eq:ATA_ATA_to_prove}
\textstyle \sum \mathrm{ev}_{14}(\til{f}_5 \otimes \til{f}_6) = q^{2m} \, \mathbf{P}_{(1432)} \sum \mathrm{ev}_{36}(\til{f}_7 \otimes \til{f}_8).
\end{align}
So the problem boils down to showing \eqref{eq:ATA_ATA_to_prove}, assuming \eqref{eq:ATA_ATA_to_assume}. Let us rewrite \eqref{eq:ATA_ATA_to_prove} by switching around some factors, so that it becomes an equation in $V_{\lambda_1} \otimes V_{\lambda_2^*} \otimes V_{\lambda_3} \otimes V_{\lambda_0^*}$:
\begin{align}
\label{eq:ATA_ATA_to_prove2}
\textstyle  \sum \mathrm{ev}_{14} (\til{f}_6 \otimes \til{f}_5)
= q^{2m} \, \mathbf{P}_{(1432)} \sum \mathrm{ev}_{36} (\til{f}_8 \otimes \til{f}_7).
\end{align}

\vs

Note
\begin{align*}
\mathrm{ev}_{14} (\til{f}_6 \otimes \til{f}_5)
& = \mathrm{ev}_{14} (D_{\lambda_4})_1^{-1} (C_{\lambda_2^*})_3^{-1} (D_{\lambda_4^*})_4 (C_{\lambda_0})_6 (\til{f}_2 \otimes \til{f}_1), \\
\mathrm{ev}_{36} (\til{f}_8 \otimes \til{f}_7)
& = \mathrm{ev}_{36} (D_{\lambda_0})_1^{-1} (C_{\lambda_5^*})_3^{-1} (D_{\lambda_2^*})_4 (C_{\lambda_5})_6 (\til{f}_4 \otimes \til{f}_3)
\end{align*}
Let's establish a little lemma:
\begin{lemma}
For any non-singular weight $\lambda=(x,y)$,
\begin{align*}
& \mbox{the map } \mathrm{ev}_{12} \circ (D_\lambda^{-1} \otimes D_{\lambda^*}) : {}^*V_\lambda \otimes V_\lambda \to V_{\lambda^*} \otimes  {}^* V_{\lambda^*} \to \mathbb{C} \mbox{ equals } \mathrm{ev}_{12} \circ A_2^{-1}, \\
& \mbox{and the map } \mathrm{ev}_{12} \circ (C_{\lambda^*}^{-1} \otimes C_\lambda) : V_\lambda \otimes V^*_\lambda \to V^*_{\lambda^*}\otimes V_{\lambda^*} \mbox{ equals } \mathrm{ev}_{12} \circ A_2.
\end{align*}
\end{lemma}
{\it Proof.} Note $\lambda^* = (x^{-1}, -yx^{-1})$. For the first assertion, observe that $e_i \otimes e_j \in {}^*V_\lambda \otimes V_\lambda$ is sent first by $D_\lambda^{-1} \otimes D_{\lambda^*}$ to $q^{-2im_{y,x^{-1}}} q^{-i(i+1)} e_{-i} \otimes q^{-2jm_{-yx^{-1},x}} q^{j(j-1)} e_{-j}$, and then by $\mathrm{ev}_{12}$ to the number $q^{-2im_{y,x^{-1}}} q^{-2jm_{-yx^{-1},x}} q^{-i(i+1)} q^{j(j-1)}  \delta_{-i,-j}$ which equals $q^{-2jm_{y,x^{-1}}} q^{2jm_{-yx^{-1},x}} q^{-2j} \delta_{i,j} = q^{-2j} \delta_{i,j}$; we used 
\begin{align}
\label{eq:m_observation1}
m_{-yx^{-1},x} = m_{yx^{-1},x} = - m_{y,x^{-1}},
\end{align}
which is straightforward to show.

\vs

For the second assertion, note that $e_i \otimes e_j \in V_\lambda \otimes V_\lambda^*$ is sent by $C_{\lambda^*}^{-1} \otimes C_\lambda$ to $q^{-2im_{-yx^{-1},x}} q^{i(i+1)} e_{-i} \otimes q^{-2jm_{y,x^{-1}}} q^{-j(j-1)} e_{-j}$ then by $\mathrm{ev}_{12}$ to $q^{-2im_{-yx^{-1},x}} q^{-2jm_{y,x^{-1}}} q^{i(i+1)} q^{-j(j-1)} \delta_{-i,-j}$, which simplifies to $q^{2j} \delta_{i,j}$. \qed

\vs

Using this lemma, we have
\begin{align*}
\mathrm{ev}_{14} (\til{f}_6 \otimes \til{f}_5) & = \mathrm{ev}_{14} (C_{\lambda_2^*})_3^{-1} (C_{\lambda_0})_6 \, A_4^{-1} (\til{f}_2 \otimes \til{f}_1), \\
\mathbf{P}_{(1432)} \mathrm{ev}_{36} (\til{f}_8 \otimes \til{f}_7)
& = \mathbf{P}_{(1432)} \mathrm{ev}_{36} (D_{\lambda_0})_1^{-1} (D_{\lambda_2^*})_4 \, A_6 (\til{f}_4 \otimes \til{f}_3).
\end{align*}
Apply the invertible linear map $(C_{\lambda_2^*})_2 \, A_3^{-1} (D_{\lambda_0})_4$ from left on both:
\begin{align*}
(C_{\lambda_2^*})_2 \, A_3^{-1} (D_{\lambda_0})_4 \, \mathrm{ev}_{14} \, (\til{f}_6 \otimes \til{f}_5) & = (C_{\lambda_2^*})_2 \, A_3^{-1} (D_{\lambda_0})_4 \, \mathrm{ev}_{14} \, (C_{\lambda_2^*})_3^{-1} (C_{\lambda_0})_6 \, A_4^{-1} (\til{f}_2 \otimes \til{f}_1), \\
& = \mathrm{ev}_{14} \, \cancel{ (C_{\lambda_2^*})_3 }  \, A_5^{-1} \ul{ (D_{\lambda_0})_6 } \, \cancel{ (C_{\lambda_2^*})_3^{-1} } \, \ul{ (C_{\lambda_0})_6 } \, A_4^{-1} (\til{f}_2 \otimes \til{f}_1), \\
& \stackrel{\eqref{eq:DC_and_CD}}{=} \mathrm{ev}_{14} \, A_4^{-1} A_5^{-1} A_6 \, (\til{f}_2 \otimes \til{f}_1), \\
(C_{\lambda_2^*})_2 \, A_3^{-1} (D_{\lambda_0})_4 \mathbf{P}_{(1432)} \, \mathrm{ev}_{36}\, (\til{f}_8 \otimes \til{f}_7)
& = (C_{\lambda_2^*})_2 \, A_3^{-1} (D_{\lambda_0})_4  \mathbf{P}_{(1432)} \, \mathrm{ev}_{36}\, (D_{\lambda_0})_1^{-1} (D_{\lambda_2^*})_4 \, A_3 (\til{f}_4 \otimes \til{f}_6) \\
& = \mathbf{P}_{(1432)} (C_{\lambda_2^*})_3 \, A_4^{-1} (D_{\lambda_0})_1 \mathrm{ev}_{36}\, (D_{\lambda_0})_1^{-1} (D_{\lambda_2^*})_4 \, A_6 (\til{f}_4 \otimes \til{f}_3) \\
& = \mathbf{P}_{(1432)} \, \mathrm{ev}_{36} \, \ul{ (C_{\lambda_2^*})_4 } \,A_5^{-1} \cancel{ (D_{\lambda_0})_1 } \cancel{ (D_{\lambda_0})_1^{-1} } \, \ul{ (D_{\lambda_2^*})_4 } \, A_6 (\til{f}_4 \otimes \til{f}_3) \\
& \stackrel{\eqref{eq:DC_and_CD}}{=} \mathbf{P}_{(1432)} \mathrm{ev}_{36} \, A_4^{-1} A_5^{-1} A_6 (\til{f}_4 \otimes \til{f}_3).
\end{align*}
As $\til{f}_1 \in \mathrm{Inv}(V_{\lambda_4} \otimes V_{\lambda_3} \otimes V_{\lambda_0}^*)$, the element $X\in \mathcal{W}$ acts as counit $\epsilon(X)=1$; since $(\Delta \otimes \mathrm{id})\circ \Delta(X) = X\otimes X\otimes X$, the $X$-action is $\mu_{\lambda_4}(X) \otimes \mu_{\lambda_3}(X) \otimes \mu^*_{\lambda_0}(X) = x_4^{1/N} x_3^{1/N} x_0^{-1/N} A \otimes A\otimes A^{-1}$ (see \eqref{eq:left_dual_of_cyclic_irreducible} for $\mu^*_\lambda$), so $A_4^{-1} A_5^{-1} A_6(\til{f}_2 \otimes \til{f}_1) = x_4^{1/N} x_3^{1/N} x_0^{-1/N} (\til{f}_2 \otimes \til{f}_1)$.  Likewise, since $\til{f}_3 \in \mathrm{Inv}(V_{\lambda_2} \otimes V_{\lambda_3} \otimes V_{\lambda_5}^*)$, the $X$-action $\mu_{\lambda_2}(X) \otimes \mu_{\lambda_3}(X) \otimes \mu^*_{\lambda_5}(X) = x_2^{1/N} x_3^{1/N} x_5^{-1/N} A \otimes A\otimes A^{-1}$ on $\til{f}_3$ is just multiplication by $\epsilon(X)=1$, so that $A_4^{-1} A_5^{-1} A_6 (\til{f}_4 \otimes \til{f}_3) = x_2^{1/N} x_3^{1/N} x_5^{-1/N} (\til{f}_4 \otimes \til{f}_3)$. 

\vs

So, we showed that applying the invertible linear map $(C_{\lambda_2^*})_2 \, A_3^{-1} (D_{\lambda_0})_4$ from left to both sides of the sought-for equation \eqref{eq:ATA_ATA_to_prove2} yields
$$
\textstyle x_4^{1/N} x_3^{1/N} x_0^{-1/N} \, \sum \mathrm{ev}_{14} (\til{f}_2 \otimes \til{f}_1)
= q^{2m} \, x_2^{1/N} x_3^{1/N} x_5^{-1/N} \, \mathbf{P}_{(1432)} \, \sum \mathrm{ev}_{36} (\til{f}_4 \otimes \til{f}_3),
$$
which is same as the assumption \eqref{eq:ATA_ATA_to_assume}, provided that $x_4^{1/N} x_3^{1/N} x_0^{-1/N} = q^{2m} \, x_2^{1/N} x_3^{1/N} x_5^{-1/N}$, which is the last remaining thing to check. First, cancel $x_3^{1/N}$. From \eqref{eq:lambda_1_to_5}, we have $x_4 = x_1x_2$, $x_5 = x_2 x_3$, $x_0 = x_1x_2x_3$. So,
\begin{align*}
\textstyle \frac{ x_4^{1/N} x_5^{1/N} }{ x_2^{1/N} x_0^{1/N} }
= \frac{ (x_1x_2)^{1/N} (x_2x_3)^{1/N} }{ x_2^{1/N} (x_1x_2x_3)^{1/N} }
= \frac{ q^{2m_{x_1,x_2}} x_1^{1/N} x_2^{1/N} (x_2x_3)^{1/N}}{ x_2^{1/N} q^{2m_{x_1,x_2x_3}} x_1^{1/N} (x_2x_3)^{1/N} } = q^{2(m_{x_1,x_2} - m_{x_1,x_2x_3})} = q^{2m},
\end{align*}
finishing the proof. \qed

\subsection{Second consistency relation involving ${\bf T}$ and ${\bf A}$}

Let us use same notation for $\lambda,\lambda',\lambda''$, $\lambda_i$, $i=0,1,\ldots,5$, from the previous subsection. The second consistency relation of ${\bf T}$ and ${\bf A}$ corresponding to the diagram in Fig.\ref{fig:TAT_AAP} is formulated and shown as follows. Again, we give a representation theoretic proof.
\begin{proposition}[${\bf TAT} = {\bf AAP}$ relation]
\label{prop:TAT_AAP}
Let $\lambda_0,\ldots,\lambda_5$ be as above. Consider the two ${\bf T}$ maps
\begin{align*}
{\bf T}_{\lambda_1,\lambda_2,\lambda_3} : M_{\lambda_4,\lambda_3}^{\lambda_0} \otimes M_{\lambda_1,\lambda_2}^{\lambda_4} \to M_{\lambda_2,\lambda_3}^{\lambda_5} \otimes M_{\lambda_1,\lambda_5}^{\lambda_0}, \\
{\bf T}_{\lambda_2,\lambda_3,\lambda_0^*} : M_{\lambda_5,\lambda_0^*}^{\lambda_1^*} \otimes M_{\lambda_2,\lambda_3}^{\lambda_5} \to M_{\lambda_3,\lambda_0^*}^{\lambda_4^*} \otimes M_{\lambda_2,\lambda_4^*}^{\lambda_1^*}.
\end{align*}
Then this map ${\bf T}_{\lambda_2,\lambda_3,\lambda_0^*}$ written as above makes sense as a ${\bf T}$ map, and one has
\begin{align}
\label{eq:TAT_AAP}
({\bf T}_{\lambda_2,\lambda_3,\lambda_0^*})_{12} \, ({\bf A}_{\lambda_1,\lambda_5}^{\lambda_0})_1 \, ({\bf T}_{\lambda_1,\lambda_2,\lambda_3})_{21}
= ({\bf A}_{\lambda_4,\lambda_3}^{\lambda_0})_1 \, ({\bf A}_{\lambda_1,\lambda_2}^{\lambda_4})_2 \, \mathbf{P}_{(12)},
\end{align}
as maps $M_{\lambda_1,\lambda_2}^{\lambda_4} \otimes M_{\lambda_4,\lambda_3}^{\lambda_0} \to M_{\lambda_3,\lambda_0^*}^{\lambda_4^*} \otimes M_{\lambda_2,\lambda_4^*}^{\lambda_1^*}$.
\end{proposition}
One can symbolically write the above relation \eqref{eq:TAT_AAP} as ${\bf T}_{12} {\bf A}_1 {\bf T}_{21} = {\bf A}_1 {\bf A}_2 {\bf P}_{(12)}$.

\vs

{\it Proof.} Again, we prove the Hom version. The Hom versions of the two ${\bf T}$ maps are
\begin{align*}
(\mathbf{T}_{\lambda_1,\lambda_2,\lambda_3}^\mathrm{Hom})_{21} & : \mathrm{Hom}_\mathcal{W}(V_{\lambda_4}, V_{\lambda_1} \otimes V_{\lambda_2}) \otimes \mathrm{Hom}_\mathcal{W}(V_{\lambda_0},V_{\lambda_4}\otimes V_{\lambda_3}) \\
& \quad \to \mathrm{Hom}_\mathcal{W}(V_{\lambda_0}, V_{\lambda_1} \otimes V_{\lambda_5}) \otimes \mathrm{Hom}_\mathcal{W}(V_{\lambda_5}, V_{\lambda_2} \otimes V_{\lambda_3}), \\
(\mathbf{T}_{\lambda_2,\lambda_3,\lambda_0^*}^\mathrm{Hom})_{12} & : \mathrm{Hom}_\mathcal{W}(V_{\lambda_1^*}, V_{\lambda_5} \otimes V_{\lambda_0^*}) \otimes \mathrm{Hom}_\mathcal{W}(V_{\lambda_5}, V_{\lambda_2} \otimes V_{\lambda_3}) \\
& \quad \to \mathrm{Hom}_\mathcal{W}(V_{\lambda_4^*}, V_{\lambda_3} \otimes V_{\lambda_0^*}) \otimes \mathrm{Hom}_\mathcal{W}(V_{\lambda_1^*}, V_{\lambda_2} \otimes V_{\lambda_4^*}).
\end{align*}
Pick any element $\sum f_1 \otimes f_2 \in \mathrm{Hom}_\mathcal{W}(V_{\lambda_4}, V_{\lambda_1} \otimes V_{\lambda_2}) \otimes \mathrm{Hom}_\mathcal{W}(V_{\lambda_0}, V_{\lambda_4} \otimes V_{\lambda_3})$. We shall apply the Hom version of both sides of the sought-for equation ${\bf TAT} = {\bf AAP}$ to this element. First, the RHS, which is easier; properly understood, we get $\sum {\bf A}^\mathrm{Hom} f_2 \otimes {\bf A}^\mathrm{Hom} f_1$. Now, the LHS. As in Lem.\ref{lem:T_Hom}, the element $\mathbf{T}^\mathrm{Hom}_{21}(\sum f_1 \otimes f_2) = \sum f_3 \otimes f_4$ of $\mathrm{Hom}_\mathcal{W}(V_{\lambda_0}, V_{\lambda_1} \otimes V_{\lambda_5}) \otimes \mathrm{Hom}_\mathcal{W}(V_{\lambda_5}, V_{\lambda_2} \otimes V_{\lambda_3})$ can be described completely by the equation
\begin{align}
\label{eq:TAT_AAP_to_assume1}
\sum (f_1\otimes \mathrm{id})\circ f_2 = \sum (\mathrm{id}\otimes f_4) \circ f_3
\end{align}
of elements in $\mathrm{Hom}_\mathcal{W}(V_{\lambda_0}, V_{\lambda_1}\otimes V_{\lambda_2} \otimes V_{\lambda_3})$. We then apply ${\bf T}_{12}^\mathrm{Hom} {\bf A}_1^\mathrm{Hom}$ to $\sum f_3 \otimes f_4$. So the task is to show the equation 
$\textstyle {\bf T}^\mathrm{Hom}_{12} \sum (\mathbf{A}^\mathrm{Hom} f_3) \otimes f_4
= \sum {\bf A}^\mathrm{Hom} f_2 \otimes {\bf A}^\mathrm{Hom} f_1
$ 
 of elements of $\mathrm{Hom}_\mathcal{W}(V_{\lambda_4^*}, V_{\lambda_3} \otimes V_{\lambda_0^*}) \otimes \mathrm{Hom}_\mathcal{W}(V_{\lambda_1^*}, V_{\lambda_2} \otimes V_{\lambda_4^*})$ under the assumption \eqref{eq:TAT_AAP_to_assume1}, where each ${\bf T}^\mathrm{Hom}$ and ${\bf A}^\mathrm{Hom}$ should be understood appropriately. Lem.\ref{lem:T_Hom} says that this sought-for equation is equivalent to the equation
\begin{align}
\label{eq:TAT_AAP_to_show1}
\textstyle \sum (f_4 \otimes \mathrm{id}) \circ ({\bf A}^\mathrm{Hom} f_3) = \sum (\mathrm{id}\otimes ({\bf A}^\mathrm{Hom} f_2)) \circ (\mathbf{A}^\mathrm{Hom} f_1)
\end{align}
of elements in $\mathrm{Hom}_\mathcal{W}(V_{\lambda_1^*}, V_{\lambda_2} \otimes V_{\lambda_3} \otimes V_{\lambda_0^*})$. Hence, the problem is to show \eqref{eq:TAT_AAP_to_show1} using \eqref{eq:TAT_AAP_to_assume1}.

\vs

We now go to Inv versions. Define the following elements via the left and right canonical maps \eqref{eq:Hom_as_tensor_products}:
\begin{align*}
& \til{f}_1 := J^\mathrm{L} f_1 \in \mathrm{Inv}(V_{\lambda_1} \otimes V_{\lambda_2} \otimes V_{\lambda_4}^*), \qquad
\til{f}_2 := J^\mathrm{L} f_2 \in \mathrm{Inv}(V_{\lambda_4}\otimes V_{\lambda_3} \otimes V_{\lambda_0}^*), \\
& \til{f}_3 := J^\mathrm{L} f_3 \in \mathrm{Inv}(V_{\lambda_1} \otimes V_{\lambda_5} \otimes V_{\lambda_0}^*), \qquad
\til{f}_4 := J^\mathrm{L} f_4 \in \mathrm{Inv}(V_{\lambda_2} \otimes V_{\lambda_3} \otimes V_{\lambda_5}^*), \\
& \til{f}_5 := J^\mathrm{R} {\bf A}^\mathrm{Hom} f_1 \in \mathrm{Inv}({}^*V_{\lambda_1^*} \otimes V_{\lambda_2} \otimes V_{\lambda_4^*}), \qquad
\til{f}_6 := J^\mathrm{R}{\bf A}^\mathrm{Hom} f_2 \in \mathrm{Inv}({}^*V_{\lambda_4^*} \otimes V_{\lambda_3} \otimes V_{\lambda_0^*}), \\
& \til{f}_7 := J^\mathrm{R} {\bf A}^\mathrm{Hom} f_3 \in \mathrm{Inv}({}^*V_{\lambda_1^*} \otimes V_{\lambda_5} \otimes V_{\lambda_0^*}).
\end{align*}
Definition \eqref{eq:A_as_Hom} of the ${\bf A}^\mathrm{Hom}$ map says
$$
\til{f}_5 = (D_{\lambda_1^*})_1 (C_{\lambda_4})_3 \til{f}_1, \qquad
\til{f}_6 = (D_{\lambda_4^*})_1 (C_{\lambda_0})_3 \til{f}_2, \qquad
\til{f}_7 = (D_{\lambda_1^*})_1 (C_{\lambda_0})_3 \til{f}_3.
$$
By the same philosophy as in the previous subsection (I omit a detailed argument this time), the equation \eqref{eq:TAT_AAP_to_assume1} can be translated to the equation
\begin{align}
\label{eq:TAT_AAP_to_assume2}
\textstyle \sum \mathrm{ev}_{34} (\til{f}_1 \otimes \til{f}_2)
= \mathbf{P}_{(243)} \sum \mathrm{ev}_{26} (\til{f}_3 \otimes \til{f}_4)
\end{align}
of elements in $V_{\lambda_1} \otimes V_{\lambda_2} \otimes V_{\lambda_3} \otimes V_{\lambda_0}^*$, while the equation \eqref{eq:TAT_AAP_to_show1} to
\begin{align}
\label{eq:TAT_AAP_to_show2}
\textstyle {\bf P}_{(243)} \sum \mathrm{ev}_{26}(\til{f}_7 \otimes \til{f}_4) =  \sum \mathrm{ev}_{34}(\til{f}_5 \otimes \til{f}_6)
\end{align}
of elements in ${}^*V_{\lambda_1^*} \otimes V_{\lambda_2} \otimes V_{\lambda_3} \otimes V_{\lambda_0^*}$. So the Inv version of the problem is to assume \eqref{eq:TAT_AAP_to_assume2} and show \eqref{eq:TAT_AAP_to_show2}. Apply the invertible linear map $(D_{\lambda_1^*})_1^{-1} (C_{\lambda_0})_4^{-1}$ from left to both sides of \eqref{eq:TAT_AAP_to_show2}, so the sought-for equation is now
\begin{align}
\label{eq:TAT_AAP_to_show3}
\textstyle (D_{\lambda_1^*})_1^{-1} (C_{\lambda_0})_4^{-1} {\bf P}_{(243)} \sum \mathrm{ev}_{26}(\til{f}_7 \otimes \til{f}_4) = (D_{\lambda_1^*})_1^{-1} (C_{\lambda_0})_4^{-1} \sum \mathrm{ev}_{34}(\til{f}_5 \otimes \til{f}_6),
\end{align}
which is an equation in $V_{\lambda_1} \otimes V_{\lambda_2} \otimes V_{\lambda_3} \otimes V_{\lambda_0}^*$. Note
\begin{align*}
(D_{\lambda_1^*})_1^{-1} (C_{\lambda_0})_4^{-1} \, {\bf P}_{(243)} \, \mathrm{ev}_{26} (\til{f}_7 \otimes \til{f}_4) & = (D_{\lambda_1^*})_1^{-1} (C_{\lambda_0})_4^{-1} {\bf P}_{(243)} \, \mathrm{ev}_{26} \, (D_{\lambda_1^*})_1 (C_{\lambda_0})_3 \, (\til{f}_3 \otimes \til{f}_4), \\
& = {\bf P}_{(243)} \, (D_{\lambda_1^*})_1^{-1} (C_{\lambda_0})_2^{-1} \, \mathrm{ev}_{26} \, (D_{\lambda_1^*})_1 (C_{\lambda_0})_3 \, (\til{f}_3 \otimes \til{f}_4), \\
& = {\bf P}_{(243)}  \, \mathrm{ev}_{26} \, (D_{\lambda_1^*})_1^{-1} (C_{\lambda_0})_3^{-1} (D_{\lambda_1^*})_1 (C_{\lambda_0})_3 \, (\til{f}_3 \otimes \til{f}_4), \\
& = {\bf P}_{(243)} \, \mathrm{ev}_{26} \, (\til{f}_3 \otimes \til{f}_4), \\
(D_{\lambda_1^*})_1^{-1} (C_{\lambda_0})_4^{-1} \, \mathrm{ev}_{34} (\til{f}_5 \otimes \til{f}_6) & = (D_{\lambda_1^*})_1^{-1} (C_{\lambda_0})_4^{-1} \, \mathrm{ev}_{34} \, (D_{\lambda_1^*})_1 (C_{\lambda_4})_3 (D_{\lambda_4^*})_4 (C_{\lambda_0})_6 (\til{f}_1 \otimes \til{f}_2) \\
& = \mathrm{ev}_{34} \, (D_{\lambda_1^*})_1^{-1} (C_{\lambda_0})_6^{-1} (D_{\lambda_1^*})_1 (C_{\lambda_4})_3 (D_{\lambda_4^*})_4 (C_{\lambda_0})_6 (\til{f}_1 \otimes \til{f}_2) \\
& = \mathrm{ev}_{34} \, (C_{\lambda_4})_3 (D_{\lambda_4^*})_4 (\til{f}_1 \otimes \til{f}_2).
\end{align*}
A little lemma:
\begin{lemma}
For any non-singular weight $\lambda=(x,y)$,
$$
\mbox{the map $\mathrm{ev}_{12} \circ (C_\lambda \otimes D_{\lambda^*}) : V^*_\lambda \otimes V_\lambda \to V_{\lambda^*} \otimes {}^* V_{\lambda^*} \to \mathbb{C}$ equals $\mathrm{ev}$.}
$$
\end{lemma}
{\it Proof.} Note $\lambda^* = (x^{-1}, -yx^{-1})$. Observe $e_i \otimes e_j \in V^*_\lambda \otimes V_\lambda$ is sent via $C_\lambda \otimes D_{\lambda^*}$ to $q^{-2im_{y,x^{-1}}} q^{-i(i-1)} e_{-i} \otimes q^{-2jm_{-yx^{-1},x}} q^{j(j-1)} e_{-j} \in V_{\lambda^*} \otimes {}^* V_{\lambda^*}$, then via $\mathrm{ev}_{12}$ to $q^{-2im_{y,x^{-1}}} q^{-i(i-1)} q^{-2jm_{-yx^{-1},x}} q^{j(j-1)} \delta_{-i,-j}$ which equals $q^{-2i(m_{y,x^{-1}}+m_{-yx^{-1},x})}  \delta_{i,j} = \delta_{i,j}$, where we used \eqref{eq:m_observation1}. \qed

\vs

Hence $(D_{\lambda_1^*})_1^{-1} (C_{\lambda_0})_4^{-1} \, \mathrm{ev}_{34} (\til{f}_5 \otimes \til{f}_6) = \mathrm{ev}_{34} (\til{f}_1 \otimes \til{f}_2)$. So we just showed that the sought-for equation \eqref{eq:TAT_AAP_to_show3} is equivalent to the assumption equation \eqref{eq:TAT_AAP_to_assume2}. \qed

\subsection{Groupoid of triangulations of $n$-gon}
\label{subsec:main}

We shall now assemble all the results obtained so far and formulate our main result. We first recall from the literature the notion of an `ideal triangulation' of bordered surfaces with marked points on the boundary:
\begin{definition}
\label{def:ideal_triangulation}
A \ul{\em bordered surface with marked points on the boundary} is an oriented topological surface $S$ possibly with boundary, possibly with distinguished `marked' points on the boundary. Boundary minus the (possibly empty) set of marked points may have several connected components; a component not homeomorphic to a circle is called a \ul{\em boundary arc}, and any component is called a \ul{\em boundary piece}.

An \ul{\em ideal triangulation} of $S$ is defined as follows. An \ul{\em ideal edge} of $S$ is a nontrivial homotopy class of unoriented paths in $S$ which intersect with the boundary of $S$ only with its endpoints at boundary pieces of $S$, where the homotopy allows each endpoint to move inside one boundary piece. A \ul{\em boundary edge} of $S$ is a nontrivial homotopy class of unoriented paths contained in the boundary of $S$ and containing exactly one marked point, whose endpoints lie in boundary pieces and the homotopy allows each endpoint to move inside one boundary piece. An \ul{\em ideal triangulation} is a maximal collection of distinct non-intersecting ideal edges and boundary edges.
\end{definition}
For a more precise and detailed treatment, see Fock-Goncharov \cite{FG06}, Fomin-Shapiro-Thurston \cite{FST}, or references therein. An ideal triangulation divides the surface into \ul{\em (ideal) triangles}. We shall only deal with the following case:
\begin{definition}
For a positive integer $n$, an \ul{\em $n$-gon} refers to a surface $S$ as in Def.\ref{def:ideal_triangulation}, of genus $0$ with one boundary component, with $n$ marked points on the boundary.
\end{definition}
So one can think of an `$n$-gon' as a usual Euclidean polygon with $n$ sides, and an ideal triangulation of it as a triangulation with straight line segments. Ideal triangulations (possibly with labeling for triangles) are what are used in Chekhov-Fock(-Goncharov)'s quantization of Teichm\"uller spaces, while Kashaev's quantization of Teichm\"uller spaces uses the following.
\begin{definition}
A \ul{\em dotted ideal triangulation} of a bordered surface with marked points on the boundary $S$ is an ideal triangulation of $S$, together with 1) the choice of a distinguished corner for each triangle, depicted with the dot $\bullet$, and with 2) the choice of (bijective) labeling for the triangles by some index set $I$.
\end{definition}
We usually take $I = \{1,2,\ldots,r\}$, where $r$ is the number of triangles, which is an invariant. For example, for an $n$-gon, we have $r=n-2$. For the present paper we upgrade to the following notion, as could have been predicted in \S\ref{subsec:graphical_encoding}:
\begin{definition}
A \ul{\em LD-triangulation} of $S$ (or a \ul{\em labeled dotted ideal triangulation} of $S$) is a dotted ideal triangulation of $S$ together with the choice of labeling of each edge by a non-singular weight, called \ul{\em edge-labels}.

An LD-triangulation of $S$ is said to be \ul{\em sane} if each triangle is sane, in the sense of Def.\ref{def:sane}.
\end{definition}

In Kashaev's quantum Teichm\"uller theory, choice of a dotted ideal triangulation leads to quantum coordinate operators on a certain version of Teichm\"uller space on a Hilbert space of states. Different ideal triangulation leads to different quantum coordinate operators on a different Hilbert space, and to each change of dotted ideal triangulation Kashaev assigns a unitary map between the Hilbert spaces, intertwining the quantum operators. In the classical limit, this unitary intertwiner recovers the classical coordinate change formulas, and the assignment of these intertwining operators is consistent, in the sense that the composition of changes of dotted ideal triangulations is preserved by the corresponding unitary operators, up to constants. One way of formulating this is by Kashaev's groupoid of dotted ideal triangulations. For convenience, we first come up with:
\begin{definition}
A \ul{\em full groupoid based on a set} $C$ is the category whose set of objects is $C$, in which there is exactly one morphism from any object to any object.
\end{definition}

\begin{definition}[groupoids of (changes of) triangulations]
\label{def:three_groupoids}
The \ul{\em Ptolemy groupoid} of $S$ is the full groupoid based on the set of all ideal triangulations of $S$.

The \ul{\em Kashaev groupoid} of $S$ is the full groupoid based on the set of all dotted ideal triangulations of $S$.

The \ul{\em LD-groupoid} of $S$ is the full groupoid based on the set of all LD-triangulations of $S$.

\end{definition}
The Chekhov-Fock(-Goncharov) quantization of Teichm\"uller space of a surface $S$ can be thought of as a projective functor from the Ptolemy groupoid of $S$ to the category $\mathrm{Hilb}$ of Hilbert spaces whose morphisms are unitary maps, while the Kashaev quantization as a projective functor from the Kashaev groupoid of $S$ to $\mathrm{Hilb}$. Here, functor is said to be `projective' if the composition of morphisms is preserved up to multiplicative constants. These results are written in terms of generators, or generating morphisms, of the corresponding groupoids. We first view each element of a groupoid in Def.\ref{def:three_groupoids} as a `change' of corresponding triangulations.
\begin{definition}[elementary morphisms]
In a Ptolemy groupoid, a morphism is called a \ul{\em flip} if it connects two ideal triangulations differing only by one edge.

In a Kashaev groupoid: a morphism is denoted by ${\bf A}_t$ if it alters only one triangle labeled by $t$ by moving its dot to the counterclockwise next one. A morphism is denoted by ${\bf T}_{st}$ if triangles labeled by $s,t$ are adjacent with the dots configured as in Fig.\ref{fig:action_on_dotted_ideal_triangulations} in the first triangulation, and the second triangulation is the same as the first one outside these triangles, on which it looks like in Fig.\ref{fig:action_on_dotted_ideal_triangulations}. A morphism is denoted by ${\bf P}_\gamma$ for a permutation $\gamma$ of index set for triangle labels, if it changes the label $t$ for each triangle in the first triangulation to the label $\gamma(t)$ in the second triangulation.

In an LD-groupoid: a morphism is denoted by ${\bf A}^{\epsilon_I}$, for a sequence $(\epsilon_t)_{t\in I} \in \{\pm1\}^I$ of signs, if it is $\prod_{t\in I} {\bf A}_t^{\epsilon_t}$ as in the Kashaev groupoid when forgetting the edge labels, and the edge labels change as in Fig.\ref{fig:A}. A morphism is denoted by ${\bf T}_{st}$ if it is as in the Kashaev groupoid, and the edge labels change as in Fig.\ref{fig:T}. A morphism is denoted by ${\bf P}_\gamma$ if it is as in the Kashaev groupoid.

Call these morphisms and their inverses the \ul{\em elementary morphisms} of each groupoid.
\end{definition}

\begin{theorem}[``Whitehead's classical fact''; see e.g. \cite{Penner}]
The Ptolemy groupoid of $S$ is generated by flips, i.e. any morphism is a sequence of flips, and any algebraic relation of flips is generated by, i.e. is a consequence of, the following ones: 1) flip on the same edge twice is identity, 2) for two edges appearing as sides of exactly one triangle, the alternating sequence of five flips on these two edges is the identity, and 3) for two edges not appearing as sides of any triangle, the flips on these edges commute with each other.
\end{theorem}

\begin{theorem}[\cite{Kash00},\cite{Teschner},\cite{Kim16b}]
\label{thm:Kashaev_presentation}
The Kashaev groupoid of $S$ is generated by elementary morphisms, and any algebraic relation among them is generated by the ones in eq.\eqref{eq:Kashaev_relations} and the `trivial ones' mentioned after eq.\eqref{eq:Kashaev_relations}.
\end{theorem}
In Chekhov-Fock(-Goncharov) and Kashaev quantizations, they assign unitary operators to the elementary morphisms, and verify that the relations corresponding to the ones above are satisfied by these operators up to constants. For our situation, we would like to consider a subgroupoid of LD-groupoid consisting of sane LD-triangulations. However, instead of full subgroupoid based on the set of all sane LD-triangulations, we consider the subgroupoid generated by elementary morphisms. 
\begin{definition}
The \ul{\em sane LD-groupoid} of $S$ is the subgroupoid of the LD-groupoid of $S$ with the set of objects being the set of all sane LD-triangulations, in which there is a morphism from an object to an object if and only if these two are connected by a sequence of elementary morphisms in the LD-groupoid. 
\end{definition}
In particular, not every two objects in the sane LD-groupoid are connected by a morphism. In practice, it is probably wise to take a `connected component' of the sane LD-groupoid which is of interest; one may look for a characterization of connected components, in terms of the condition on the edge labels. From Thm.\ref{thm:Kashaev_presentation} we get:
\begin{corollary}
Any algebraic relation of elementary morphisms of an LD-groupoid is generated by the ones shown in Propositions \ref{prop:TT_TTT}, \ref{prop:AAA}, \ref{prop:ATA_ATA}, \ref{prop:TAT_AAP} (i.e. the ones symbolized by the relations in eq.\ref{eq:Kashaev_relations}), and the `trivial ones'. \qed
\end{corollary}
Now we are ready to formulate our main result in the following final form:
\begin{theorem}[main result: representation of the sane LD-groupoid, from representation theory of the quantum torus algebra $\mathcal{W}$]
\label{thm:main}
Let $\mathrm{Vec}$ be the category of finite dimensional complex vector spaces, where morphisms are linear maps. Let $S$ be an $n$-gon for a positive integer $n\ge 3$, and fix an index set $I$ for triangles of ideal triangulations of $S$. The following assignment yields a well-defined projective functor from the sane LD-groupoid of $S$ to $\mathrm{Vec}$.

Consider a sane LD-triangulation of $S$. For each triangle labeled by $t\in I$, denote the edge-labels of the sides of $t$ by $\lambda_{t_1}, \lambda_{t_2}, \lambda_{tx_3}$ in this counterclockwise order, so that $\lambda_{t_3}$ is for the side facing the dot of $t$. To this sane LD-triangulation we assign the vector space
$$
\bigotimes_{t\in I} M_{\lambda_{t_1}, \lambda_{t_2}}^{\lambda_{t_3}}
$$
For the elementary morphism ${\bf A}^{\epsilon_I}$ we assign the operator $\bigotimes_{t\in I} {\bf A}_t^{\epsilon_t}$, where ${\bf A}_t^{+1}$ stands for the map ${\bf A}_{\lambda_1,\lambda_2}^{\lambda_3} : M_{\lambda_{t_1}, \lambda_{t_2}}^{\lambda_{t_3}} \to M_{\lambda_{t_2}, \lambda_{t_3}^*}^{\lambda_{t_1}^*}$ and ${\bf A}_t^{-1}$ stands for the map $({\bf A}_{\lambda_{t_3}^*, \lambda_{t_1}}^{\lambda_{t_2}^*})^{-1} : M_{\lambda_{t_1},\lambda_{t_2}}^{\lambda_{t_3}} \to M_{\lambda_{t_3}^*,\lambda_{t_1}}^{\lambda_{t_2}^*}$. For the elementary morphism ${\bf T}_{st}$ we assign the operator $({\bf T}_{\lambda_{s_1}, \lambda_{s_2}, \lambda_{t_2}})_{st}$ which is the operator ${\bf T}_{\lambda_{s_1}, \lambda_{s_2}, \lambda_{t_2}} : M_{\lambda_{t_1},\lambda_{t_2}}^{\lambda_{t_3}} \otimes M_{\lambda_{s_1},\lambda_{s_2}}^{\lambda_{s_3}} \to M_{\lambda_{s_2},\lambda_{t_2}}^{\lambda_{s,t}} \otimes M_{\lambda_{s_1},\lambda_{s,t}}^{\lambda_{t_3}}$ applied on the $s$-th and $t$-th tensor factors, where $\lambda_{s,t} := \lambda_{s_2} \lambda_{t_2}$ is the edge-label for the new diagonal edge. For the elementary morphism ${\bf P}_\gamma$ we assign the factor permuting operator ${\bf P}_\gamma : \bigotimes_{t\in I} M_{\lambda_{t_1},\lambda_{t_2}}^{\lambda_{t_3}} \to \bigotimes_{t\in I} M_{\lambda_{\gamma^{-1}(t)_1}, \lambda_{\gamma^{-1}(t)_2}}^{\lambda_{\gamma^{-1}(t),3}}$ sending each $t$-th factor to the $\gamma(t)$-th factor. \qed
\end{theorem}
We may say that our final result is the construction of families of finite dimensional projective representations of the sane LD-groupoid of an $n$-gon, which is the genus $0$ bordered surface with one boundary component with $n$ marked points on the boundary, from the representation theory of the quantum torus algebra $\mathcal{W}=\mathcal{W}^q$ at root of unity $q$. 

\vs

As an anonymous referee pointed out to the author, the existence of these projective representations of the sane LD-groupoid was first asserted in Prop.9 of \cite[\S2.1.2]{K98} without proof, and then developed more concretely in \cite{Kash00b}; in particular, an explicit formula for the order three operator ${\bf A}$ is not written in \cite{K98} but appears in \cite{Kash00b}. The referee also mentioned that in a later work \cite{GKT}, the `tetrahedral symmetries' of $6$j-symbols were worked out in full generality in the case of the Borel subalgebra of $\mathcal{U}_q(\mathfrak{sl}_2)$, i.e. the quantum torus algebra $\mathcal{W}$. The 6j-symbols there correspond to our ${\bf T}$ operator, and I believe that our consistency relations involving ${\bf A}$ and ${\bf T}$ are somehow encoded in their `tetrahedral symmetries'.

\section{Discussion}
\label{sec:discussion}

\subsection{Going to higher genus: once-punctured torus}
\label{subsec:OPT}

I first claim that in the case when $S$ is an $n$-gon which is dealt with in the present paper, there exists a `big enough' connected component of the sane LD-groupoid of $S$, in the sense that for any ideal triangulation of the $n$-gon, this connected component contains an object (i.e. a sane LD-triangulation) whose underlying ideal triangulation is this one. For investigation of a similar question for higher genus cases, here we consider an example of the once-punctured torus.

\vs

\begin{figure}[htbp!]
\centering
\begin{pspicture}[showgrid=false](0,0.5)(8.4,3.5)
\rput[bl](0;0){
\PstSquare[unit=1.8,PolyName=P]
\pcline(P1)(P2)\ncput*{$\lambda_1$}
\pcline(P2)(P3)\ncput*{$\lambda_2$}
\pcline(P3)(P4)\ncput*{$\lambda_1$}
\pcline(P4)(P1)\ncput*{$\lambda_2$}
\pcline(P1)(P3)\ncput*{$\lambda_3$}
\rput[l]{-45}(P2){\hspace{1,4mm}$\bullet$}
\rput[l]{22}(P3){\hspace{2,4mm}$\bullet$}
}
\rput[l](2.5,1.3){$s$}
\rput[l](1.2,2.5){$t$}
\rput[bl](4.5;0){
\PstSquare[unit=1.8,PolyName=P]
\pcline(P1)(P2)\ncput*{$\lambda_1$}
\pcline(P2)(P3)\ncput*{$\lambda_2$}
\pcline(P3)(P4)\ncput*{$\lambda_1$}
\pcline(P4)(P1)\ncput*{$\lambda_2$}
\pcline(P2)(P4)\ncput*{$\lambda_4$}
\rput[l]{-23}(P2){\hspace{2,4mm}$\bullet$}
\rput[l]{41}(P3){\hspace{1,4mm}$\bullet$}
}
\rput[l](5.8,1.3){$s$}
\rput[l](6.8,2.5){$t$}
\rput[l](3.7,2){\pcline{->}(0,0)(1;0)\Aput{${\bf T}_{st}$}}
\end{pspicture}
\caption{Flip at an edge labeled by $\lambda_3$, for a once-punctured torus}
\label{fig:OPT}
\end{figure}
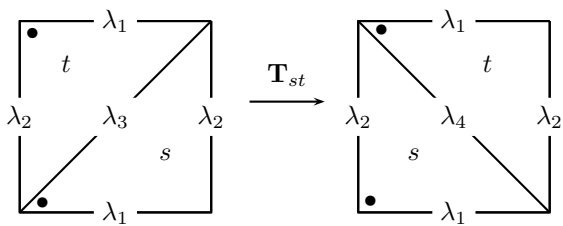

First, notice that once-punctured torus can be obtained by identifying the two pairs of the `parallel' outer edges of a $4$-gon, as in the LHS of Fig.\ref{fig:OPT}. At the moment, one may view this as being a special case of the $4$-gon situation; in particular, the identified edges of the $4$-gon are labeled by a same weight. The first question is, do there exist weights $\lambda_1=(x_1,y_1)$ and $\lambda_2=(x_2,y_2)$ so that the LD-triangulation in the LHS of Fig.\ref{fig:OPT} is sane? The condition is that $\lambda_1,\lambda_2,\lambda_3 = \lambda_1\lambda_2 = (x_1x_2,y_1x_2 + y_2)$ must be non-singular, and that $\lambda_1 \lambda_2 \lambda_1 = \lambda_2$, which as one can see by computation is equivalent to the condition $x_2(x_1^2-1)=0$ and $y_1(x_1x_2+1) + y_2(x_1-1)=0$. From non-singularity $x_1,x_2,y_1,y_2$ are nonzero, so $x_1^2-1=0$, hence $x_1 = \pm 1$. For $x_1=1$, it must be that $y_1 (x_2+1)=0$ hence $x_2=-1$, while for $x_1=-1$ it must be that $y_1(-x_2+1)-2y_2=0$ hence $x_2 = -\frac{2y_2}{y_1} +1$. For the moment, let us choose to work with the case $x_1=1$, $x_2=-1$. Then the non-singularity of $\lambda_3$ says $y_1\neq - y_2$. We would like to be able to flip at the diagonal labeled by $\lambda_3$; in order for the new LD-triangulation after the flip to be sane, it must be that $\lambda_4 = \lambda_2 \lambda_1 = (-1,-y_2+y_1)$ is non-singular, i.e. $y_2\neq y_1$. 

\vs

\begin{figure}[htbp!]
\centering
\begin{pspicture}[showgrid=false](0,-0.5)(8.4,4.5)
\rput[bl](2,1.5){
\pcline(-1.35,0.53)(1.2,0.53)\ncput*{$\lambda_1$}
\pcline(1.2,0.53)(1.2,3.08)\ncput*{$\lambda_2$}
\pcline(1.2,3.08)(-1.35,0.53)\ncput*{$\lambda_3$}
\pcline(-1.35,0.53)(-1.35,-2.02)\ncput*{$\lambda_2$}
\pcline(-1.35,-2.02)(1.2,0.53)\ncput*{$\lambda_3$}
\rput[l](0.3,1.3){$s$}
\rput[l](-0.7,-0.2){$t$}
\rput[l]{23}(-1.35,0.53){\hspace{2,4mm}$\bullet$}
\rput[l]{-45}(-1.35,0.53){\hspace{1,4mm}$\bullet$}
}
\rput[bl](6.5,1.5){
\pcline(-1.35,0.53)(1.2,0.53)
\pcline(1.2,0.53)(1.2,3.08)\ncput*{$\lambda_2$}
\pcline(1.2,3.08)(-1.35,0.53)\ncput*{$\lambda_3$}
\pcline(-1.35,0.53)(-1.35,-2.02)\ncput*{$\lambda_2$}
\pcline(-1.35,-2.02)(1.2,0.53)\ncput*{$\lambda_3^*$}
\rput[l](-0.2,0.8){$\lambda_1$}
\rput[l](-0.3,0.25){$\lambda_1^*$}
\rput[l](0.45,1.5){$s$}
\rput[l](-0.8,-0.35){$t$}
\rput[l]{23}(-1.35,0.53){\hspace{2,4mm}$\bullet$}
\rput[l]{67}(-1.35,-2.02){\hspace{3,4mm}$\bullet$}
}
\rput[l](3.7,2){\pcline{->}(0,0)(1;0)\Aput{${\bf A}_{t}$}}
\end{pspicture}
\caption{Preparing the flip at an edge labeled by $\lambda_1$, for a once-punctured torus}
\label{fig:OPT3}
\end{figure}
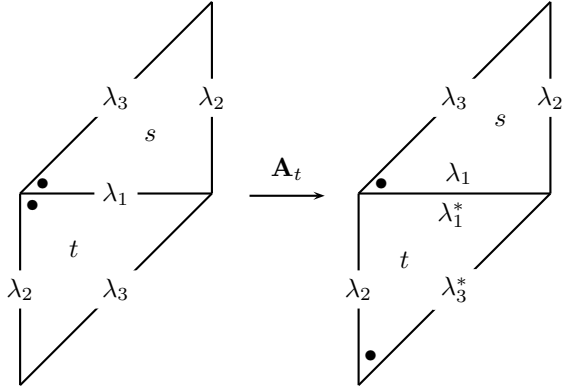

To make the connected component for this sane LD-triangulation bigger, we'd better be able to flip at other two edges of the triangulation too. To flip at the edge labeled by $\lambda_2$, let us redraw the once-punctured torus in the LHS of Fig.\ref{fig:OPT} as in the LHS of Fig.\ref{fig:OPT2}, and flip at the diagonal labeled by $\lambda_2$. For the newly obtained LD-triangulation to be sane, it must be that $\lambda_5 = \lambda_1 \lambda_3$ is non-singular, which one can compute to be equivalent to the condition $2y_1 \neq y_2$. What about flipping at the edge labeled by $\lambda_1$? We first redraw the once-punctured torus like in the LHS of Fig.\ref{fig:OPT3}. To be able to flip, we should first change the dots by some ${\bf A}$ move, for example to make the picture like in the RHS of Fig.\ref{fig:OPT3}. In order for this picture to make sense, we must have $\lambda_1 = \lambda_1^*$ and $\lambda_3=\lambda_3^*$. One can verify that $\lambda_1=\lambda_1^*$ is impossible because of $x_1=1$ (and, in case we chose $x_1=-1$ in the beginning, one can show that $\lambda_3 = \lambda_3^*$ is impossible). So the flip at the edge labeled by $\lambda_1$ is impossible.

\begin{figure}[htbp!]
\centering
\begin{pspicture}[showgrid=false](0,0.5)(8.4,3.5)
\rput[bl](0;0){
\pcline(-1.35,0.53)(1.2,0.53)\ncput*{$\lambda_1$}
\pcline(1.2,0.53)(3.75,3.08)\ncput*{$\lambda_3$}
\pcline(3.75,3.08)(1.2,3.08)\ncput*{$\lambda_1$}
\pcline(1.2,3.08)(-1.35,0.53)\ncput*{$\lambda_3$}
\pcline(1.2,0.53)(1.2,3.08)\ncput*{$\lambda_2$}
\rput[l]{23}(-1.35,0.53){\hspace{2,4mm}$\bullet$}
\rput[l]{-45}(1.2,3.08){\hspace{1,4mm}$\bullet$}
}
\rput[l](0.3,1.3){$s$}
\rput[l](1.8,2.5){$t$}
\rput[bl](5.5;0){
\pcline(-1.35,0.53)(1.2,0.53)\ncput*{$\lambda_1$}
\pcline(1.2,0.53)(3.75,3.08)\ncput*{$\lambda_3$}
\pcline(3.75,3.08)(1.2,3.08)\ncput*{$\lambda_1$}
\pcline(1.2,3.08)(-1.35,0.53)\ncput*{$\lambda_3$}
\pcline(-1.35,0.53)(3.75,3.08)\ncput*{$\lambda_5$}
\rput[l]{13}(-1.35,0.53){\hspace{4,3mm}$\bullet$}
\rput[l]{-66}(1.2,3.08){\hspace{1,0mm}$\bullet$}
}
\rput[l](6.4,1.1){$s$}
\rput[l](6.7,2.4){$t$}
\rput[l](3.7,2){\pcline{->}(0,0)(1;0)\Aput{${\bf T}_{st}^{-1}$}}
\end{pspicture}
\caption{Flip at an edge labeled by $\lambda_2$, for a once-punctured torus}
\label{fig:OPT2}
\end{figure}
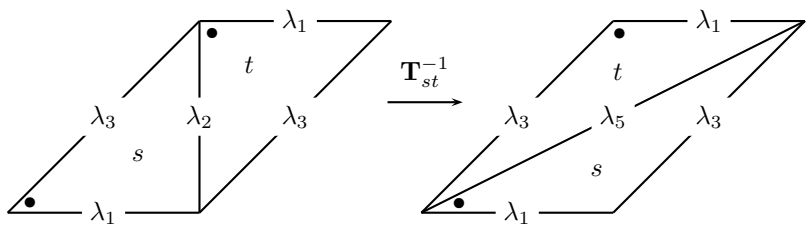

After a flip at the edge labeled by $\lambda_2$ or $\lambda_3$, are we able to go on flipping at other edges, to obtain new LD-triangulations? For this, the above obtained conditions $x_1=1,x_2=-1,y_1\neq \pm y_2, 2y_1\neq y_2$ must hold for the `new' non-singular weights $\lambda_1', \lambda_2'$ playing the role of the previous $\lambda_1,\lambda_2$ (more precisely, $\lambda_1'=\lambda_1$, $\lambda_2'=\lambda_4$, or $\lambda_1'=\lambda_1$,$\lambda_2'=\lambda_3$). One can come up with such an example; if in the initial picture we put $\lambda_1=(1,\alpha)$, $\lambda_2= (-1,\beta)$ for any non-zero real numbers such that $\alpha/\beta$ is irrational, then one can easily show that these conditions are preserved all the time, so in particular, this connected component of sane LD-groupoid has infinitely many objects; lots of elements of the mapping class group $\mathrm{SL}(2,\mathbb{Z})$, perhaps the elements of a subgroup isomorphic to $\mathbb{Z}$, can be realized as morphisms (or, paths) inside this connected component. One could investigate similar questions for other higher genus surfaces, to see which subgroups (and how big they are) of the mapping class group can be represented in a connected component of the sane LD-groupoid.

\vs

A referee pointed out that the case of genus three with one puncture is dealt with in \cite{Kash00b}, or in its reference \cite{K98b}, but an explicit sane LD-triangulation of such a surface is not found there; so one might try to construct a sane LD-triangulation for this case, using these works \cite{Kash00b} \cite{K98b}.

\subsection{On representation theory of other Hopf algebras}

For any Hopf algebra, one might study the left and right duals and the left and right Hom representations as done in the present paper, to see if we get anything new, i.e. a new representation of the Kashaev-type groupoid. For the case when the tensor product of two irreducibles decompose into the direct sum of a single irreducible as in the present paper, one might hope to get something like we obtained in this paper, i.e. the ${\bf A}, {\bf T}, {\bf P}$ operators, giving a representation of the LD-groupoid. A natural candidate is a Borel subalgebra of the higher rank quantum groups $\mathcal{U}_q(\mathfrak{sl}_n)$, for $q$ a root of unity. Or, maybe one can try (Borel subalgebras of) affine quantum groups too. See also \cite{GKT}.

\subsection{On cyclic quantum dilogarithm}
\label{subsec:discussion_on_Phi_C}

In the present paper, we avoided going too deeply into properties of the cyclic quantum dilogarithm, which we would need e.g. if we want to express the action of the operator $F$ \eqref{eq:F} on the basis vectors explicitly and neatly, and to compute the formula for ${\bf A}$ more directly using its definition. Properties that I suggest to establish is a `Fourier transform formula' and the `reflection formula', as an analog of the counterpart \cite{Volkov} \cite{Ruij} for the compact or non-compact quantum dilogarithm functions, defined for $q$ not a root of unity. Vaguely saying, the former says that the Fourier transform of a quantum dilogarithm is also a quantum dilogarithm. 

\vs

Let us be more explicit. What can be called a {\em Fourier transform} on $\mathbb{C}^N$ is the operator
$$
\mathcal{F} : \mathbb{C}^N \to \mathbb{C}^N, \quad e_i \mapsto \sum_{j=0}^{N-1} q^{-2ij} e_j, \quad \forall i =0,\ldots,N-1,
$$
which has the properties: $\mathcal{F}^{-1}(e_i) = \sum_j q^{2ij} e_j$, $\mathcal{F}^2(e_i) = e_{-i}$, $\mathcal{F}^4 = \mathrm{id}$, and
$$
\mathcal{F} \, A \, \mathcal{F}^{-1} = B,\qquad
\mathcal{F} \, B \, \mathcal{F}^{-1} = A^{-1}.
$$
A generalization of such operator might look like
$$
\mathcal{G} : \mathbb{C}^N \to \mathbb{C}^N, \quad e_i \mapsto \sum_{j=0}^{N-1} q^{ a i^2 + bij + c j^2 + d i + f j} e_j, \quad \forall i=0,\ldots,N-1,
$$
with some numbers $a,b,c,d,f$; such operator was referred to in \S\ref{subsec:computation_of_A} as an analog of Fourier transform, or exponential of a quadratic expression in Schr\"odinger operators. As seen in Prop.\ref{prop:computation_of_A}, our ${\bf A}$ operator is an example. These operators appear as `metaplectic representations', in the sense that
$$
\mathcal{G} \, A \, \mathcal{G}^{-1} = q^{2m_1} A^{g_1} B^{g_2}, \qquad
\mathcal{G} \, B \, \mathcal{G}^{-1} = q^{2m_2} A^{g_3} B^{g_4},
$$
for some integers $m_1,m_2,g_1,g_2,g_3,g_4$. So, conjugation by such operator on a cyclic quantum dilogarithm might help finding the explicit action of a cyclic quantum dilogarithm on basis vectors. One suggestion to be studied for a sought-for `Fourier transform formula for cyclic quantum dilogarithm' is to investigate the operator $\mathcal{F} \, \Phi^q_{{\bf a}, {\bf c}}(A)$, and see if it can be written almost as $\Phi^q_{{\bf a}, {\bf c}}(C)$ for some operator $C$. I think a basic question is to investigate the expression $\sum_{j=0}^{N-1} q^{-2ij} w({\bf a}, {\bf c}|i)$ and see if it equals something like $q^{aj^2 + bj} \, w({\bf a}', {\bf c'}|j)$. 

\vs

I also expect that there is a `reflection identity' that might look something like $\Phi^q_{{\bf a}, {\bf c}}(C) \, \Phi^q_{{\bf a}, {\bf c}}(C^{-1}) = $ the exponential of a quadratic expression in Schr\"odinger operators, or with the parameters ${\bf a},{\bf c}$ of the second factor modified appropriately; more precise form could be deduced from the identity in Prop.\ref{prop:TAT_AAP}. This problem is easier to tackle than the Fourier transform formula, for a version of a reflection identity is established already in \cite{Kash00b} (without proof).



\end{document}